\documentclass[11pt,reqno]{amsart}
\usepackage[T1]{fontenc}
\usepackage{lmodern}
\usepackage[a4paper,margin=27mm]{geometry}
\usepackage{amsmath,amssymb,mathtools}
\usepackage{graphicx}
\usepackage{tikz}
\usepackage{standalone}
\usetikzlibrary{arrows.meta}
\usepackage{microtype}
\usepackage{enumitem}
\usepackage{needspace}
\usepackage{xcolor}
\usepackage{hyperref}
\hypersetup{colorlinks=true,linkcolor=blue!50!black,citecolor=blue!50!black,
  urlcolor=blue!50!black}
\setlist[enumerate]{topsep=5pt plus 2pt minus 1pt,
  itemsep=3pt plus 1pt minus 1pt,parsep=0pt}
\AtBeginDocument{\raggedbottom}
\makeatletter
\renewcommand{\subsection}{\@startsection{subsection}{2}%
  \z@{.8\linespacing\@plus.3\linespacing}{.35\linespacing}%
  {\normalfont\bfseries}}
\makeatother
\renewcommand{\MR}[1]{}
\allowdisplaybreaks[2]
\numberwithin{equation}{section}
\newtheorem{theorem}{Theorem}[section]
\newtheorem{citedtheorem}{Theorem}
\newtheorem{proposition}[theorem]{Proposition}
\newtheorem{lemma}[theorem]{Lemma}
\newtheorem{corollary}[theorem]{Corollary}
\theoremstyle{remark}
\newtheorem{remark}[theorem]{Remark}
\newtheorem{example}[theorem]{Example}
\newtheorem{question}{Question}
\newcommand{\eps}{\varepsilon}
\newcommand{\cE}{\mathcal E}
\newcommand{\D}{D}
\newcommand{\Var}{\operatorname{Var}}
\newcommand{\vol}{\operatorname{Vol}}
\newcommand{\1}{\mathbf 1}
\newcommand{\dd}{\,\mathrm d}
\newcommand{\R}{\mathbb R}
\newcommand{\N}{\mathbb N}
\newcommand{\norm}[1]{\left\lVert#1\right\rVert}

\newcommand{\ip}[2]{\left\langle#1,#2\right\rangle_m}
\newcommand{\CD}{\mathrm{CD}}
\DeclareMathOperator{\dist}{dist}
\DeclareMathOperator{\Pois}{Pois}
\newcommand{\E}{\mathbb{E}}

\title[Wang's log-Harnack inequality and its applications]{Wang's log-Harnack inequality\\ for Markov chains and its applications}

\author{Yanlong Ding}
\address{\textnormal{Yanlong Ding:}\endgraf\nobreak\indent School of Mathematical Sciences, University of Science and Technology of China, Hefei 230026, China.}
\email{dylustc@mail.ustc.edu.cn}

\author{Shiping Liu}
\address{\textnormal{Shiping Liu:}\endgraf\nobreak\indent School of Mathematical Sciences, University of Science and Technology of China, Hefei 230026, China.}
\email{spliu@ustc.edu.cn}

\author{Chiyu Zhou}
\address{\textnormal{Chiyu Zhou:}\endgraf\nobreak\indent School of Mathematical Sciences, University of Science and Technology of China, Hefei 230026, China.}
\email{dovong@mail.ustc.edu.cn}

\date{Revised 30 October 2026}

\begin{document}
\begin{abstract}
We consider irreducible reversible Markov chains whose nonzero
 off-diagonal transition probabilities are bounded below by $\eps>0$
and whose reversible measure satisfies $m_-\le m(x)\le\kappa m_-$.
Under the Bakry--\'Emery curvature condition $\CD(K,\infty)$, with
$K\in\R$, we establish Wang's log-Harnack inequality
\[
 P_t(\log f)(x)\le\log P_tf(y)
       +\frac{C\dist(x,y)^2}{\int_0^t e^{2Ks}\dd s},
\]
for bounded positive $f$ bounded away from zero and all $x,y\in V$, where $C$ is a constant dependent only on $\varepsilon$. Building upon this inequality, we prove volume doubling property, which extends the work of Guo-Huang-Huang to general weighted setting. We further prove a sharp volume growth and Li-Yau gradient estimates under $\CD(0,\infty)$ .
\end{abstract}
\maketitle

\section{Introduction}

\subsection{Background and previous results}
On a complete \(n\)-dimensional Riemannian manifold \(M\) with nonnegative Ricci curvature, the Bishop--Gromov comparison theorem (see, e.g., \cite{Petersen2016}) states that, for every \(o\in M\) and \(0<r\le R\),
\[
\frac{\operatorname{vol}(B(o,R))}{\operatorname{vol}(B(o,r))}
\le \left(\frac{R}{r}\right)^n.
\]
Consequently, \(M\) has polynomial volume growth of order at most \(n\) and satisfies the volume doubling property with doubling constant \(2^n\). Together with local Poincar\'e inequalities, these volume
estimates provide the foundation for heat-kernel bounds, Harnack inequalities, and
Liouville theorems \cite{LiYau1986,Grigoryan1992,SaloffCoste1992}.

In this paper, we study volume estimates for graphs and Markov chains under Bakry--\'Emery curvature bounds. The proof of the Bishop--Gromov comparison theorem does not readily carry over to this discrete setting. Even establishing polynomial volume growth or uniform volume doubling from such local curvature bounds requires a different approach; see~\cite{Munch2019CD,RussPajot,BlacharPajotSalez2026,GuoHuangHuang}.
%For graphs or Markov chains, there is no good analogue of  Bishop--Gromov theorem under assuming Bakry--\'Emery curvature bounds.

%\subsection{Progress on polynomial growth and volume doubling}

The Bakry--\'Emery curvature-dimension condition \(\mathrm{CD}(K,\infty)\) and the associated \(\Gamma\)-calculus provide an alternative formulation of Ricci curvature lower bounds in terms of the Laplace operator~\cite{BakryEmery1985}. Specifically, for the Laplace--Beltrami operator on a Riemannian manifold \((M,g)\), the condition \(\mathrm{CD}(K,\infty)\) is equivalent to \(\operatorname{Ric}_g\ge Kg\).
%For the Laplace--Beltrami operator on a Riemannian manifold $(M,g)$, the condition $\CD(K,\infty)$ is equivalent to $\operatorname{Ric}_g\ge Kg$.
In the discrete setting, Bakry--\'Emery theory was initiated in~\cite{elworthy,michael,LinYau2010} and has since become an effective tool for studying geometric, analytic, and combinatorial properties of graphs and Markov chains; see, for instance,~\cite{FathiShu2018,HornJGT,HM2024,KMY2021,MuenchRose2020,SalezJEMS}.

\begin{figure}[htbp]
 \centering
 \resizebox{0.90\textwidth}{!}{%
   % TikZ figure fragment, included by main.tex.
% Dashed segments indicate continuation of the infinite graphs.
\begingroup

\tikzset{
  model edge/.style={draw=black,line width=0.55pt,line cap=round},
  continuation/.style={model edge,dash pattern=on 4pt off 3pt},
  bridge edge/.style={model edge,
    preaction={draw=white,line width=2.1pt,line cap=round}},
  bridge continuation/.style={continuation,
    preaction={draw=white,line width=2.1pt,line cap=round}}
}
\newcommand{\VertexRadius}{2pt}
\newcommand{\LatticeExtent}{1.65}

\begin{tikzpicture}

% Q_3: eight vertices and twelve edges, without infinite continuations.
\begin{scope}[shift={(0cm,0cm)},
  x={(1.00cm,0cm)},y={(0.48cm,0.42cm)},z={(0cm,1.08cm)}]
  \foreach \k in {-0.85,0.85} {
    \foreach \j in {-0.85,0.85} {
      \draw[model edge] (-0.85,\j,\k) -- (0.85,\j,\k);
    }
    \foreach \i in {-0.85,0.85} {
      \draw[model edge] (\i,-0.85,\k) -- (\i,0.85,\k);
    }
  }
  \foreach \i in {-0.85,0.85} {
    \foreach \j in {-0.85,0.85} {
      \draw[bridge edge] (\i,\j,-0.85) -- (\i,\j,0.85);
    }
  }
  \foreach \i in {-0.85,0.85} {
    \foreach \j in {-0.85,0.85} {
      \foreach \k in {-0.85,0.85} {
        \fill (\i,\j,\k) circle[radius=\VertexRadius];
      }
    }
  }
\end{scope}

% Z^3: the block {-1,0,1}^3 and continuations in all six directions.
% The projection has no coincident vertices or vertices inside other edges.
\begin{scope}[shift={(4.8cm,0cm)},
  x={(0.82cm,0cm)},y={(0.32cm,0.29cm)},z={(0cm,0.84cm)}]
  \foreach \k in {-1,0,1} {
    \foreach \j in {-1,0,1} {
      \draw[model edge] (-1,\j,\k) -- (1,\j,\k);
      \draw[continuation] (-1,\j,\k) -- (-\LatticeExtent,\j,\k);
      \draw[continuation] (1,\j,\k) -- (\LatticeExtent,\j,\k);
    }
    \foreach \i in {-1,0,1} {
      \draw[model edge] (\i,-1,\k) -- (\i,1,\k);
      \draw[continuation] (\i,-1,\k) -- (\i,-\LatticeExtent,\k);
      \draw[continuation] (\i,1,\k) -- (\i,\LatticeExtent,\k);
    }
  }
  \foreach \i in {-1,0,1} {
    \foreach \j in {-1,0,1} {
      \draw[bridge edge] (\i,\j,-1) -- (\i,\j,1);
      \draw[bridge continuation] (\i,\j,-1) -- (\i,\j,-\LatticeExtent);
      \draw[bridge continuation] (\i,\j,1) -- (\i,\j,\LatticeExtent);
    }
  }
  \foreach \i in {-1,0,1} {
    \foreach \j in {-1,0,1} {
      \foreach \k in {-1,0,1} {
        \fill (\i,\j,\k) circle[radius=\VertexRadius];
      }
    }
  }
\end{scope}

% T_3: the radius-three ball, with two outgoing dashed edges at each leaf.
% The root has three children; each other internal vertex has two children.
\begin{scope}[shift={(9.6cm,-0.17cm)}]
  \coordinate (root) at (0,0);
  \foreach \a [count=\arm] in {90,210,330} {
    \coordinate (a\arm) at (\a:0.75);
    \draw[model edge] (root) -- (a\arm);
    \foreach \b [count=\child] in {-32,32} {
      \pgfmathsetmacro{\angleB}{\a+\b}
      \coordinate (b\arm\child) at
        ([shift={(\angleB:0.58)}]a\arm);
      \draw[model edge] (a\arm) -- (b\arm\child);
      \foreach \c [count=\leaf] in {-16,16} {
        \pgfmathsetmacro{\angleC}{\angleB+\c}
        \coordinate (c\arm\child\leaf) at
          ([shift={(\angleC:0.40)}]b\arm\child);
        \draw[model edge] (b\arm\child) -- (c\arm\child\leaf);
        \foreach \out in {-14,14} {
          \pgfmathsetmacro{\angleOut}{\angleC+\out}
          \draw[continuation] (c\arm\child\leaf)
            -- ++(\angleOut:0.31);
        }
        \fill (c\arm\child\leaf) circle[radius=\VertexRadius];
      }
      \fill (b\arm\child) circle[radius=\VertexRadius];
    }
    \fill (a\arm) circle[radius=\VertexRadius];
  }
  \fill (root) circle[radius=\VertexRadius];
\end{scope}

\node[font=\normalsize,anchor=north] at (0,-2.18) {$Q_3$};
\node[font=\normalsize,anchor=north] at (4.8,-2.18) {$\mathbb Z^3$};
\node[font=\normalsize,anchor=north] at (9.6,-2.18) {$T_3$};

\end{tikzpicture}
\endgroup%
 }
 \caption{The hypercube $Q_3$, the cubic lattice $\mathbb Z^3$, and
 the $3$-regular tree $T_3$.}
 \label{fig:curvature-models}
\end{figure}
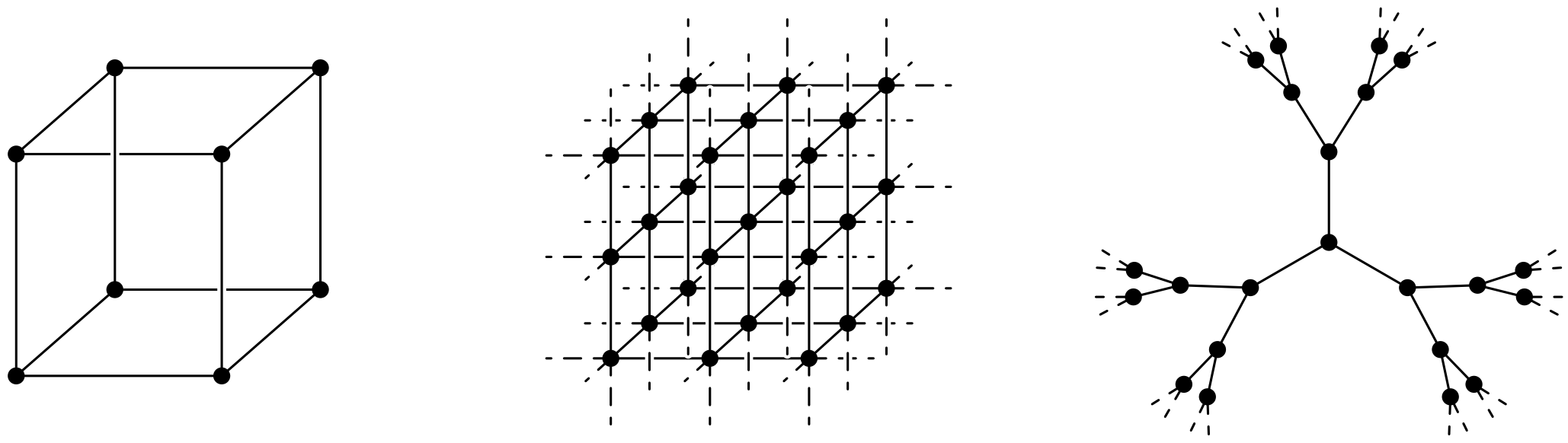

Three basic examples are the \(n\)-dimensional hypercube \(Q_n\), the integer lattice \(\mathbb{Z}^n\), and the \(d\)-regular tree \(T_d\) with \(d\ge3\). With respect to the unnormalized graph Laplacian, these graphs satisfy \(\mathrm{CD}(K,\infty)\) with \(K=2\), \(K=0\), and \(K=2-d\), respectively~\cite{CushingLiuPeyerimhoff2020}; see Figure~\ref{fig:curvature-models}. The hypercube \(Q_n\) has \(2^n\) vertices, whereas the lattice \(\mathbb{Z}^n\) and the regular tree \(T_d\) exhibit polynomial and exponential volume growth, respectively. More precisely, for \(R>1\),
\[
|B_{\mathbb{Z}^n}(o,R)|\asymp R^n,
\qquad
|B_{T_d}(o,R)|\asymp(d-1)^R.
\]
Notably, the polynomial volume growth exponent \(n\) of \(\mathbb{Z}^n\) is exactly half its vertex degree \(2n\).
\iffalse
In the discrete setting, Bakry--\'Emery theory was initiated in \cite{elworthy,michael,LinYau2010}
and has since become an effective tool for studying geometric, analytic, and combinatorial aspects of graphs and Markov chains, see,
for instance, \cite{FathiShu2018,HornJGT,HM2024, KMY2021, MuenchRose2020,SalezJEMS}.  The following three combinatorial graphs, hypercube $Q_n$, the lattice
$\mathbb Z^n$, and the $d$-regular tree $T_d$ ($d\ge3$), satisfy $\CD(K,\infty)$,  with $K=2,\ 0$ and $2-d$ respectively~\cite{CushingLiuPeyerimhoff2020}; see Figure~\ref{fig:curvature-models}.
The hypercube has $2^n$ vertices, while the lattice and the regular
tree exhibit polynomial and exponential volume growth, respectively: For $R>1$,
\[
 |B_{\mathbb Z^n}(o,R)|\asymp R^n,\qquad
 |B_{T_d}(o,R)|\asymp(d-1)^R.
\]
Notably, the polynomial volume growth exponent $n$ of $\mathbb Z^n$
is exactly half its vertex degree $2n$.
\fi

One approach to volume estimates on graphs is based on Li--Yau gradient estimates and the resulting heat kernel bounds. However, the failure of the usual chain rule prevents a direct adaptation of the standard Li--Yau argument and has motivated the introduction of several nonlinear curvature-dimension conditions. Bauer, Horn, Lin, Lippner, Mangoubi, and Yau~\cite{BauerEtAl2015} introduced the exponential curvature dimension condition \(\mathrm{CDE}(K,N)\) and obtained Li--Yau estimates and polynomial volume growth on bounded-degree graphs under \(\mathrm{CDE}(0,N)\) with finite $N$. Horn, Lin, Liu, and Yau~\cite{HornLinLiuYau2019} subsequently established Li--Yau estimates and volume doubling under the stronger condition \(\mathrm{CDE}'(0,N)\).
%We mention that \(\mathrm{CDE}'(N,K)\) also implies the Bakry--\'Emery condition (\mathrm{CD}(K,N))~\cite{}.

M\"unch~\cite{Munch2018Nonlinear} developed the \(\mathrm{CD}\psi\) framework and obtained various Li--Yau-type estimates on finite graphs under \(\mathrm{CD}\psi(0,N)\). For further Li--Yau-type estimates on finite graphs, see ~\cite{DierKassmannZacher2021,WeberZacherJFA,WeberZacher2023HeatKernel}.
%Dier, Kassmann, and Zacher~\cite{DierKassmannZacher2021} established further Li--Yau-type estimates on finite graphs under \(\mathrm{CD}(F;0)\). 
More recently, Zhang~\cite{Zhang2026CDpsiDoubling} obtained a family of Li--Yau estimates on locally finite weighted graphs under \(\mathrm{CD}\psi(-K,N)\), \(K\ge0\), and derived volume doubling estimates under additional assumptions on \(\psi\).

We note that these nonlinear curvature-dimension conditions are stronger
than the Bakry--\'Emery condition \(\mathrm{CD}(K,\infty)\).
Indeed, M\"unch~\cite[Theorem~3.2 and Corollary~3.3]{MunchRemarks}
proved that both \(\mathrm{CDE}'(K,N)\) and \(\mathrm{CD}\psi(K,N)\)
imply \(\mathrm{CD}(K,\infty)\), under suitable assumptions on \(\psi\)
in the latter case.
The \(\mathrm{CDE}(K,N)\) condition also implies \(\mathrm{CD}(K,\infty)\), as can be seen by adding a sufficiently large constant
to the test function and passing to the limit.\footnote{By locality, it suffices
to consider bounded test functions \(g\). Fix a vertex \(x\).
If \(\Delta g(x)<0\), apply \(\mathrm{CDE}(K,N)\) at \(x\) to
\(C+g\) and let \(C\to\infty\).
The term
\(\Gamma\bigl(g,\Gamma(g)/(C+g)\bigr)(x)\) tends to zero,
whereas \(\Gamma_2(g)\), \(\Gamma(g)\), and \(\Delta g\)
are unchanged by the addition of a constant.
The case \(\Delta g(x)>0\) follows by replacing \(g\) with \(-g\),
and the case \(\Delta g(x)=0\) follows by continuity.}

For the Bakry-\'Emery curvature condition, Cushing, Liu, and
Peyerimhoff \cite[Conjecture~9.10]{CushingLiuPeyerimhoff2020} asked
whether bounded-degree graphs have polynomial volume growth, with
constants depending only on the degree bound. Under
$\CD(0,N)$ with finite $N$, M\"unch \cite{Munch2019CD} proved doubling
on finite graphs using a modified nonlinear heat equation, and Russ
and Pajot \cite{RussPajot} extended this method to infinite weighted
graphs of bounded geometry. Blachar, Pajot, and Salez
\cite{BlacharPajotSalez2026} obtained doubling under $\CD(0,\infty)$
for edge-regular graphs through finite-dimensional self-improvement.

Recently, Guo, Huang, and Huang~\cite{GuoHuangHuang} made a breakthrough
in this direction by proving volume doubling and polynomial volume growth
for unweighted graphs of bounded degree satisfying $\CD(0,\infty)$.
\Needspace{12\baselineskip}
Their results can be stated as follows.
\begin{citedtheorem}[{\cite[Theorem~1.1, Corollary~1.2]{GuoHuangHuang}}]
\label{thm:ghh-volume-doubling}
Let $G$ be a connected simple graph of maximum degree $d<\infty$
whose non-normalized Laplacian satisfies $\CD(0,\infty)$.
Then there is a constant $K(d)<\infty$ such that, for every $o\in V$ and every integer $R\ge1$, 
\[
 |B(o,2R)|\le K(d)\,|B(o,R)|,
\]
and, consequently, 
\[
 |B(o,R)|\le (1+R)^{\lceil2\log_2 K(d)\rceil}.
\]
\end{citedtheorem}

In this paper, we extend the results of Guo, Huang, and
Huang~\cite{GuoHuangHuang} to a general reversible Markov chain setting.
We establish volume doubling under arbitrary
lower curvature bounds and sharp polynomial volume growth bounds
under nonnegative curvature. Using these volume estimates, we prove a Li-Yau gradient estimate under the Bakry--\'Emery codition $CD(0,\infty)$.

\subsection{Main results}
Let $P=(p(x,y))_{x,y\in V}$ be an irreducible Markov kernel on a
countable state space $V$, reversible with respect to a positive
measure $m$. The associated Markov generator is $L=P-I$. The support graph of $P$ has vertex set $V$ and edge set
$E=\{\{x,y\}:x\ne y,\ p(x,y)>0\}$. We write $x\sim y$ when
$\{x,y\}\in E$. For any $A\subseteq V$, define its volume by
$\vol(A)=\sum_{x\in A}m(x)$.

Throughout the paper, assume that, for fixed $0<\eps\le1$, $m_->0$,
and $\kappa\ge1$,
\begin{equation}\label{eq:geometry}
 p(x,y)\ge\eps\quad(x\sim y),\qquad
 m_-\le m(x)\le\kappa m_-\quad(x\in V).
\end{equation}

Define the Laplacian associated with $P$ by
\[
 Lf(x)=\sum_y p(x,y)\bigl(f(y)-f(x)\bigr).
\]
Let $P_t=e^{tL}$ be its heat semigroup and $p_t(x,y)$ be its heat kernel. Next, we give our main results:

\begin{theorem}[Volume doubling]
\label{thm:main}
Let $P$ be a Markov kenel satisfying \eqref{eq:geometry} and  $\CD(K,\infty)$ for
some $K\in\R$. Put $K_-:=\max\{-K,0\}$. There are constants
$C_0\ge1$ and $0\le C_1<\infty$, depending only on $\eps,\kappa$, such that
\begin{equation}\label{eq:doubling-K}
 \vol(B(o,2R))\le C_0e^{C_1\sqrt{K_-}\,R}\vol(B(o,R)).
\end{equation}
for any $o\in V$ and $R>0$.
\end{theorem}

A key step in the proof of Theorem~\ref{thm:main}
is to show the following Wang's log-Harnack inequality.
For any $K\in \mathbb{R}$, define
\begin{equation}\label{eq:curvature-time}
 I_K(t)=\int_0^t e^{2Ks}\dd s
 =\begin{cases}
 (e^{2Kt}-1)/(2K),&K\ne0,\\
 t,&K=0,
 \end{cases}
\end{equation}
\begin{theorem}[Wang's log-Harnack inequality]\label{thm:logharnack}
Let $P$ be a Markov kernel satisfying \eqref{eq:geometry} and  $\CD(K,\infty)$ for
some $K\in\R$. There exists a constant $C_{\mathrm{LH}}<\infty$ dependent only on $\varepsilon$ such that, for every
$t>0$ and $x,y\in V$,
\begin{equation}\label{eq:logharnack}
 \sum_{z\in V}p_t(x,z)\log \frac{p_t(x,z)}{p_t(y,z)}
 \le\frac{C_{\mathrm{LH}}}{I_K(t)}\dist(x,y)^2.
\end{equation}
Equivalently, every bounded $f$ with $\inf f>0$
satisfies
\[
 P_t(\log f)(x)\le\log P_tf(y)
                +\frac{C_{\mathrm{LH}}}{I_K(t)}\dist(x,y)^2.
\]
\end{theorem}

The study of the above type of log-Harnack inequality was initiated by Feng-Yu Wang~\cite{Wang2010LogHarnack}. He used it to
characterize lower bounds of Bakry--\'Emery curvature on Riemannian manifolds. The log-Harnack
inequality for stochastic differential equations in Hilbert
spaces were established by R\"ockner and Wang~\cite{RocknerWang2010LogHarnack}.
Its connections with heat flow and optimal transportation were further
explored by Bakry, Gentil, and
Ledoux~\cite{BakryGentilLedoux2015Harnack}.

Li-Yau gradient estimate \cite{LiYau1986} is a fundamental tool in geometric analysis \cite{LiXu2011DifferentialHarnack,BakryBolleyGentil2017LiYau}; for a more comprehensive introduction, see \cite{SchoenYau1994Lectures}.  As discussed above, Li-Yau type gradient estimates for the heat equation on graphs have been studied in several works under various nonlinear curvature assumptions~\cite{BauerEtAl2015,HornLinLiuYau2019,Munch2018Nonlinear,DierKassmannZacher2021,WeberZacherJFA,WeberZacher2023HeatKernel,Zhang2026CDpsiDoubling}.
%For the standard linear heat equation, Bauer et al.~\cite{BauerEtAl2015} obtained such estimates under $\mathrm{CDE}(n,0)$, and Horn, Lin, Liu, and Yau~\cite{HornLinLiuYau2019} under $\mathrm{CDE}'(n,0)$. On finite graphs, M\"unch~\cite{Munch2018Nonlinear} worked under $\mathrm{CD}^{\log}(n,0)$, and Dier, Kassmann, and Zacher~\cite{DierKassmannZacher2021} obtained Li--Yau estimates under $\CD(F;0)$. Weber and Zacher~\cite{WeberZacher2023HeatKernel} subsequently obtained a sharp Li--Yau estimate on unweighted complete graphs. More recently, Zhang~\cite{Zhang2026CDpsiDoubling} established a family of Li--Yau type estimates under $\mathrm{CD}\psi(n,-K)$.

For a modified nonlinear heat flow, M\"unch~\cite{Munch2019CD} obtained
Li--Yau estimates under $\CD(0,N)$ on finite graphs, and Russ and
Pajot~\cite{RussPajot} extended these estimates to infinite weighted graphs
of bounded geometry under the same curvature assumption.
Guo, Huang, and Huang~\cite{GuoHuangHuang} obtained analogous estimates
under $\CD(0,\infty)$ on unweighted graphs.

Despite these advances, establishing Li--Yau estimates for the
standard linear heat equation on graphs satisfying
\(\CD(0,\infty)\) has remained an open problem.
Building on the volume doubling property established above,
we prove the following Li--Yau gradient estimate.

\begin{theorem}[Li--Yau gradient estimate]
\label{thm:li-yau} Let $P$ be a Markov kernel satisfying \eqref{eq:geometry} and  $\CD(0,\infty)$.
There are constants
$\alpha_0>1$ and $C_{LY}<\infty$, depending only on $\eps,\kappa$,
such that if $f$ is a bounded, nonnegative, nonzero function, then
\[
 \begin{gathered}
 \Gamma(\log P_tf)(x)-\alpha\partial_t\log P_tf(x)
 \le\frac{C_{LY}\alpha}{t},
 \end{gathered}
\]
for $\alpha\geq \alpha_0,\ t>0$ and $x\in V$.
\end{theorem}

We also establish the following sharp polynomial volume growth bound
in terms of the maximum vertex degree.
\begin{theorem}[Sharp volume growth]
\label{thm:sharp-volume}
Let $P$ be a Markov kernel satisfying \eqref{eq:geometry} and  $\CD(0,\infty)$.
Set $d=\sup_{x\in V}\deg(x)$. There is a constant
$C=C(\eps,\kappa)<\infty$ such that, for any $o\in V$ and $R\geq 1$, 
\begin{equation}\label{eq:sharp-volume-main}
 \vol(B(o,R))\le m_-CR^{\left\lfloor d/2\right\rfloor}.
\end{equation}
 The exponent
$\lfloor d/2\rfloor$ is optimal. 
\end{theorem}

Recall that, for the lattice $\mathbb Z^n$, the volume growth $|B(o,R)|\asymp R^n$, with $n$ being equal to half of its degree.

\subsection{Proof strategy}
An approach to proving volume doubling or polynomial volume
growth on graphs is to establish Li--Yau gradient estimates; see \cite{BauerEtAl2015, HornLinLiuYau2019, GuoHuangHuang, Zhang2026CDpsiDoubling}. Under the
condition $\CD(K,\infty)$, however, obtaining such
estimates is particularly challenging. In this article,  we develop an alternative approach to volume doubling based on Wang's log-Harnack inequality.

To prove the volume doubling, we first establish and then apply Wang's log-Harnack inequality (Theorem~\ref{thm:logharnack})  to derive useful entropy estimates. These entropy estimates led to a preliminary polynomial volume growth bound, which provides the integrability needed for the analytic arguments that follow. We then combine gradient estimates for the  modified heat flow \(\partial_su=Lu+\Gamma u\), introduced by M\"unch \cite{Munch2019CD}, with weighted time integrals of an auxiliary flow to control the spatial spread of the heat flow, see \eqref{eq:curved-Qbounds} below. Subsequently, we employ  Wang's log-Harnack inequality again to convert this heat flow estimate into volume doubling.

For the Li--Yau gradient estimates under $\CD(0,\infty)$,  we first apply volume doubling to get a local Poincar\'e inequality.  Then by results of Delmotte \cite{Delmotte1998Harmonic},  we derive the parabolic Harnack inequalities (Theorem~\ref{thm:parabolic-harnack}). The Harnack inequality could bound the negative time derivative of the logarithm of the solution of heat equation (Lemma~\ref{lem:ly-time}). Combining this bound with local gradient
estimates for nonnegative solutions of heat equation
(Lemma~\ref{lem:sharp-positive-regularity}) yields a square-root
gradient estimate (Lemma~\ref{lem:ly-sqrt}). 
We then bound the ratios of values of a heat solution at adjacent
vertices, allowing us to pass from square-root to logarithmic
gradient estimates.

To establish the sharp volume growth, we apply the parabolic Harnack inequalities to yield the Liouville property for positive harmonic functions, see Corollary~\ref{cor:positive-liouville} below. We show that the linear space of harmonic functions with constant gradient norm has dimension at most \(\lfloor d/2\rfloor\), where \(d\) is the maximum combinatorial degree. This dimension bound turns out to be crucial in establishing a sharp logarithmic entropy bound. Wang's log-Harnack inequality then translates this entropy bound into the sharp volume growth estimate \eqref{eq:sharp-volume-main}.

\subsection*{Use of artificial intelligence}
The authors used generative artificial intelligence tools to help
explore and refine the entropy methods, heat semigroup techniques, and
harmonic function arguments presented in this article.
The authors are responsible for independently verifying all mathematical
arguments and computations and for the final content of the article.

\section{Preliminaries}
\subsection{Graph and the heat semigroup}
Let $E,F$ be countable sets. A \emph{Markov kernel} from $E$
to $F$ is a function $M:E\times F\to[0,1]$ such that
\[
 \sum_{b\in F}M(a,b)=1\qquad(a\in E).
\]

Let $P=(p(x,y))_{x,y\in V}$ be a Markov kernel on an
at most countable set $V$, reversible with respect to a positive
measure $m$. Thus
\[
 m(x)p(x,y)=m(y)p(y,x).
\]
$p(x,x)>0$ is allowed.

Its \emph{support graph} is the connected simple graph $G=(V,E)$ with
\[
 E=\{\{x,y\}:x\ne y,\ p(x,y)>0\}.
\]
We write $x\sim y$ when $\{x,y\}\in E$.  Set $ \deg(x)=|\{y:y\sim x\}|$. Define the Laplacian associated with $P$ by
\[
 Lf(x)=\sum_y p(x,y)\bigl(f(y)-f(x)\bigr).
\]
Under \eqref{eq:geometry}, $\deg(x)\le\eps^{-1}$, and all metric balls
are finite. 

Write $\dist$ for the combinatorial distance, equivalently
\[
 \dist(x,y)=\min\{k\ge0:P^k(x,y)>0\}.
\]
We call $P$ irreducible if $\dist(x,y)<\infty$ for any $x,y\in V$.

For $A\subseteq V$ and $R\ge0$, set
\begin{equation*}
 \vol(A)=\sum_{z\in A}m(z),\qquad
 B(o,R)=\{x\in V:\dist(o,x)\le R\}.
\end{equation*}

Let $ \ell^p(V,m)=\left\{f:V\to \mathbb{R}:\ \sum_{x\in V}m(x)|f(x)|^p<\infty\right\}$ and $\ell^\infty(V)$ be the set of bounded function on $V$.

For $f,g\in \ell^2(V,m)$, set
\begin{equation*}
 \begin{aligned}
 \cE(f,g)=\sum_{\{x,y\}}m(x)p(x,y)(f(x)-f(y))(g(x)-g(y))=-\ip{f}{Lg}.
 \end{aligned}
\end{equation*}
Also denote $\cE(f)=\cE(f,f)$. For other symmetric edge weights $a$, the
corresponding energy is denoted by $\cE_a$.

\Needspace{8\baselineskip}
Write $P_t=e^{tL}$ for its minimal heat semigroup on $\ell^2(V,m)$,
with its natural extensions to $\ell^1(V,m)$ and $\ell^\infty(V)$.
Write $p_t(x,y)=(P_t\1_{\{y\}})(x)$ for its transition probabilities.
For every $t\ge0$, the heat semigroup satisfies:
\begin{enumerate}[label=\textup{(\roman*)},leftmargin=2em]
\item conservative: $\displaystyle P_t\mathbf{1}=\mathbf{1}$.
\item positive preserving: for any nonnegative $f\in\ell^\infty(V)$, $P_tf\ge0$.
\item heat equation: for every $f\in\ell^\infty(V)$,
\[
 \partial_tP_tf=LP_tf=P_tLf,\qquad P_0f=f.
\]
\item self-adjoint: $\displaystyle \ip{P_tf}{g}=\ip{f}{P_tg}$ for $f,g\in\ell^2(V,m)$.
\item contraction: $\displaystyle \norm{P_tf}_{\ell^p(V,m)}\le\norm{f}_{\ell^p(V,m)}$ for $p\in\{1,\infty\},\ f\in\ell^p(V,m)$.
\end{enumerate}
For background on Dirichlet forms and reversible Markov semigroups, see
\cite{KellerLenz2012,LinLiu2018,Wojciechowski2009}.

We recall Duhamel's formula.
\begin{theorem}[{\cite[Chapter~VI, Section~7]{EngelNagel2000}}]
Let $A:D(A)\subseteq X\to X$ generate a strongly continuous semigroup
$(T(t))_{t\ge0}$ on a Banach space $X$. Let $T_0>0$ and
$v\in C^1([0,T_0],X)$ satisfy $v_t\in D(A)$ for every $t\in[0,T_0]$,
and assume that $t\mapsto Av_t$ is continuous.
Put $g_t=v'_t-Av_t$. Then
\begin{equation}\label{eq:duhamel}
 v_t=T(t-s)v_s+\int_s^tT(t-r)g_r\dd r
 \qquad(0\le s\le t\le T_0),
\end{equation}
where the integral is a Bochner integral in $X$.
\end{theorem}

\subsection{Bakry-E\'mery Curvature}
Define
\begin{equation*}
      \begin{aligned}
             \Gamma(f,g)(x)&=\frac12\sum_y p(x,y)(f(y)-f(x))(g(y)-g(x)),\\
              \Gamma_2(f,g)&=\frac12 \bigl(L\Gamma(f,g)-\Gamma(f,Lg)-\Gamma(g,Lf)\bigr).
      \end{aligned}
\end{equation*}
We abbreviate $\Gamma(f,f)=\Gamma f$ and
$\Gamma_2(f,f)=\Gamma_2 f$. For $K\in\R$, the Bakry--\'Emery curvature-dimension condition
$\CD(K,\infty)$ at a vertex $x$ means that
\begin{equation*}
 \Gamma_2(f)(x)\ge K\Gamma(f)(x)
 \quad\text{for every }f:V\to\R.
\end{equation*}
We say that the chain satisfies $\CD(K,\infty)$ if this holds at every vertex. 
\begin{remark}\label{rem:universal-curvature-bounds}
Every Markov chain satisfies $\CD(-1,\infty)$.
Indeed, direct expansion gives
\[
 \Gamma_2v(x)=\frac14\sum_{y,z}p(x,y)p(y,z)
 (v(z)-2v(y)+v(x))^2+\frac12(Lv(x))^2-\Gamma v(x).
\]

If $V$ contains at least two vertices, irreducibility gives
$a_x=1-p(x,x)>0$. For $e_x=\1_{\{x\}}$, we have
\[
 \Gamma e_x(x)=a_x/2,\qquad
 \Gamma_2e_x(x)=a_x^2/4+\frac34\sum_{y\ne x}p(x,y)p(y,x)\le a_x.
\]
Testing $\CD(K,\infty)$ with $e_x$ therefore gives
$K\le\Gamma_2e_x(x)/\Gamma e_x(x)\le2$.
\end{remark}

We record the fundamental semigroup inequalities under a
Bakry--\'Emery curvature lower bound; see \cite{LinLiu2018}.
For gradient estimates with unbounded graph Laplacians, see
\cite{HuaLin2017,KellerMunch2024}. Recall $I_K$ from~\eqref{eq:curvature-time}.
\begin{theorem}[Gradient and variance estimates]\label{thm:gradient}
Let $L$ be a bounded operator on $\ell^\infty(V)$ generating a
reversible continuous-time Markov chain, and assume $\CD(K,\infty)$
for some $K\in\R$. For every bounded real function
$f$ and every $t>0$,
\begin{equation}\label{eq:gradient-K}
 \Gamma P_tf\le e^{-2Kt}P_t\Gamma f,
\end{equation}
and
\begin{equation}\label{eq:variance-K}
 2I_K(t)\Gamma P_tf
 \le P_tf^2-(P_tf)^2
 \le2I_{-K}(t)P_t\Gamma f.
\end{equation}
\end{theorem}
\subsection{Information theory}
Let $\mu,\nu$ be probability measures on $V$. Their total variation
distance is defined by
\[
 \|\mu-\nu\|_{\mathrm{TV}}
 =\frac12\sum_{z\in V}|\mu(z)-\nu(z)|.
\]

The entropy of $\mu$ is defined by
\[
 H(\mu)=-\sum_{z\in V}\mu(z)\log\mu(z).
\]
By convention, $0\log0=0$.

If $\nu$ has full support, the relative entropy
(or Kullback--Leibler divergence) of $\mu$ with respect to $\nu$ is defined by
\[
 \D(\mu\Vert\nu)
 =\sum_{z\in V}\mu(z)\log\frac{\mu(z)}{\nu(z)}.
\]
By Jensen's inequality, $\D(\mu\Vert\nu)\ge0$.

The map $(\mu,\nu)\mapsto\D(\mu\Vert\nu)$ is convex
\cite[Theorem~5.1]{PolyanskiyWu2025}: for probability measures
$\mu_0,\mu_1,\nu_0,\nu_1$ on $V$, with $\nu_0,\nu_1$ having full support,
and $\lambda\in[0,1]$,
\begin{equation}\label{eq:entropy-joint-convexity}
 \begin{aligned}
 &\D((1-\lambda)\mu_0+\lambda\mu_1\Vert
       (1-\lambda)\nu_0+\lambda\nu_1)\\
 &\qquad\le(1-\lambda)\D(\mu_0\Vert\nu_0)
          +\lambda\D(\mu_1\Vert\nu_1).
 \end{aligned}
\end{equation}

 We write $\Var_\mu f$ and $\E_\mu f$ for the
variance and expectation of $f$ with respect to the probability measure
$\mu$, respectively.

We recall the variational characterization of relative entropy
\cite[Eq.~(4.11)]{PolyanskiyWu2025},
as well as the data processing inequality~\eqref{lem:data-processing}
\cite[Section~2.5]{PolyanskiyWu2025}. For background on information theory, see
\cite{CoverThomas2006, PolyanskiyWu2025}.

\begin{lemma}[Variational principle]\label{lem:entropy-variational}
Let $\mu,\nu$ be probability measures on $V$, with $\nu$ having full support.
The relative entropy admits the Donsker--Varadhan representation
\begin{equation}\label{eq:entropy-variational}
 \D(\mu\Vert \nu)
 =\sup_{g\in\ell^\infty(V)}
       \bigl(\E_\mu g-\log(\E_\nu e^g)\bigr).
\end{equation}
\end{lemma}

We also use relative entropy on arbitrary at most countable sets, with
the conventions $0\log(0/b)=0$ for $b\ge0$ and
$a\log(a/0)=+\infty$ for $a>0$.
Let $M:E\times F\to [0,1]$ be a Markov kernel. For a probability measure $\alpha$ on $E$, the distribution
$\alpha M$ on $F$ is defined by
\[
 (\alpha M)(b)=\sum_{a\in E}\alpha(a)M(a,b).
\]

\begin{lemma}[Data processing inequality]\label{lem:data-processing}
Let $M:E\times F\to [0,1]$ be a Markov kernel. For any probability measures $\alpha,\beta$
on $E$,
\begin{equation*}
 \D(\alpha M\Vert\beta M)\le\D(\alpha\Vert\beta).
\end{equation*}
\end{lemma}

\begin{corollary}
      Let $\mu,\nu$ be probability measures on an at most countable set $E$,
and let $B\subseteq E$. Then
\begin{equation}\label{eq:binary-data-processing}
 \D(\mu\|\nu)
 \ge\mu(B)\log\frac1{\nu(B)}-\log2.
\end{equation}
\end{corollary}
\begin{proof}
      Indeed, consider the Markov kernel $M:E\times\{0,1\}\to[0,1]$
defined by
\[
 M(x,1)=\1_B(x),\qquad M(x,0)=\1_{B^c}(x).
\]
Writing $p=\mu(B)$ and $q=\nu(B)$, we have
$\mu M=(1-p,p)$ and $\nu M=(1-q,q)$.
Lemma~\ref{lem:data-processing} gives
\[
 \begin{aligned}
 \D(\mu\|\nu)
 &\ge\D(\mu M\|\nu M)\\
 &=p\log\frac pq+(1-p)\log\frac{1-p}{1-q}.
 \end{aligned}
\]
Now \eqref{eq:binary-data-processing} follows from
$-p\log p-(1-p)\log(1-p)\le\log2$ and
$(1-p)\log(1/(1-q))\ge0$.
\end{proof}

\section{Log-Harnack inequality}
\label{sec:preliminary}

Boundedness of $L$ gives the
uniformisation formula
\[
 P_t=e^{-t}\sum_{j=0}^\infty\frac{t^j}{j!}P^j.
\]
Consequently, the probability measure $p_t(x,\cdot)$ has full support
for every $x\in V$ and $t>0$.

\subsection{Log-Harnack inequality}

\begin{lemma}
\label{lem:differential-harnack}
Assume that $p(x,y)\ge\eps$ whenever $x\sim y$.
Let $c:[a,b]\to(0,\infty)$ be continuous,
and $F:[a,b]\times V\to\R$ be such that $F(\cdot,z)$ is continuously differentiable on $[a,b]$
for every $z\in V$. Suppose that $\partial_rF_r(z)\ge c(r)\Gamma F_r(z)$ for $r\in(a,b),\ z\in V$. Then, for every $a\le s<t\le b$ and $x,y\in V$,
\begin{equation}\label{eq:differential-harnack}
 F_s(x)-F_t(y)
 \le\frac{2\dist(x,y)^2}{\eps\int_s^t c(r)\dd r}.
\end{equation}
\end{lemma}
\begin{proof}
By continuity, it suffices to consider $a<s<t<b$.
Since $\partial_rF_r\ge0$, the conclusion for $x=y$ is immediate.

For adjacent vertices $v\sim w$ and $s\le\sigma<\tau\le t$, set
\[
 g(r)=F_\sigma(v)-F_r(w)\qquad(\sigma\le r\le\tau).
\]
The function $g$ is nonincreasing. If $g(\tau)>0$, then
$g(r)>0$ and $F_r(v)-F_r(w)\ge g(r)$ throughout $[\sigma,\tau]$.
Using $p(w,v)\ge\eps$, we obtain
\[
 \begin{aligned}
 -g'(r)=\partial_rF_r(w)
 &\ge c(r)\Gamma F_r(w)\\
 &\ge\frac{\eps c(r)}2\bigl(F_r(v)-F_r(w)\bigr)^2
 \ge\frac{\eps c(r)}2g(r)^2.
 \end{aligned}
\]
Integrating $(1/g)'=-g'/g^2$ yields
\[
 \frac1{g(\tau)}
 \ge\frac1{g(\sigma)}+\frac\eps2\int_\sigma^\tau c(r)\dd r
 \ge\frac\eps2\int_\sigma^\tau c(r)\dd r.
\]
If $g(\tau)\le0$, the following upper bound is immediate as well.
Thus in all cases
\begin{equation}
      \label{equ:single-edge-bound}
       F_\sigma(v)-F_\tau(w)
 \le\frac2{\eps\int_\sigma^\tau c(r)\dd r}.
\end{equation}

For $x\ne y$, take a shortest path
$x=x_0\sim x_1\sim\cdots\sim x_n=y$, where $n=\dist(x,y)$.
Since $c$ is positive and continuous, we can choose
$s=t_0<t_1<\cdots<t_n=t$ such that
\[
 \int_{t_{i-1}}^{t_i}c(r)\dd r
 =\frac1n\int_s^t c(r)\dd r
 \qquad(1\le i\le n).
\]
Applying \eqref{equ:single-edge-bound} and summing gives
\[
 \begin{aligned}
 F_s(x)-F_t(y)
 &=\sum_{i=1}^n
   \bigl(F_{t_{i-1}}(x_{i-1})-F_{t_i}(x_i)\bigr)\\
 &\le\sum_{i=1}^n\frac2{\eps\int_{t_{i-1}}^{t_i}c(r)\dd r}
 =\frac{2n^2}{\eps\int_s^t c(r)\dd r}.
 \end{aligned}
\]
This proves \eqref{eq:differential-harnack}.
\end{proof}

For a bounded real-valued function $v$, put
\begin{equation}\label{eq:Psi}
 \Psi(v)(x)=\sum_y p(x,y)(e^{v(y)-v(x)}-1-(v(y)-v(x))),
\end{equation}
so that \[e^{-v}Le^v=Lv+\Psi(v).\]
This nonlinear operator and its associated calculus have been developed
and studied in~\cite{DierKassmannZacher2021,WeberZacherJFA,WeberZacher2023HeatKernel}. Fix $\tau=1/3$. 

\begin{lemma}\label{lem:shortvariance}
There exists $C=C(\eps)\in(0,\infty)$ such that  every bounded real $v$ satisfies
\begin{equation}\label{eq:shortvariance}
 \Var_{p_\tau(x,\cdot)}v
 \le P_\tau((v-v(x))^2)(x)\le C P_\tau\Psi(v)(x).
\end{equation}
\end{lemma}
\begin{proof}
Let $a_n=e^{-\tau}\tau^n/n!$, and let
$(Y_j)_{j\ge0}$ be a random walk with transition kernel $P$ and $Y_0=x$.
Weighted Cauchy--Schwarz gives
\[
 (v(Y_n)-v(x))^2\le\frac32\sum_{j=1}^n3^{n-j}
                   (v(Y_j)-v(Y_{j-1}))^2.
\]
Taking expectation with respect to the random walk and using
$a_j^{-1}\sum_{n\ge j}3^{n-j}a_n
=\sum_{k\ge0}j!/(j+k)!\le e-1$ for $j\ge1$, we obtain
\[
 \begin{gathered}
 \begin{aligned}
 P_\tau((v-v(x))^2)(x)
 &=\sum_{n\ge1}a_n\E[(v(Y_n)-v(x))^2]\\
 &\le\frac32\sum_{n\ge1}a_n\sum_{j=1}^n3^{n-j}
             \E[(v(Y_j)-v(Y_{j-1}))^2]\\
 &=\frac32\sum_{j\ge1}\left(\sum_{n\ge j}3^{n-j}a_n\right)
             \E[(v(Y_j)-v(Y_{j-1}))^2]\\
 &\le\frac{3(e-1)}2\sum_{j\ge1}a_j
                \E[(v(Y_j)-v(Y_{j-1}))^2]\\
 &=\frac{3(e-1)}2\sum_{u,w}
       \left(p(u,w)\sum_{j\ge1}a_jP^{j-1}(x,u)\right)
       (v(w)-v(u))^2.
 \end{aligned}
 \end{gathered}
\]
Let $c_{uw}=p(u,w)\sum_{j\ge1}a_jP^{j-1}(x,u)$. The inequalities $a_j\le\tau a_{j-1}$
give $c_{uw}\le\tau p_\tau(x,u)p(u,w)$ and $c_{uw}\le p_\tau(x,w)$.

On an unordered edge, set $A=p_\tau(x,u)p(u,w)$,
$B=p_\tau(x,w)p(w,u)$ and $S=c_{uw}+c_{wu}$.
Since $p(u,w),p(w,u)\ge\eps$, we have
\[
 S\le(\eps^{-1}+\tau)\min\{A,B\}.
\]
Writing $\Upsilon(b)=e^b-1-b$, we have
$\Upsilon(b)+\Upsilon(-b)\ge b^2$. Consequently,
\[
 Sb^2\le(\eps^{-1}+\tau)
            \bigl(A\Upsilon(b)+B\Upsilon(-b)\bigr).
\]
Group the quadratic sum by unordered edges and apply this bound with
$b=v(w)-v(u)$. The resulting exponential sum is $P_\tau\Psi(v)(x)$,
which proves the upper bound for some $C$ depending only on
$\eps$.
\end{proof}

\begin{lemma}\label{lem:poisson-insertion}
For $n=\dist(x,y)\ge1$ and $t>0$,
\begin{equation}\label{eq:shorttime-row-entropy}
 \D(p_t(x,\cdot)\|p_t(y,\cdot))
 \le\sum_{j=1}^n\log\frac{t+j}{\eps t}.
\end{equation}
\end{lemma}
\begin{proof}
A path from $y$ to $x$ gives $P^n(y,x)\ge\eps^n$.
Set $a_j=e^{-t}t^j/j!$. Uniformisation implies
\[
 p_t(y,z)\ge\sum_{j\ge0}a_jP^j(x,z)\frac{(\eps t)^n j!}{(j+n)!}.
\]

Since $\log(1+k/j)\le k$ for $j\ge1$, the identities
$\sum_{k\ge0}a_k=1$ and $\sum_{k\ge0}ka_k=t$ imply
\[
 \begin{aligned}
 &\sum_z\sum_{k\ge0}a_kP^k(x,z)
       \sum_{j=1}^n\left|\log\frac{k+j}{\eps t}\right|\\
 &=\sum_{k\ge0}a_k\sum_{j=1}^n
       \left|\log\frac{k+j}{\eps t}\right|
 \le nt+\sum_{j=1}^n\left|\log\frac{j}{\eps t}\right|<\infty.
 \end{aligned}
\]
Thus Fubini's theorem applies to the logarithmic sums below.
The negative part of the relative-entropy sum is also summable:
$u\log(v/u)\le v/e$ for $0<u\le v$, and
$\sum_zp_t(y,z)=1$.

For each $z$, apply Jensen to $-\log$ with the probability weights
$a_kP^k(x,z)/p_t(x,z)$. Summing in $z$ and then applying Jensen
to $\log$ gives
\[
 \begin{aligned}
 \D(p_t(x,\cdot)\|p_t(y,\cdot))
 &\le\sum_z\sum_{k\ge0}a_kP^k(x,z)
                  \log\frac{(k+n)!}{k!(\eps t)^n}\\
 &=\sum_{k\ge0}a_k\log\frac{(k+n)!}{k!(\eps t)^n}\\
 &=\sum_{j=1}^n\sum_{k\ge0}a_k\log\frac{k+j}{\eps t}\\
 &\le\sum_{j=1}^n\log\frac{\sum_{k\ge0}ka_k+j}{\eps t}\\
 &=\sum_{j=1}^n\log\frac{t+j}{\eps t}.
 \end{aligned}
\]
\end{proof}

\begin{proof}[Proof of Theorem~\ref{thm:logharnack}]
For $x=y$, both sides of \eqref{eq:logharnack} are zero. For $x\ne y$, the
variational formula~\eqref{eq:entropy-variational} reduces the proof
to the function inequality. Theorem~\ref{thm:gradient} and
Lemma~\ref{lem:shortvariance} give
\begin{equation}\label{eq:averagedPsi}
 \Gamma P_rv\le\frac{C}{2I_K(\tau)}e^{-2K(r-\tau)}P_r\Psi(v)
 \qquad(r\ge\tau).
\end{equation}

For $t>\tau$, put $F_s=P_{t-s}\log P_sf$,
where $f$ is bounded and $\inf f>0$. Then, for $0\le s\le t-\tau$,
\[
 \partial_sF_s=P_{t-s}\Psi(\log P_sf)\ge c(s)\Gamma F_s,
 \qquad c(s)=\frac{2I_K(\tau)}{C}e^{2K(t-\tau-s)}>0.
\]
Lemma~\ref{lem:differential-harnack} and
$F_{t-\tau}=P_\tau\log P_{t-\tau}f\le\log P_tf$ give
\[
 \begin{aligned}
 P_t\log f(x)-\log P_tf(y)
 &\le F_0(x)-F_{t-\tau}(y)\\
 &\le\frac{2\dist(x,y)^2}{\eps\int_0^{t-\tau}c(s)\dd s}\\
 &=\frac{C\dist(x,y)^2}{\eps I_K(\tau)I_K(t-\tau)}.
 \end{aligned}
\]

By Remark~\ref{rem:universal-curvature-bounds}, the chain satisfies
$\CD(-1,\infty)$ and $K\le2$. Since $I_K(t)\le I_{-1}(t)$
for $K<-1$, it suffices to consider
$-1\le K\le2$.
For $t\ge2\tau$, the identity for $I_K(t)$ gives
\[
 \begin{aligned}
 \frac{I_K(t)}{I_K(\tau)I_K(t-\tau)}
 &=\frac1{I_K(t-\tau)}+\frac{e^{2K\tau}}{I_K(\tau)}\\
 &\le2K\coth(K/3)\le4\coth(2/3).
 \end{aligned}
\]
Here the middle bound has value $6$ at $K=0$, and
$z\coth z$ is even and increasing for $z>0$. The preceding bound
gives \eqref{eq:logharnack} for
$t\ge2\tau$, with a constant depending only on $\eps$.

For $0<t\le2/3$, Lemma~\ref{lem:poisson-insertion} gives
$\D(p_t(x,\cdot)\|p_t(y,\cdot))\le\dist(x,y)^2/(\eps t)$, while
\[
 \frac{I_K(t)}t\le\frac{e^{8/3}-1}{8/3}.
\]
Increasing $C_{\mathrm{LH}}=C_{\mathrm{LH}}(\eps)$ covers these times.
\end{proof}

\subsection{Volume bounded by entropy of heat flow}

In this subsection, We construct the relationship between entropy and volume. As an application, we prove a coarsed polynomial growth under non-negative  Bakry-\'Emery Curvature, which will be used later. Write
\[
 \mathcal H_t(x)=-\sum_zp_t(x,z)\log\frac{p_t(x,z)}{m(z)}.
\]
We allow $\mathcal{H}_t(x)=\infty$.

Let $M$ be a Markov kernel on $V$ preserving $m$, that is,
\[
 \sum_{x\in V}m(x)M(x,z)=m(z)\qquad(z\in V).
\]
Then $q_t=P_tM$ satisfies
\[
 \sum_zq_t(x,z)=1,\qquad
 \sum_xm(x)q_t(x,z)=m(z).
\]

\begin{lemma}\label{lem:entropy-volume}
Let $M$ be a Markov kernel preserving $m$, and set $q_t=P_tM$ and
\[
 \mathcal H_t^M(o)=-\sum_zq_t(o,z)\log\frac{q_t(o,z)}{m(z)}.
\]
For every $t>0$ and $R\ge0$,
\begin{equation}\label{eq:curved-volume-entropy-transfer}
 \log\vol(B(o,R))\le\mathcal H_t^M(o)+\frac{C_{\mathrm{LH}}R^2}{I_K(t)}.
\end{equation}
\end{lemma}
\begin{proof}
Set $V_R=\vol(B(o,R))$ and let
\[
 \nu_t(z)=V_R^{-1}\sum_{x\in B(o,R)}m(x)q_t(x,z).
\]
Since both $P_t$ and $M$ preserve $m$, so does $q_t$; hence
$\nu_t(z)\le m(z)/V_R$. Convexity in the second argument of relative
entropy, the data processing inequality~\eqref{lem:data-processing},
and the log-Harnack bounds give
\[
 \begin{aligned}
 \log V_R-\mathcal H_t^M(o)
 &\le\D(q_t(o,\cdot)\|\nu_t)\\
 &\le V_R^{-1}\sum_{x\in B(o,R)}m(x)
              \D(q_t(o,\cdot)\|q_t(x,\cdot))\\
 &\le V_R^{-1}\sum_{x\in B(o,R)}m(x)
              \D(p_t(o,\cdot)\|p_t(x,\cdot))
 \le\frac{C_{\mathrm{LH}}R^2}{I_K(t)}.
 \end{aligned}
\]
\end{proof}

\begin{lemma}\label{lem:heatentropy}
Let $M$ be a Markov kernel on $V$ preserving $m$. For every $x\in V$ and $t\ge0$,
\begin{equation}\label{eq:curved-entropy}
 \mathcal H_t^M(x)\le P_t\mathcal H_0^M(x)
       +C_\eps+C_{\mathrm{LH}}\log(1+t)+2C_{\mathrm{LH}}K_-t.
\end{equation}
where $C_\varepsilon$ is a constant dependent only on $\varepsilon$.
If $q_t(x,\cdot)$ has finite entropy for every $x\in V$ and
$t>0$, then
\begin{equation}\label{eq:heat-entropy-evolution}
 (\partial_t-L)\mathcal H_t^M(x)
 =\sum_yp(x,y)\D\bigl(q_t(y,\cdot)\|q_t(x,\cdot)\bigr)
 \qquad(t>0).
\end{equation}
\end{lemma}
\begin{proof}
Fix a finite partition $\pi=\{A_1,\ldots,A_q\}$ of $V$ into nonempty
sets. For $\rho>0$ define the Markov kernel
\[
 M^\prime_{\pi,\rho}(z,i)=\frac{\1_{A_i}(z)+\rho}{1+q\rho},
\]
and write
\[
 v_i^\rho(t,x)=(q_tM^\prime_{\pi,\rho})(x,i)
     =\frac{\sum_{z\in A_i}q_t(x,z)+\rho}{1+q\rho},\qquad
 H_t^{\pi,\rho}(x)=-\sum_{i=1}^qv_i^\rho(t,x)\log v_i^\rho(t,x).
\]
These probabilities solve $\partial_tv_i^\rho=Lv_i^\rho$, sum to
one, and are bounded below by $\rho/(1+q\rho)$. Thus the entropy can
be differentiated in $\ell^\infty(V)$. Expanding $L$ gives
\[
 \begin{aligned}
 F_t^{\pi,\rho}(x)
 &:=(\partial_t-L)H_t^{\pi,\rho}(x)\\
 &=-\sum_i(1+\log v_i^\rho(t,x))Lv_i^\rho(t,x)-LH_t^{\pi,\rho}(x)\\
 &=\sum_yp(x,y)\sum_i v_i^\rho(t,y)
                         \log\frac{v_i^\rho(t,y)}{v_i^\rho(t,x)}\\
 &=\sum_yp(x,y)\D\bigl(v^\rho(t,y)\|v^\rho(t,x)\bigr).
 \end{aligned}
\]

For adjacent vertices and $0<t\le1$, the bound
\eqref{eq:shorttime-row-entropy} in Lemma~\ref{lem:poisson-insertion} gives
\[
 \D\bigl(p_t(y,\cdot)\|p_t(x,\cdot)\bigr)
 \le\log\frac{t+1}{\eps t}
 \le t\log2+\log\frac2{\eps t}.
\]
Together with the data processing inequality~\eqref{lem:data-processing}
and Theorem~\ref{thm:logharnack}, this
implies, for $x\sim y$ and $t>0$,
\[
 \begin{aligned}
 \D\bigl(v^\rho(t,y)\|v^\rho(t,x)\bigr)
 &\le\D\bigl(q_t(y,\cdot)\|q_t(x,\cdot)\bigr)\\
 &\le\D\bigl(p_t(y,\cdot)\|p_t(x,\cdot)\bigr)
 \le b_{\eps,K_-}(t),
 \end{aligned}
\]
where
\[
 b_{\eps,K_-}(t)=
 \begin{cases}
 t\log2+\log\dfrac{2}{\eps t},&0<t\le1,\\[3pt]
 C_{\mathrm{LH}}\bigl(t^{-1}+2K_-\bigr),&t>1.
 \end{cases}
\]
The long-time bound uses Theorem~\ref{thm:logharnack} under
$\CD(K,\infty)$ and the bound
$C_{\mathrm{LH}}/I_K(t)\le C_{\mathrm{LH}}(t^{-1}+2K_-)$.
Since the $\D\bigl(v^\rho(t,x)\|v^\rho(t,x)\bigr)$ vanishes and $\sum_{y\ne x}p(x,y)\le1$,
\[
 0\le F_t^{\pi,\rho}(x)\le b_{\eps,K_-}(t).
\]

Since $0\le H_0^{\pi,\rho}\le\log q$, Duhamel's formula gives
\[
 H_t^{\pi,\rho}(x)
 =P_tH_0^{\pi,\rho}(x)+\int_0^tP_{t-s}F_s^{\pi,\rho}(x)\dd s
 \le P_tH_0^{\pi,\rho}(x)+\int_0^t b_{\eps,K_-}(s)\dd s.
\]
The short-time bound is integrable and
$\int_0^1b_{\eps,K_-}(s)\dd s=C_\eps$.

Letting $\rho\downarrow0$, continuity of entropy and dominated
convergence give
\[
 \begin{aligned}
 H_t^\pi(x)
 &:=-\sum_{A\in\pi}\left(\sum_{z\in A}q_t(x,z)\right)
                  \log\left(\sum_{z\in A}q_t(x,z)\right)\\
 &\le P_tH_0^\pi(x)+\int_0^t b_{\eps,K_-}(s)\dd s.
 \end{aligned}
\]
As finite partitions $\pi_N$ refine monotonically to singletons,
we have, since $q_0=M$,
\[
 H_t^{\pi_N}(x)\uparrow H(q_t(x,\cdot)),\qquad
 H_0^{\pi_N}(y)\uparrow H(M(y,\cdot)).
\]
Monotone convergence therefore gives
\[
 H(q_t(x,\cdot))
 \le\sum_yp_t(x,y)H(M(y,\cdot))
       +\int_0^tb_{\eps,K_-}(s)\dd s.
\]
Adding $P_tM(\log m)$ and integrating the long-time
bound proves \eqref{eq:curved-entropy}.

For the evolution identity \eqref{eq:heat-entropy-evolution}, assume
$q_t(x,\cdot)$ has finite entropy for every $x\in V$ and $t>0$,
and fix $0<a<b<T$. Uniformisation gives, for $a\le s\le b$,
\[
 q_s(x,z)\ge e^{a-b}q_a(x,z),\qquad
 |\partial_sq_s(x,z)|\le C_{a,b,T}q_T(x,z).
\]

For fixed $x$, define a Markov kernel from $\{0,1,\ldots\}$ to $V$ by
\[
 \widetilde M_x(j,z)=(P^jM)(x,z).
\]
It satisfies $\Pois(t)\widetilde M_x=q_t(x,\cdot)$, where
$\Pois(t)$ denotes the Poisson distribution with mean $t$:
\[
 \Pois(t)(j)=e^{-t}\frac{t^j}{j!}\qquad(j=0,1,\ldots).
\]
The data processing inequality~\eqref{lem:data-processing}, applied with
$\widetilde M_x$, also gives
\[
 \D(q_T(x,\cdot)\|q_a(x,\cdot))
 \le\D(\Pois(T)\|\Pois(a))<\infty.
\]

Since $q_a(x,z)\le1$, finite entropy of $q_T(x,\cdot)$ gives
\[
 \sum_zq_T(x,z)|\log q_a(x,z)|
 =H(q_T(x,\cdot))
  +\D(q_T(x,\cdot)\|q_a(x,\cdot))<\infty.
\]
Moreover, the uniformisation estimates imply, for $a\le s\le b$,
\[
 \begin{aligned}
 \left|\partial_s\bigl(q_s(x,z)\log q_s(x,z)\bigr)\right|
 &=|1+\log q_s(x,z)|\,|\partial_sq_s(x,z)|\\
 &\le C_{a,b,T}q_T(x,z)
       \bigl(1+b-a+|\log q_a(x,z)|\bigr).
 \end{aligned}
\]
The right-hand side is summable in $z$ and independent of $s$,
so differentiation may be interchanged with summation on $[a,b]$.

Finally,
\[
 \mathcal H_t^M(x)=H(q_t(x,\cdot))+q_t(\log m)(x).
\]
Since $\log m$ is bounded, so is $M(\log m)$, and
\[
 (\partial_t-L)q_t(\log m)
 =(\partial_t-L)P_t(M(\log m))=0.
\]
The preceding finite-dimensional calculation therefore proves
\eqref{eq:heat-entropy-evolution}.
\end{proof}

\begin{proposition}[Coarse volume growth]\label{prop:coarse-volume-growth}
Assume \eqref{eq:geometry} and $\CD(K,\infty)$.
There are constants $C,D<\infty$, depending only on $\eps,\kappa$,
such that, for every $o\in V$ and $R\ge0$,
\begin{equation}\label{eq:curved-coarse-volume}
 \vol(B(o,R))\le m_-C(1+R)^D
 \qquad\text{whenever }K_-R^2\le4.
\end{equation}
\end{proposition}
\begin{proof}
For $R\ge1$, apply Lemma~\ref{lem:entropy-volume} with $M=I$ and
$t=R^2$. Since $\mathcal H_0^I=\log m$ and $m\le\kappa m_-$,
\[
 P_t\mathcal H_0^I(o)=P_t(\log m)(o)\le\log(\kappa m_-).
\]
Thus \eqref{eq:curved-entropy} and
\eqref{eq:curved-volume-entropy-transfer}, together with
$C_{\mathrm{LH}}/I_K(t)\le C_{\mathrm{LH}}(t^{-1}+2K_-)$, give
\[
 \log\frac{\vol(B(o,R))}{m_-}
 \le C_\eps+\log\kappa+C_{\mathrm{LH}}\log(1+R^2)
                     +C_{\mathrm{LH}}(1+4{K_-} R^2).
\]
If ${K_-} R^2\le4$, this proves \eqref{eq:curved-coarse-volume} with
$D=2C_{\mathrm{LH}}$ and $C=\kappa e^{C_\eps+17C_{\mathrm{LH}}}$.

For $0\le R<1$, the ball is the singleton $\{o\}$ and
\eqref{eq:curved-coarse-volume} follows from $m(o)\le\kappa m_-$.
\end{proof}

\Needspace{11\baselineskip}
\section{Averaged heat displacement and volume doubling}
Throught this section, we always assume \eqref{eq:geometry} and $\CD(K,\infty)$ for $K\in\R$. Let $K_-=\max\{-K,0\}$.
\label{sec:modified}

\subsection{An auxiliary flow}
We recall the bounded maximum principle of Guo, Huang, and Huang
\cite[Lemma~A.1]{GuoHuangHuang}.

\begin{lemma}[Bounded maximum principle]
\label{lem:bounded-maximum-principle}
Let $T>0$ and $h\in C^1([0,T],\ell^\infty(V))$. Suppose that
nonnegative coefficients $a_t(x,y)$ satisfy
\[
 \sup_{\substack{0\le t\le T\\x\in V}}\sum_y a_t(x,y)<\infty,
\]
and that, for every $x\in V$ and $0\le t\le T$,
\[
 h'_t(x)\le\sum_y a_t(x,y)\bigl(h_t(y)-h_t(x)\bigr).
\]
Then $t\mapsto\sup_{x\in V}h_t(x)$ is nonincreasing on $[0,T]$.
\end{lemma}

We use a gradient estimate. Its proof
is given in Appendix~\ref{sec:moment}.

\begin{proposition}
\label{prop:moment}\label{prop:curved-signed-moment}
Assume \eqref{eq:geometry} and $\CD(K,\infty)$ for $K\in\R$.
There are constants $0<\delta<1$ and $C<\infty$, depending
only on $\eps,\kappa$, with the following property. Let $h>0$ belong to
$\ell^1(m)$, let $\beta\ge0$, and set
\[
 a=\frac{(\beta-L)h}{h},\qquad
 M=\sum_xm(x)h_x|a_x|^3,
 \qquad S=\sum_xm(x)h_x(\Gamma\log h(x))^3.
\]
The function $a$ may have either sign. If
$\sup_{x\sim y}|\log h_x-\log h_y|\le\delta$, then
\begin{equation}\label{eq:moment}
 S\le C(M+{K_-}^3\norm{h}_1).
\end{equation}
Consequently, for every bounded $1$-Lipschitz function $f$,
\begin{equation}\label{eq:source-flux}
 \begin{aligned}
 |\ip h{Lf}|=|\cE(h,f)|
 &\le\frac1{\sqrt2}\norm{h}_1^{5/6}S^{1/6}\le C\norm{h}_1^{5/6}(M+{K_-}^3\norm{h}_1)^{1/6}.
 \end{aligned}
\end{equation}
\end{proposition}
\begin{lemma}
\label{prop:initial-density}\label{lem:curved-initial-density}
There are $c_0,\alpha_0>0$, depending only on $\eps,\kappa$, such that
if $0<\alpha\le\alpha_0$ and ${K_-}\le\alpha$, then, for every $o\in V$,
\[
 \sigma=c_0\alpha^{1/3},\qquad
 Z_\sigma=\sum_xm(x)e^{-\sigma\dist(o,x)},\qquad
 \phi_\alpha(x)=Z_\sigma^{-1}e^{-\sigma\dist(o,x)}
\]
are well defined and satisfy
\begin{equation}\label{eq:curved-initial-moment}
 \norm{\phi_\alpha}_1=1,\qquad
 \ip{\phi_\alpha}{\dist(o,\cdot)}\le C\alpha^{-1/2}.
\end{equation}
We may choose $c_0=\sqrt\eps\,\delta/(4e)$ and $\alpha_0\le1$.
\end{lemma}
\begin{proof}
Fix this value of $c_0$ and let $\ell=\lceil\alpha^{-1/2}\rceil$.
Since ${K_-}\ell^2\le4$, \eqref{eq:curved-coarse-volume} implies
\[
 \sup_x|B(x,\ell)|\le B_\alpha:=\max\{2,C(1+\ell)^D\}.
\]
Splitting a geodesic into pieces of length at most the integer $\ell$
gives $|B(o,n\ell)|\le B_\alpha^n$ for every integer $n\ge1$.
As $\alpha\downarrow0$,
\[
 \sigma\ell\ge c_0\alpha^{-1/6},\qquad
 \log B_\alpha=O_{\eps,\kappa}(\log(2/\alpha)).
\]
Thus $\alpha_0$ can be chosen so that
$\sigma\ell/2\ge3\log B_\alpha$ whenever $\alpha\le\alpha_0$.
Shrinking $\alpha_0$ further if necessary, we also have $\sigma\le1$.

Fix $\lambda\in\{\sigma,\sigma/2\}$, so that $0<\lambda\le1$.
Split the sum defining $Z_\lambda$ over $B(o,\ell)$ and the sets
\[
 A_n=\{x\in V:n\ell<\dist(o,x)\le(n+1)\ell\},\qquad n\ge1.
\]

For the inner ball, write
\[
 e^{-\lambda\dist(o,x)}
 =\lambda\int_{\dist(o,x)}^\infty e^{-\lambda s}\dd s.
\]
Summing over $x\in B(o,\ell)$ and interchanging the nonnegative sum and
integral, we obtain
\[
 \sum_{x\in B(o,\ell)}m(x)e^{-\lambda\dist(o,x)}
 =\lambda\int_0^\infty e^{-\lambda s}
       \vol\bigl(B(o,\min\{s,\ell\})\bigr)\dd s.
\]
Since $K_-\min\{s,\ell\}^2\le K_-\ell^2\le4$,
\eqref{eq:curved-coarse-volume} bounds the volume in the integrand by
$Cm_-(1+\min\{s,\ell\})^D\le Cm_-(1+s)^D$.
Thus, changing variables $u=\lambda s$ and using $\lambda\le1$, we obtain
\[
 \begin{aligned}
 \sum_{x\in B(o,\ell)}m(x)e^{-\lambda\dist(o,x)}
 &\le Cm_-\lambda\int_0^\infty e^{-\lambda s}(1+s)^D\dd s\\
 &=Cm_-\lambda^{-D}\int_0^\infty e^{-u}(\lambda+u)^D\dd u\\
 &\le Cm_-\lambda^{-D}\int_0^\infty e^{-u}(1+u)^D\dd u
 \le Cm_-\lambda^{-D}.
 \end{aligned}
\]
Here the last integral depends only on $D$.

For each $n\ge1$, $m(x)\le\kappa m_-$ and
$|A_n|\le|B(o,(n+1)\ell)|\le B_\alpha^{n+1}$ imply
\[
 \sum_{x\in A_n}m(x)e^{-\lambda\dist(o,x)}
 \le\kappa m_-B_\alpha^{n+1}e^{-\lambda n\ell}.
\]
Moreover, $\lambda\ell\ge\sigma\ell/2\ge3\log B_\alpha$, hence
$e^{-\lambda n\ell}\le B_\alpha^{-3n}$. Summing over $n\ge1$ and
using $B_\alpha\ge2$ therefore gives
\[
 \begin{aligned}
 \sum_{n\ge1}\sum_{x\in A_n}m(x)e^{-\lambda\dist(o,x)}
 &\le\kappa m_-\sum_{n\ge1}B_\alpha^{1-2n}\\
 &=\kappa m_-\frac{B_\alpha^{-1}}{1-B_\alpha^{-2}}
 \le\frac23\kappa m_-\le Cm_-.
 \end{aligned}
\]
Combining the two bounds and using $\lambda\le1$ yields
$Z_\lambda\le Cm_-(\lambda^{-D}+1)\le Cm_-\lambda^{-D}$.
The term $x=o$ gives $Z_\lambda\ge m(o)\ge m_-$.

Since $0<Z_\sigma<\infty$, the definition of $\phi_\alpha$ gives
$\sum_xm(x)\phi_\alpha(x)=1$. For every integer $N\ge1$, Jensen's
inequality gives
\[
 \frac\sigma2\ip{\phi_\alpha}{\dist(o,\cdot)\wedge N}
 \le\log\left(\sum_xm(x)\phi_\alpha(x)
       e^{\frac\sigma2(\dist(o,x)\wedge N)}\right).
\]
Using $\dist(o,x)\wedge N\le\dist(o,x)$ and the definition of
$\phi_\alpha$, the sum inside the logarithm satisfies
\[
 \begin{aligned}
 \sum_xm(x)\phi_\alpha(x)e^{\frac\sigma2(\dist(o,x)\wedge N)}
 &\le\frac1{Z_\sigma}\sum_xm(x)e^{-\frac\sigma2\dist(o,x)}=\frac{Z_{\sigma/2}}{Z_\sigma}.
 \end{aligned}
\]
Taking $N\to\infty$ and using monotone convergence therefore yields
\[
 \ip{\phi_\alpha}{\dist(o,\cdot)}
 \le\frac2\sigma\log\frac{Z_{\sigma/2}}{Z_\sigma}.
\]

The bounds on $Z_\lambda$, with $C\ge1$, imply
\[
 \log\frac{Z_{\sigma/2}}{Z_\sigma}
 \le\log\bigl(C(\sigma/2)^{-D}\bigr)
 =\log C+D\log(2/\sigma)
 \le C'\log(2/\alpha),
\]
where we used $\sigma=c_0\alpha^{1/3}$ and the fact that $c_0,D$
depend only on $\eps,\kappa$. Consequently,
\[
 \ip{\phi_\alpha}{\dist(o,\cdot)}
 \le C\alpha^{-1/3}\log(2/\alpha).
\]
Finally, $\log t\le6t^{1/6}$ for $t\ge1$, so that
\[
 \alpha^{-1/3}\log(2/\alpha)
 \le6\,2^{1/6}\alpha^{-1/2}
 \qquad(0<\alpha\le1).
\]
\end{proof}

\begin{lemma}[An auxiliary flow]
\label{prop:auxiliary}\label{lem:curved-auxiliary}
Suppose $0<\alpha\le\alpha_0$, ${K_-}\le\alpha$, and put $T=\alpha^{-1}$.
There are a positive curve $q\in C^1([0,T],\ell^1(m))$ with
$q_0=\phi_\alpha$, and a potential $W_t$, such that
\begin{equation}\label{eq:qflow}
 q_t'=Lq_t+W_tq_t,\qquad 0\le W_t\le\alpha/8,\qquad
 q_t\ge P_t\phi_\alpha,\qquad \norm{q_t}_1\le e^{\alpha t/8}.
\end{equation}
For every $0\le t\le T$,
$\sup_{x\sim y}|\log q_t(x)-\log q_t(y)|\le\delta$.
\end{lemma}
\begin{proof}
Set $u_0=\log\phi_\alpha$. Since $u_0$ has bounded edge
increments, the map $u\mapsto Lu+\Gamma u$ is locally Lipschitz
from $u_0+\ell^\infty(V)$ to $\ell^\infty(V)$, so
$u'=Lu+\Gamma u$ has a unique local solution. Set $G_t=\Gamma u_t$.
Differentiating and using $u_t'=Lu_t+G_t$ and the definition of
$\Gamma_2$, we obtain
\[
 \begin{aligned}
 \partial_tG_t
 &=2\Gamma(u_t,u_t')\\
 &=2\Gamma(u_t,Lu_t)+2\Gamma(u_t,G_t)\\
 &=LG_t-2\Gamma_2(u_t)+2\Gamma(u_t,G_t).
 \end{aligned}
\]
The curvature condition $\CD(K,\infty)$ and $G_t\ge0$ give
\[
 -2\Gamma_2(u_t)\le-2KG_t\le2K_-G_t.
\]
Moreover, the definitions of $L$ and $\Gamma$ imply
\[
 \begin{aligned}
 LG_t(x)+2\Gamma(u_t,G_t)(x)=\sum_y p(x,y)\bigl(1+u_t(y)-u_t(x)\bigr)
                          \bigl(G_t(y)-G_t(x)\bigr).
 \end{aligned}
\]
Combining these identities and inequalities yields
\[
 \partial_tG_t(x)\le
 \sum_yp(x,y)(1+u_t(y)-u_t(x))(G_t(y)-G_t(x))+2{K_-} G_t(x).
\]

On any interval on which $|u_t(y)-u_t(x)|<1/2$ on edges, we have
\[
 \frac12p(x,y)
 \le p(x,y)\bigl(1+u_t(y)-u_t(x)\bigr)
 \le\frac32p(x,y).
\]
Thus these coefficients are nonnegative and their sum over $y$ is
at most $3/2$.
Applying Lemma~\ref{lem:bounded-maximum-principle} to
$e^{-2{K_-} t}G_t$ gives
\[
 \begin{aligned}
 \norm{G_t}_\infty&\le e^{2{K_-} t}\sigma^2/2,\quad\sup_{x\sim y}|u_t(y)-u_t(x)|
 \le e^{{K_-} t}\sigma/\sqrt\eps
 \leq\frac{\delta\alpha^{1/3}}4
 \quad(0\le t\le T).
 \end{aligned}
\]
This closes the bootstrap; boundedness of $u_t'$ continues the
solution to $T$.

Define $q_t=e^{u_t+\alpha t/16}$ and
$W_t=\alpha/16+\Gamma u_t-\Psi(u_t)$. Taylor's formula gives
\[
 |\Psi(u_t)-\Gamma u_t|
 \le\frac{e^{1/2}}6
       \left(\frac{\delta\alpha^{1/3}}4\right)^3
 \le\alpha/16.
\]
The preceding estimate and the definition of $W_t$ give
$0\le W_t\le\alpha/8$.

To justify the equation in $\ell^1(m)$, write
\[
 q_t=\phi_\alpha\exp\bigl(u_t-u_0+\alpha t/16\bigr).
\]
Since $u-u_0\in C^1([0,T],\ell^\infty(V))$, the exponential factor
also belongs to $C^1([0,T],\ell^\infty(V))$.
Multiplication by $\phi_\alpha$ is bounded from $\ell^\infty(V)$ to
$\ell^1(m)$, because
\[
 \norm{\phi_\alpha g}_1
 \le\norm{g}_\infty\norm{\phi_\alpha}_1
 =\norm{g}_\infty.
\]
Consequently, $q\in C^1([0,T],\ell^1(m))$, with
$q_0=\phi_\alpha$ and $q_t'=(u_t'+\alpha/16)q_t$.
The bounded edge increments of $u_t$ allow us to apply the exponential
identity pointwise, giving
\[
 \begin{aligned}
 \frac{q_t'}{q_t}-\frac{Lq_t}{q_t}
 &=u_t'+\alpha/16-Lu_t-\Psi(u_t)\\
 &=\alpha/16+\Gamma u_t-\Psi(u_t)=W_t.
 \end{aligned}
\]
Since $L=P-I$ is bounded on $\ell^1(m)$ and $W_t$ is bounded,
this proves $q_t'=Lq_t+W_tq_t$ in $\ell^1(m)$.

In particular, $t\mapsto W_tq_t=q_t'-Lq_t$ is continuous in
$\ell^1(m)$, so Duhamel's formula applies and yields
\[
 q_t=P_t\phi_\alpha+\int_0^tP_{t-s}(W_sq_s)\dd s
 \ge P_t\phi_\alpha,
\]
where the inequality follows from $W_sq_s\ge0$ and positivity of
$P_{t-s}$.

Moreover, $q_t>0$ and mass conservation give
\[
 \frac{\dd}{\dd t}\norm{q_t}_1
 =\ip{q_t'}{\1}
 =\ip{W_tq_t}{\1}
 \le\frac\alpha8\norm{q_t}_1.
\]
Thus $e^{-\alpha t/8}\norm{q_t}_1$ is nonincreasing, and
$\norm{q_0}_1=1$ gives $\norm{q_t}_1\le e^{\alpha t/8}$.
This proves \eqref{eq:qflow}.
Finally, $\log q_t(y)-\log q_t(x)=u_t(y)-u_t(x)$, so the preceding
edge estimate also gives the required bound by $\delta$.
\end{proof}

\subsection{Averaged heat displacement}

\begin{lemma}
\label{prop:integrated-flow}\label{lem:curved-laplace-mixture}
Under the assumptions of Lemma~\ref{lem:curved-auxiliary}, put
\[
 \chi(t)=(1-t/T)^3\quad(0\le t\le T),\qquad
 H_\beta=\int_0^T t^3\chi(t)e^{-\beta t}q_t\dd t
 \quad(\beta\ge\alpha).
\]
For every bounded $1$-Lipschitz $f$,
\begin{equation}\label{eq:Hbounds}
 \norm{H_\beta}_1\le C\beta^{-4},\qquad
 |\ip{H_\beta}{Lf}|\le C\beta^{-7/2}.
\end{equation}
Furthermore, if $U_\alpha=\int_0^T\chi(t)e^{-\alpha t}q_t\dd t$, then
\begin{equation}\label{eq:curved-Ubounds}
 |\ip{U_\alpha}{Lf}|\le C\alpha^{-1/2},\qquad
 \ip{U_\alpha}{\dist(o,\cdot)}\le C\alpha^{-3/2}.
\end{equation}
\end{lemma}
\begin{proof}
$\norm{H_\beta}_1\le C\beta^{-4}$ follows from \eqref{eq:qflow}. For $x\sim y$,
since $e^{-\delta}q_t(x)\le q_t(y)\le e^\delta q_t(x)$, we have
\[
 \begin{gathered}
 |\log H_\beta(y)-\log H_\beta(x)|\le\delta.
 \end{gathered}
\]

Integration by parts in $\ell^1(m)$ gives
\begin{equation}\label{eq:curved-Hsource}
 (\beta-L)H_\beta
 =\int_0^T t^3\chi(t)e^{-\beta t}
       \left(\frac3t-\frac3{T-t}+W_t\right)q_t\dd t.
\end{equation}
 Writing
$a_\beta=(\beta-L)H_\beta/H_\beta$, Jensen's inequality for $|z|^3$
and then Tonelli's theorem imply
\[
 \begin{aligned}
 &\sum_xm(x)H_\beta(x)|a_\beta(x)|^3\\
 &\le C_1\int_0^T e^{-(\beta-\alpha/8)t}
       t^3\chi(t)\bigl(t^{-3}+(T-t)^{-3}+\alpha^3\bigr)\dd t\\
 &\le C_2\int_0^T e^{-(\beta-\alpha/8)t}(1+\alpha^3t^3)\dd t
 \le C_3\beta^{-1}.
 \end{aligned}
\]
Both endpoint singularities are integrable because the cutoff vanishes
to order three at each relevant endpoint. Also
${K_-}^3\norm{H_\beta}_1\le C\beta^{-1}$, since
${K_-}\le\alpha\le\beta$. Proposition~\ref{prop:moment}
therefore yields the bound in \eqref{eq:Hbounds}.

For $t>0$, $\frac12\int_\alpha^\infty(\beta-\alpha)^2
 t^3e^{-\beta t}\dd\beta=e^{-\alpha t}$. Hence
\begin{equation}\label{eq:curved-positive-mixture}
 U_\alpha=\frac12\int_\alpha^\infty
                   (\beta-\alpha)^2H_\beta\dd\beta.
\end{equation}
The equality holds in $\ell^1(m)$ by positivity and $\norm{H_\beta}_1\leq C\beta^{-4}$.
For the $|\ip{U_\alpha}{Lf}|$, Fubini follows from $|Lf|\le1$ and
\eqref{eq:Hbounds}. Thus
\[
 |\ip{U_\alpha}{Lf}|
 \le C_4\int_\alpha^\infty(\beta-\alpha)^2\beta^{-7/2}\dd\beta
 \le C_5\alpha^{-1/2}.
\]

A second integration by parts gives
\[
 (\alpha-L)U_\alpha
 =\phi_\alpha+\int_0^T e^{-\alpha t}
                         (\chi'(t)+\chi(t)W_t)q_t\dd t.
\]
Pair this identity with $\dist(o,\cdot)\wedge N\ge0$.
Since $\chi'\le0$ and $W_t\le\alpha/8$,
\[
 \begin{aligned}
 &\frac{7\alpha}{8}\ip{U_\alpha}{\dist(o,\cdot)\wedge N}\\
 \le&\ip{\phi_\alpha}{\dist(o,\cdot)\wedge N}+|\ip{U_\alpha}{L(\dist(o,\cdot)\wedge N)}|\le C\alpha^{-1/2}.
 \end{aligned}
\]
Letting $N\to\infty$ and using monotone convergence proves
\eqref{eq:curved-Ubounds}.
\end{proof}

For $\alpha>0$ define the conservative reversible Markov kernel
\begin{equation}\label{eq:curved-Qdefinition}
 Q_\alpha=4\alpha\int_{1/(4\alpha)}^{1/(2\alpha)}P_t\dd t.
\end{equation}
\begin{proposition}
\label{prop:curved-averaged-heat}
There is $M=M(\eps,\kappa)\ge1$ such that, with $A_*=4C_{\mathrm{LH}}$,
for every positive $\alpha\ge{K_-}$,
\begin{equation}\label{eq:curved-Qbounds}
 \begin{split}
 Q_\alpha[\dist(o,\cdot)](o)&\le M\alpha^{-1/2},\\
 \D(Q_\alpha(x,\cdot)\|Q_\alpha(y,\cdot))
   &\le A_*\alpha\dist(x,y)^2.
 \end{split}
\end{equation}
\end{proposition}
\begin{proof}
First take ${K_-}\le\alpha\le\alpha_0$, use the flow above, and write
$r_N=\dist(o,\cdot)\wedge N$,
$m_0=\ip{\phi_\alpha}{\dist(o,\cdot)}$.
For $0\le t\le T=\alpha^{-1}$, \eqref{eq:gradient-K} implies
$\Gamma P_tr_N\le e^2/2$, so
\[
 |P_tr_N(x)-P_tr_N(o)|\le e\eps^{-1/2}\dist(o,x).
\]
Integrate against $m\phi_\alpha$, and use reversibility and
$q_t\ge P_t\phi_\alpha$, to obtain
\[
 P_tr_N(o)\le\ip{q_t}{r_N}+e\eps^{-1/2}m_0.
\]
Multiplication by $\alpha\chi(t)e^{-\alpha t}$ and integration,
followed by \eqref{eq:curved-initial-moment} and
\eqref{eq:curved-Ubounds}, gives
\[
 \alpha\int_0^T\chi(t)e^{-\alpha t}P_tr_N(o)\dd t
 \le C\alpha^{-1/2}.
\]
On $[T/4,T/2]$, $\chi(t)e^{-\alpha t}\ge e^{-1/2}/8$.
Letting $N\to\infty$ therefore proves the first bound in
\eqref{eq:curved-Qbounds} for small $\alpha$.

For $\alpha\ge\alpha_0$, uniformisation gives
$P_t[\dist(o,\cdot)](o)\le t$ and hence
\[
 Q_\alpha[\dist(o,\cdot)](o)\le\frac{3}{8\alpha}
 \le\frac{3}{8\sqrt{\alpha_0}}\alpha^{-1/2}.
\]
Increasing $M$ covers these parameters as well.

For the entropy bound, since $K\ge-K_-$, we have $ I_K(t)\geq I_{-K_-}(t)$.
Convexity of relative entropy and Theorem~\ref{thm:logharnack} give
\[
 \begin{aligned}
 \D(Q_\alpha(x,\cdot)\|Q_\alpha(y,\cdot))
 &\le4\alpha\int_{1/(4\alpha)}^{1/(2\alpha)}
          \D(p_t(x,\cdot)\|p_t(y,\cdot))\dd t\\
 &\le4\alpha C_{\mathrm{LH}}\dist(x,y)^2
          \int_{1/(4\alpha)}^{1/(2\alpha)}\frac{\dd t}{I_{-K_-}(t)}\\
 &=4\alpha C_{\mathrm{LH}}\log\bigl(1+e^{{K_-}/(2\alpha)}\bigr)
                      \dist(x,y)^2\\
 &\le4C_{\mathrm{LH}}\alpha\dist(x,y)^2.
 \end{aligned}
\]
The integral identity also holds at ${K_-}=0$, when its value is
$\log2$. The last inequality uses ${K_-}\le\alpha$ and
$\log(1+\sqrt e)<1$. Thus one may take $A_*=4C_{\mathrm{LH}}$ in
\eqref{eq:curved-Qbounds}.
\end{proof}

\subsection{Volume doubling}
\label{subsec:lower-curvature-doubling}

\begin{proof}[Proof of Theorem~\ref{thm:main}]
Under $\CD(K,\infty)$, we first prove local doubling and then pass to all radii.

\emph{Local doubling.}
Let $A=4M^2$, with $M$ from
Proposition~\ref{prop:curved-averaged-heat}. If ${K_-} R^2\le1$,
set $\alpha=A/R^2\ge{K_-}$. By Markov's inequality and
\eqref{eq:curved-Qbounds},
\[
 Q_\alpha(o,B(o,R))\ge\frac12.
\]
For $y\in B(o,2R)$ the relative entropy in
\eqref{eq:curved-Qbounds} is at most $4A_*A=16C_{\mathrm{LH}} A$. Applying
\eqref{eq:binary-data-processing} to the event $B(o,R)$ gives
\[
 Q_\alpha(y,B(o,R))\ge c_*:=\tfrac14e^{-32C_{\mathrm{LH}} A}.
\]
Since $Q_\alpha$ preserves $m$,
\[
 c_*\vol(B(o,2R))
 \le\sum_y m(y)Q_\alpha(y,B(o,R))=\vol(B(o,R)).
\]
Choose $C_V=c_*^{-1}=4e^{32C_{\mathrm{LH}} A}$.
This argument is valid for every $o\in V$ and $r>0$, and gives
\begin{equation}\label{eq:curved-local-doubling}
 \vol(B(o,2r))\le C_V\vol(B(o,r))
 \qquad\text{if }K_-r^2\le1.
\end{equation}

\emph{Passage to arbitrary radii.}
Suppose $0<{K_-}\le1/16$ and let $a=(4\sqrt{K_-})^{-1}\ge1$.
For $s\ge0$, take a maximal set of points $x_i\in B(o,s)$ whose
pairwise distances are strictly greater than $a$.
The balls $B(x_i,a/2)$ are disjoint and contained in $B(o,s+a/2)$.
They satisfy
\[
 \vol(B(x_i,4a))\le C_V^3\vol(B(x_i,a/2)),
\]
by three applications of local doubling. These applications
are valid since their largest base radius is $2a$ and
${K_-}(2a)^2=1/4$.

Every vertex of $B(o,s+a)$ can be joined to a vertex of $B(o,s)$
by a geodesic of length at most $a+1$, and maximality places
that vertex within distance $a$ of some $x_i$. Since
$2a+1\le4a$, the balls $B(x_i,4a)$ cover $B(o,s+a)$.
Consequently
\begin{equation}\label{eq:curved-volume-step}
 \vol(B(o,s+a))\le C_V^3\vol(B(o,s+a/2)).
\end{equation}

For $R\ge a/2$, iterate this inequality in increments $a/2$,
starting at radius $R$. At most $\lceil2R/a\rceil$ steps give
\[
 \frac{\vol(B(o,2R))}{\vol(B(o,R))}
 \le C_V^{3\lceil2R/a\rceil}
 \le C_V^3\exp\bigl(24(\log C_V)\sqrt{K_-}\,R\bigr).
\]
For $R<a/2$, use local doubling directly.

Finally, when ${K_-}>1/16$, put $d_0=\lfloor\eps^{-1}\rfloor$.
A geodesic covering and the degree bound give
\[
 \frac{\vol(B(o,2R))}{\vol(B(o,R))}
 \le\kappa(d_0+1)^{\lceil R\rceil+1}
 \le\kappa(d_0+1)^2e^{R\log(d_0+1)}.
\]
Indeed, $B(o,2R)$ is covered by the balls of radius
$\lceil R\rceil$ centred at vertices of $B(o,R)$; the measure
comparison bounds the volume of each such ball by its cardinality
bound times $\kappa m_-$. Since $\sqrt{K_-}>1/4$, this is bounded by
$C_0e^{C_1\sqrt{K_-} R}$, with constants depending only on $\eps,\kappa$.
Together with the preceding cases this proves \eqref{eq:doubling-K}.
\end{proof}

\subsection{Applications}
We prove two consequences of volume doubling and discuss related comparison questions.

For a ball, define the first Dirichlet eigenvalue of $-L$ by
\[
 \lambda_1^D(B(o,R))
 :=\inf_{\substack{0\ne f\in\ell^2(V,m)\\
                   \operatorname{supp}f\subseteq B(o,R)}}
       \frac{\cE(f,f)}{\norm f_2^2},
\]
where $\cE$ is the energy of $L$.

\begin{theorem}\label{thm:spectral}
Under \eqref{eq:geometry} and $\CD(K,\infty)$, put
$K_-:=\max\{-K,0\}$. For every $o\in V$ and $R\ge1$,
\[
 \lambda_1^D(B(o,R))\le4C\max\{R^{-2},K_-\},
\]
where $C=C_0e^{C_1/2}$ and $C_0,C_1$ are the constants in
Theorem~\ref{thm:main}.
\end{theorem}

\begin{proof}[Proof of Theorem~\ref{thm:spectral}]
Put $\rho=R$ if $K_-=0$, and
$\rho=\min\{R,K_-^{-1/2}\}$ otherwise. Then
$\rho\le R$, $K_-\rho^2\le1$, and
$\rho^{-2}=\max\{R^{-2},K_-\}$.
Theorem~\ref{thm:main}, applied at radius $\rho/2>0$, gives
\[
 \frac{\vol(B(o,\rho))}{\vol(B(o,\rho/2))}
 \le C_0e^{C_1\sqrt{K_-}\,\rho/2}
 \le C_0e^{C_1/2}=C.
\]

Use $\chi(x)=(1-\dist(o,x)/\rho)_+$. The cutoff $\chi$ is supported
in $B(o,\rho)\subseteq B(o,R)$ and satisfies
$\chi\ge1/2$ on $B(o,\rho/2)$. Moreover, $\chi(x)\ne\chi(y)$
only if at least one of $x$ and $y$ lies in $B(o,\rho)$. Consequently,
\[
 \cE(\chi,\chi)\le\rho^{-2}\vol(B(o,\rho)),\qquad
 \norm{\chi}_2^2\ge\frac14\vol(B(o,\rho/2)).
\]
The Rayleigh principle therefore gives
\[
 \begin{aligned}
 \lambda_1^D(B(o,R))
 &\le\frac4{\rho^2}\frac{\vol(B(o,\rho))}{\vol(B(o,\rho/2))}
 \le4C\rho^{-2}=4C\max\{R^{-2},K_-\}.
 \end{aligned}
\]
\end{proof}

Conversely, under \eqref{eq:geometry} and $\CD(K,\infty)$,
a Dirichlet eigenvalue bound of the form given in
Theorem~\ref{thm:spectral} implies the corresponding volume doubling estimate.

\begin{proposition}
\label{prop:spectral-converse}
Assume \eqref{eq:geometry} and $\CD(K,\infty)$ for $K\in\R$, and put
$K_-:=\max\{-K,0\}$. Suppose that, for some constant $A\ge1$,
\begin{equation}\label{eq:spectral-converse-assumption}
 \lambda_1^D(B(o,R))\le A\max\{R^{-2},K_-\}
 \qquad(o\in V,\ R\ge1).
\end{equation}
Then there are constants $C_0\ge1$ and $C_1\ge0$, depending only on
$\eps,\kappa,A$, such that
\[
 \vol(B(o,2R))\le C_0e^{C_1\sqrt{K_-}\,R}\vol(B(o,R))
 \qquad(o\in V,\ R>0).
\]
In particular, if $A$ depends only on $\eps,\kappa$, this gives the
conclusion of Theorem~\ref{thm:main}.
\end{proposition}

\begin{proof}
First suppose that $R\ge1$ and $K_-R^2\le1$, and set $D=B(o,R)$.
Since $D$ is finite, its first Dirichlet eigenvalue
$\lambda=\lambda_1^D(D)$ has a nonnegative eigenfunction $\phi$,
extended by zero outside $D$ and normalized so that
$\max\phi=\phi(x_*)=1$ for some $x_*\in D$.
The zero extension satisfies $(L+\lambda)\phi\ge0$ on $V$.
Positivity of $P_t$ therefore gives
\[
 \frac{\dd}{\dd t}\bigl(e^{\lambda t}P_t\phi\bigr)
 =e^{\lambda t}P_t(L+\lambda)\phi\ge0.
\]
Taking $t=R^2/(2A)$ and using
\eqref{eq:spectral-converse-assumption}, we obtain
\[
 P_t\1_D(x_*)\ge P_t\phi(x_*)\ge e^{-\lambda t}
 \ge e^{-1/2}>\frac12.
\]

For every $y\in B(o,2R)$, we have $\dist(x_*,y)\le3R$ and
\[
 I_K(t)\ge te^{-2K_-t}\ge e^{-1}t.
\]
Thus Theorem~\ref{thm:logharnack} yields
\[
 \D(p_t(x_*,\cdot)\|p_t(y,\cdot))
 \le\frac{9C_{\mathrm{LH}}R^2}{I_K(t)}
 \le18eAC_{\mathrm{LH}}.
\]
Applying \eqref{eq:binary-data-processing} to the event $D$ gives
\[
 P_t\1_D(y)\ge c:=\tfrac14e^{-36eAC_{\mathrm{LH}}}
 \qquad(y\in B(o,2R)).
\]
Since $P_t$ preserves $m$, it follows that
\[
 c\vol(B(o,2R))
 \le\sum_{y\in V}m(y)P_t\1_D(y)=\vol(D).
\]

For $0<R<1$, the degree bound and \eqref{eq:geometry} give
\[
 \vol(B(o,2R))
 \le\kappa\bigl(1+\lfloor\eps^{-1}\rfloor\bigr)\vol(B(o,R)).
\]
Consequently, with
\[
 C_V=\max\bigl\{c^{-1},
             \kappa(1+\lfloor\eps^{-1}\rfloor)\bigr\},
\]
we have local doubling whenever $K_-R^2\le1$.

The geometric covering and iteration in the ``Passage to arbitrary
radii'' part of the proof of Theorem~\ref{thm:main} use only this
local doubling estimate and \eqref{eq:geometry}. Applying that
argument with the present $C_V$ proves the asserted bound at every
radius, with constants depending only on $\eps,\kappa,A$.
\end{proof}

For finite graphs, nonnegative curvature also restricts the spectral
gap. Write $\lambda_{\mathrm{gap}}(G)$ for the first nonzero eigenvalue
of $I-P_G$, where $P_G$ is the simple random-walk kernel of a finite
connected graph $G$.
Salez~\cite{Salez2022Expanders} proved that
there is no expander family with uniformly bounded degrees under
either nonnegative Ollivier--Ricci curvature or the Bakry--\'Emery
condition $\CD(0,\infty)$. Thus, for any such sequence $(G_n)$ with
$|V(G_n)|\to\infty$, one has $\lambda_{\mathrm{gap}}(G_n)\to0$.

Under nonnegative Ollivier--Ricci curvature,
M\"unch and Salez~\cite[Theorem~1]{MunchSalez2023Expansion} obtained
a quantitative conductance bound. Combining their estimate with
Cheeger's inequality gives
\[
 \lambda_{\mathrm{gap}}(G)
 \le C\sqrt{\frac{d\log d}{\log N}},
 \qquad N=|V(G)|\ge3,
\]
where $d$ is the maximum degree and $C$ is universal.

\begin{remark}
      \label{rem:eigenvalue}
      Cheng's theorem \cite[Theorem~1.1]{Cheng1975}, in the case of
nonnegative Ricci curvature, states that if $(M^n,g)$, $n\ge2$, is a
complete Riemannian manifold with $\operatorname{Ric}_g\ge0$, then for
every $o\in M$ and $R>0$,
\begin{equation}\label{eq:cheng-model-comparison}
 \lambda_1^D(B_M(o,R))
 \le\lambda_1^D(B_{\R^n}(0,R))
 =\lambda_1^D(B_{\R^n}(0,1))R^{-2}.
\end{equation}
Here the balls are open geodesic balls of the same radius.
Moreover, equality in the inequality in
\eqref{eq:cheng-model-comparison} holds if and only if $B_M(o,R)$ is
isometric to the Euclidean ball
$B_{\R^n}(0,R)$.

This motivates us asking whether the nonnegative-curvature case of
Theorem~\ref{thm:spectral} can be strengthened to a lattice comparison.
For this question, let $G$ be a simple unweighted graph with non-normalized Laplacian
\[
 \Delta_G f(x)=\sum_{y\sim x}(f(y)-f(x)).
\]
For a graph $H$, write
$\lambda_{1,H}^D(\Omega)$ for the first Dirichlet eigenvalue of
$-\Delta_H$ on a finite vertex set $\Omega$, using counting measure.
\begin{question}[A Cheng-type comparison with $\mathbb Z^n$]
\label{ques:cheng-lattice}
Let $G$ be an infinite connected simple unweighted graph satisfying
$\CD_{\Delta_G}(0,\infty)$, and set
\[
 d=\sup_x\deg(x)<\infty,\qquad n=\lceil d/2\rceil.
\]
Is it true that, for every $o\in V(G)$ and integer $R\ge1$,
\begin{equation}\label{eq:cheng-lattice-question}
 \lambda_{1,G}^D(B_G(o,R))
 \le\lambda_{1,\mathbb Z^n}^D(B_{\mathbb Z^n}(0,R))?
\end{equation}
Moreover, is $G\cong\mathbb Z^n$ equivalent to the existence of
$o\in V(G)$ such that equality in \eqref{eq:cheng-lattice-question}
holds for every integer $R\ge1$?
\end{question}

The choice of a lattice dimension $\lceil d/2\rceil$ is motivated by volume
growth. By Theorem~\ref{thm:sharp-volume},
$\vol(B_G(o,R))\le C_dR^{\lfloor d/2\rfloor}\le C_dR^n$ for $R\ge1$,
whereas balls in $\mathbb Z^n$ have volume comparable to $R^n$.

In the present nonnegative-curvature setting, applying
Theorem~\ref{thm:spectral} to $L=\Delta_G/d$ gives only
$\lambda_{1,G}^D(B_G(o,R))\le4dCR^{-2}$; it does not establish
\eqref{eq:cheng-lattice-question}.

At radius one, take $f(o)=\sqrt d$,
$f=1$ on the neighbours of $o$, and $f=0$ elsewhere. Writing
$q=\deg(o)\le d$, the Rayleigh principle gives
\[
 \begin{aligned}
 \lambda_{1,G}^D(B_G(o,1))
 &\le\frac{2q(d-\sqrt d)}{d+q}\le d-\sqrt d\\
 &\le2n-\sqrt{2n}
 =\lambda_{1,\mathbb Z^n}^D(B_{\mathbb Z^n}(0,1)),
 \end{aligned}
\]
since $d\le2n$ and $t\mapsto t-\sqrt t$ is increasing for $t\ge1$.
\end{remark}

\begin{remark}[Failure in dimension $\lfloor d/2\rfloor$]
\label{rem:cheng-degree-counterexample}
The comparison with $\mathbb Z^{\lfloor d/2\rfloor}$ can fail even
under $\CD(0,\infty)$. For every $k\ge1$, consider the unweighted
Cartesian product
\[
 G=\mathbb Z^k\mathbin{\square}K_2.
\]
Figure~\ref{fig:lattice-product-comparison} shows the case $k=2$.

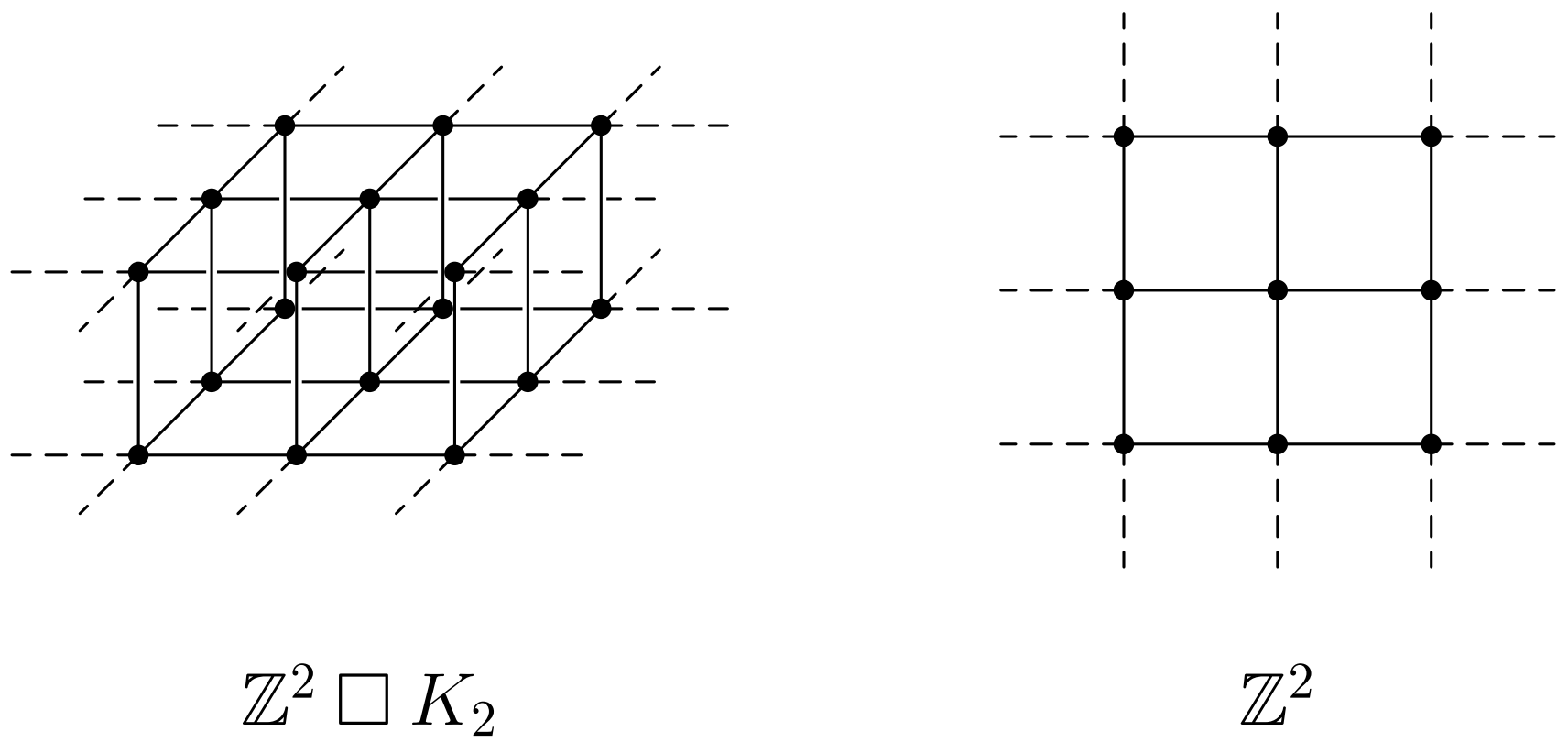
\begin{figure}[htbp]
 \centering
 \resizebox{0.75\textwidth}{!}{%
   % TikZ figure fragment, included by main.tex.
% Dashed rays indicate continuation of the infinite lattice.
\begingroup

\tikzset{
  lattice edge/.style={draw=black,line width=0.55pt,line cap=round},
  continuation/.style={lattice edge,dash pattern=on 4pt off 3pt},
  layer edge/.style={
    lattice edge,
    preaction={draw=white,line width=2.1pt,line cap=round}
  }
}
\newcommand{\VertexRadius}{2pt}
\newcommand{\LatticeExtent}{1.8}

\begin{tikzpicture}

% Left: two copies of Z^2, with exactly one edge between paired vertices.
% The three basis vectors control the axonometric projection.
\begin{scope}[shift={(0cm,0.25cm)},x={(1.08cm,0cm)},y={(0.50cm,0.50cm)},z={(0cm,1.25cm)}]
  \foreach \layer in {0,1} {
    \foreach \j in {-1,0,1} {
      \draw[lattice edge] (-1,\j,\layer) -- (1,\j,\layer);
      \draw[continuation] (-1,\j,\layer) -- (-\LatticeExtent,\j,\layer);
      \draw[continuation] (1,\j,\layer) -- (\LatticeExtent,\j,\layer);
    }
    \foreach \i in {-1,0,1} {
      \draw[lattice edge] (\i,-1,\layer) -- (\i,1,\layer);
      \draw[continuation] (\i,-1,\layer) -- (\i,-\LatticeExtent,\layer);
      \draw[continuation] (\i,1,\layer) -- (\i,\LatticeExtent,\layer);
    }
  }

  % White under-strokes distinguish crossings from vertices.
  % No dashed continuation in this direction: the factor K_2 is finite.
  \foreach \i in {-1,0,1} {
    \foreach \j in {-1,0,1} {
      \draw[layer edge] (\i,\j,0) -- (\i,\j,1);
    }
  }
  \foreach \layer in {0,1} {
    \foreach \i in {-1,0,1} {
      \foreach \j in {-1,0,1} {
        \fill (\i,\j,\layer) circle[radius=\VertexRadius];
      }
    }
  }
\end{scope}

% Right: the ordinary square lattice, centered at the same height.
\begin{scope}[shift={(6.2cm,0.875cm)},x={(1.05cm,0cm)},y={(0cm,1.05cm)}]
  \foreach \j in {-1,0,1} {
    \draw[lattice edge] (-1,\j) -- (1,\j);
    \draw[continuation] (-1,\j) -- (-\LatticeExtent,\j);
    \draw[continuation] (1,\j) -- (\LatticeExtent,\j);
  }
  \foreach \i in {-1,0,1} {
    \draw[lattice edge] (\i,-1) -- (\i,1);
    \draw[continuation] (\i,-1) -- (\i,-\LatticeExtent);
    \draw[continuation] (\i,1) -- (\i,\LatticeExtent);
  }
  \foreach \i in {-1,0,1} {
    \foreach \j in {-1,0,1} {
      \fill (\i,\j) circle[radius=\VertexRadius];
    }
  }
\end{scope}

\node[font=\Large,anchor=north] at (0,-1.55)
  {$\mathbb Z^2\mathbin{\square}K_2$};
\node[font=\Large,anchor=north] at (6.2,-1.55)
  {$\mathbb Z^2$};

\end{tikzpicture}
\endgroup%
 }
 \caption{$\mathbb Z^2\mathbin{\square}K_2$ and $\mathbb Z^2$}
 \label{fig:lattice-product-comparison}
\end{figure}

It is $(2k+1)$-regular. Since both $\mathbb Z^k$ and $K_2$ satisfy
$\CD(0,\infty)$ for their non-normalized Laplacians, the Cartesian
product rule of Cushing et~al.~\cite[Theorem~1.10]{CushingKamtueLiuPeyerimhoff2022} gives
$\CD(0,\infty)$ for $G$.

Both $G$ and $\mathbb Z^k$ are triangle-free. In a triangle-free
$d$-regular graph, the radius-one ball induces $K_{1,d}$, and its
Dirichlet matrix is $dI-A_{K_{1,d}}$. Since the largest adjacency
eigenvalue of this star is $\sqrt d$, we obtain
\[
 \begin{aligned}
 \lambda_{1,G}^D(B_G(o,1))
 &=2k+1-\sqrt{2k+1}\\
 &>2k-\sqrt{2k}
 =\lambda_{1,\mathbb Z^k}^D(B_{\mathbb Z^k}(0,1)).
 \end{aligned}
\]
The strict inequality follows from
$\sqrt{2k+1}-\sqrt{2k}<1$. Thus even the nonnegative-curvature graph
$\mathbb Z^k\mathbin{\square}K_2$, which has volume growth of order
$R^k$ and $\lfloor d/2\rfloor=k$, violates the comparison with
$\mathbb Z^k$. For $k=1$ this is the ladder graph, with
$3-\sqrt3>2-\sqrt2$.
This does not contradict Question~\ref{ques:cheng-lattice}, whose
comparison dimension here is $n=\lceil d/2\rceil=k+1$.
\end{remark}

\Needspace{10\baselineskip}
The following relative volume comparison follows by iterating
Theorem~\ref{thm:main}.

\begin{theorem}[Relative volume comparison]
\label{thm:relative-volume}
Assume \eqref{eq:geometry} and $\CD(K,\infty)$ for $K\in\R$.
Let $C_0,C_1$ be the constants in Theorem~\ref{thm:main}, and put
$K_-:=\max\{-K,0\}$ and $Q=\log_2C_0$.
For every $o\in V$ and $0<r\le R$,
\begin{equation}\label{eq:relative-volume-comparison}
 \frac{\vol(B(o,R))}{\vol(B(o,r))}
 \le C_0\left(\frac Rr\right)^Q
       e^{2C_1\sqrt{K_-}\,R}.
\end{equation}
All constants depend only on $\eps,\kappa$.
\end{theorem}

\begin{proof}[Proof of Theorem~\ref{thm:relative-volume}]
Set $k=\lceil\log_2(R/r)\rceil$. Since $R\le2^kr$,
monotonicity of ball volumes and $k$ applications of Theorem~\ref{thm:main}
give
\[
 \begin{aligned}
 \frac{\vol(B(o,R))}{\vol(B(o,r))}
 &\le\prod_{j=0}^{k-1}
       \left(C_0e^{C_1\sqrt{K_-}\,2^jr}\right)\\
 &=C_0^k e^{C_1\sqrt{K_-}\,(2^k-1)r}.
 \end{aligned}
\]
Here the product equals one when $k=0$. Since
$k<1+\log_2(R/r)$ and $2^kr<2R$, we have
\[
 C_0^k\le C_0\left(\frac Rr\right)^Q,
 \qquad e^{C_1\sqrt{K_-}\,(2^k-1)r}
       \le e^{2C_1\sqrt{K_-}\,R}.
\]
These bounds prove
\eqref{eq:relative-volume-comparison}.
\end{proof}

\begin{remark}
\label{rem:relative-volume-exponent}
In the nonnegative-curvature case, even when volume comparison is
restricted to integer radii, the optimal uniform doubling constant need
not give the optimal comparison exponent.

Fix $\eps=2/5$ and $\kappa=2$, and let $D_{\mathrm{int}}$ be the
optimal uniform doubling constant for positive integer radii over
all chains satisfying \eqref{eq:geometry} and $\CD(0,\infty)$.
Every support degree is at most two. For any centre, put
$H(k)=|B(o,k)|$ and $a_j=H(j)-H(j-1)$. Then $H(0)=1$ and
$0\le a_{j+1}\le a_j$, so for integers $1\le r\le R$,
\[
 \frac{H(R)-1}{R}
 =\frac1R\sum_{j=1}^{R}a_j
 \le\frac1r\sum_{j=1}^{r}a_j
 =\frac{H(r)-1}{r}.
\]
Thus $H(R)\le1+(R/r)(H(r)-1)\le(R/r)H(r)$, and the measure
bounds give the uniform comparison
\begin{equation}\label{eq:integer-relative-volume}
 \frac{\vol(B(o,R))}{\vol(B(o,r))}
 \le2\frac{H(R)}{H(r)}\le2\frac Rr
 \qquad(1\le r\le R,\ r,R\in\N).
\end{equation}
In particular, $D_{\mathrm{int}}\le4$. The symmetric simple random
walk on $\mathbb Z$, with constant measure $m_-$, belongs to this
class; its linear volume growth shows that the exponent $1$ is optimal.

To show that $D_{\mathrm{int}}>2$, take the path $\{0,1,2,3,4\}$ and set
\[
 \theta=2^{1/4},\qquad p=\frac{\theta}{1+\theta},\qquad
 q=\frac1{1+\theta},\qquad m(i)=m_-\theta^i.
\]
Define the Markov kernel by
\[
 \begin{aligned}
 p(i,i+1)&=p\quad(0\le i<4),&
 p(i,i-1)&=q\quad(1\le i\le4),\\
 p(0,0)&=q,& p(4,4)&=p,
 \end{aligned}
\]
with all other entries zero. Reversibility follows from
$\theta q=p$; moreover, $q>2/5$ and $m_-\le m(i)\le2m_-$.

Completing squares in the definition of $\Gamma_2$ for $L=P-I$ gives
\[
 \begin{aligned}
 \Gamma_2f(0)&\ge\frac{3pq}{4}(f(1)-f(0))^2,\\
 \Gamma_2f(4)&\ge\frac{3pq}{4}(f(4)-f(3))^2,\\
 \Gamma_2f(i)&\ge\frac{pq}{2}
   (f(i-1)-2f(i)+f(i+1))^2\qquad(1\le i\le3).
 \end{aligned}
\]
Hence this chain satisfies $\CD(0,\infty)$ and belongs to the same
parameter class. Since $\theta^4=2$, its integer-radius doubling ratio is
\[
 \frac{\vol(B(0,4))}{\vol(B(0,2))}
 =\frac{1+\theta+\theta^2+\theta^3+\theta^4}{1+\theta+\theta^2}
 =2+\frac{(\theta-1)^2(\theta+1)}{1+\theta+\theta^2}>2.
\]
Thus $2<D_{\mathrm{int}}\le4$, and the dyadic argument gives the
exponent $\log_2D_{\mathrm{int}}>1$ even when only integer radii
are used. But \eqref{eq:integer-relative-volume} remains valid with
prefactor $D_{\mathrm{int}}$, since $D_{\mathrm{int}}>2$.
Consequently, with that same prefactor, the exponent can be lowered
to its optimal value $1$ throughout the parameter class.
\end{remark}

It is therefore natural to ask for the optimal uniform exponent in relative
volume comparison for the normalized Laplacian.
\begin{question}
\label{ques:degree-relative-volume-exponents}
Fix an integer $d\ge2$, and consider connected simple unweighted
graphs with $1\le\deg(x)\le d$. Consider its
normalized Laplacian. Use measure $m(x)=\deg(x)$ and write
$\vol(A)=\sum_{x\in A}m(x)$.
Assume $G$ is $\CD(0,\infty)$.
Can one determine the optimal uniform exponent $q(d)$? More precisely,
let $q(d)$ be the infimum of the exponents $q\ge0$ for which
\[
 \frac{\vol(B(o,R))}{\vol(B(o,r))}
 \le C(d)\left(\frac Rr\right)^q
 \qquad(1\le r\le R,\ r,R\in\N)
\]
holds with a finite prefactor depending only on $d$, uniformly over
all such graphs, centres, and integer radii.
\end{question}

\section{Poincar\'e and Parabolic Harnack inequalities}
\label{sec:harnack-bridge}

From this section onward, we assume \eqref{eq:geometry} and
$\CD(0,\infty)$ unless explicitly stated otherwise.
As applications of Theorem~\ref{thm:main}, we establish local
Poincar\'e and Harnack inequalities and deduce the positive-harmonic
Liouville property.

\begin{lemma}
\label{lem:second-cutoffs}
There are constants $K_\Gamma,K_L<\infty$, depending only on
$\eps,\kappa$, such that for every $o\in V$ and $R>0$ there is a
finitely supported function $0\le\chi_{o,R}\le1$ satisfying
\begin{equation}\label{eq:cutoff-support}
 \chi_{o,R}=1\ \text{on }B(o,R),\qquad
 \chi_{o,R}=0\ \text{on }\{x:\dist(o,x)\ge2R\},
\end{equation}
and
\begin{equation*}
 \Gamma\chi_{o,R}\le\frac{K_\Gamma}{R^2},\qquad
 |L\chi_{o,R}|\le\frac{K_L}{R^2}.
\end{equation*}
\end{lemma}
\begin{proof}
Let $M$ be as in Proposition~\ref{prop:curved-averaged-heat}, and put
\[
 f_R=\left(2-\frac{\dist(o,\cdot)}R\right)_+\wedge1,
 \qquad \alpha=\frac{16M^2}{R^2},\qquad v_R=Q_\alpha f_R.
\]
The function $f_R$ equals one on $B(o,R)$ and zero on
$\{x:\dist(o,x)\ge2R\}$, and is $1/R$-Lipschitz. Thus
\eqref{eq:curved-Qbounds} gives, for every $x\in V$,
\[
 |v_R(x)-f_R(x)|
 \le\frac1R Q_\alpha[\dist(x,\cdot)](x)
 \le\frac{M}{R\sqrt\alpha}=\frac14.
\]
Moreover, $\Gamma f_R\le1/(2R^2)$, so convexity of $\Gamma$ and
\eqref{eq:gradient-K} under $\CD(0,\infty)$ yield
\[
 \Gamma v_R\le Q_\alpha\Gamma f_R\le\frac1{2R^2}.
\]
By the heat equation and \eqref{eq:curved-Qdefinition},
\[
 \begin{aligned}
 Lv_R&=4\alpha\bigl(P_{1/(2\alpha)}f_R-P_{1/(4\alpha)}f_R\bigr),\\
 |Lv_R|&\le4\alpha=\frac{64M^2}{R^2},
 \end{aligned}
\]
where the last bound uses $0\le P_tf_R\le1$.

Choose a fixed smooth $F:\R\to[0,1]$ equal to zero on
$(-\infty,1/4]$ and one on $[3/4,\infty)$, and set
$\chi_{o,R}=F(v_R)$. The approximation bound implies
$v_R\ge3/4$ on $B(o,R)$ and $v_R\le1/4$ on
$\{x:\dist(o,x)\ge2R\}$. This proves \eqref{eq:cutoff-support};
the support is finite because $B(o,2R)$ is.
The Lipschitz bound and Taylor's formula on each edge give
\[
 \begin{aligned}
 \Gamma\chi_{o,R}
 &\le\|F'\|_\infty^2\Gamma v_R
       \le\frac{K_\Gamma}{R^2},\\
 |L\chi_{o,R}|
 &\le\|F'\|_\infty|Lv_R|
       +\|F''\|_\infty\Gamma v_R
       \le\frac{K_L}{R^2},
 \end{aligned}
\]
with constants depending only on $\eps,\kappa$.
\end{proof}

For a finite nonempty set $U\subseteq V$, write
\[
 \cE_U(f)=\sum_{\{x,y\}\subset U}m(x)p(x,y)(f(x)-f(y))^2,
 \qquad f_U=\frac{\sum_{x\in U}m(x)f(x)}{\vol(U)}.
\]

On $\ell^2(U,m)$, define the operators $P_U$ and $L_U$ by
\[
 \begin{aligned}
 (P_Uh)(x)&=\sum_{y\in U}p(x,y)h(y),\\
 (L_Uh)(x)&=\sum_{y\in U}p(x,y)(h(y)-h(x))
 \qquad(x\in U).
 \end{aligned}
\]
Let $I_U$ be the identity operator on this space. For $s\ge0$, set
\[
 D_s=e^{s(P_U-I_U)},\qquad M_s=e^{sL_U}.
\]
Both semigroups act on $\ell^2(U,m)$ and are positivity preserving;
the Dirichlet semigroup $D_s$ satisfies $D_s\1\le\1$, while
$M_s\1=\1$.

\begin{lemma}[Weak local Poincar\'e inequality]
\label{lem:weak-poincare}
Under \eqref{eq:geometry} and $\CD(0,\infty)$, set
\[
 \Lambda=\max\{2,\sqrt{128K_L/\eps}\},
\]
where $K_L$ is from Lemma~\ref{lem:second-cutoffs}.
Then, for every $o\in V$, $R\ge1$, and real function $f$ on
$U=B(o,\Lambda R)$,
\begin{equation*}
 \sum_{x\in B(o,R)}m(x)|f(x)-f_{B(o,R)}|^2
 \le\frac{64R^2}{\eps}\cE_U(f).
\end{equation*}
\end{lemma}
\begin{proof}
Put $B=B(o,R)$ and $t=4R^2/\eps$.
For adjacent vertices $u\sim v$ and any $A\subseteq V$, apply
\eqref{eq:variance-K} to $\1_A$. The edge lower bound gives
\[
 \begin{aligned}
 \eps t\,|p_t(u,A)-p_t(v,A)|^2
 &\le2t\Gamma(P_t\1_A)(u)\\
 &\le\Var_{p_t(u,\cdot)}\1_A\\
 &=p_t(u,A)(1-p_t(u,A))\le\tfrac14.
 \end{aligned}
\]
Taking the supremum over $A\subseteq V$, we obtain
\[
 \|p_t(u,\cdot)-p_t(v,\cdot)\|_{\mathrm{TV}}
 =\sup_{A\subseteq V}|p_t(u,A)-p_t(v,A)|
 \le\frac1{2\sqrt{\eps t}}.
\]
The triangle inequality along a geodesic then gives
\[
 \|p_t(x,\cdot)-p_t(y,\cdot)\|_{\mathrm{TV}}
 \le\frac{\dist(x,y)}{2\sqrt{\eps t}}\le\frac12
 \qquad(x,y\in B).
\]

Set $\chi=\chi_{o,\Lambda R/2}$. By Lemma~\ref{lem:second-cutoffs},
$0\le\chi\le1$, $\chi=1$ on $B$, $\chi=0$ on $U^c$, and
$\|L\chi\|_\infty\le\delta:=4K_L/(\Lambda^2R^2)$.
Since $(P_U-I_U)(\chi|_U)=(L\chi)|_U$, Duhamel's formula gives,
for $x\in B$,
\[
 \begin{aligned}
 1-D_t\1(x)
 &\le1-D_t(\chi|_U)(x)\\
 &=-\int_0^tD_s\bigl((L\chi)|_U\bigr)(x)\dd s\\
 &\le\delta\int_0^tD_s\1(x)\dd s
 \le\delta t\le\frac18.
 \end{aligned}
\]

Writing $r=\1-P_U\1\ge0$, we have
$L_U=P_U+\operatorname{diag}(r)-I_U$.
Since $P_U\le P_U+\operatorname{diag}(r)$ entrywise and
$(P_U^j)(x,z)\le P^j(x,z)$ for $x,z\in U$ and $j\ge0$,
the uniformisation expansions give, after extending
$D_t(x,\cdot)$ and $M_t(x,\cdot)$ by zero outside $U$,
\[
 D_t(x,z)\le\min\{M_t(x,z),p_t(x,z)\}
 \qquad(x\in U,\ z\in V).
\]
Both $M_t(x,\cdot)$ and $p_t(x,\cdot)$ have total mass one, so
\[
 \begin{aligned}
 \|M_t(x,\cdot)-p_t(x,\cdot)\|_{\mathrm{TV}}
 &=1-\sum_z\min\{M_t(x,z),p_t(x,z)\}\\
 &\le1-D_t\1(x)\le\frac18
 \qquad(x\in B).
 \end{aligned}
\]
For $x,y\in B$, the triangle inequality therefore shows that $\|M_t(x,\cdot)-M_t(y,\cdot)\|_{\mathrm{TV}}\leq 3/4$. This implies that
\[
 \sum_z\min\{M_t(x,z),M_t(y,z)\}\ge1/4.
\]

Set $A_x=\sum_{z\in U}M_t(x,z)|f(x)-f(z)|^2$. Integrating $|f(x)-f(y)|^2\le
2|f(x)-f(z)|^2+2|f(y)-f(z)|^2$ against this common mass yields
$|f(x)-f(y)|^2\le8(A_x+A_y)$. Hence
\[
 \begin{aligned}
 \sum_{x\in B}m(x)|f(x)-f_B|^2
 &=\frac{1}{2\vol(B)}\sum_{x,y\in B}m(x)m(y)|f(x)-f(y)|^2\\
 &\le8\sum_{x\in B}m(x)A_x
 \le16\langle f,(I_U-M_t)f\rangle_{\ell^2(U,m)}\\
 &\le16t\cE_U(f)=\frac{64R^2}{\eps}\cE_U(f).
 \end{aligned}
\]
The last two inequalities follow, respectively, from the reversibility
of $M_t$ and the bound $1-e^{-t\lambda_k}\le t\lambda_k$
for each eigenvalue $\lambda_k$ of $-L_U$.
\end{proof}

\begin{theorem}[Local Poincar\'e inequality]
\label{thm:local-poincare}
Under \eqref{eq:geometry} and $\CD(0,\infty)$, there is a constant
$C_{\mathrm P}=C_{\mathrm P}(\eps,\kappa)$ such that, for every
$o\in V$, $R\ge1$, and real function $f$ defined on $B(o,R)$,
\begin{equation}\label{eq:local-poincare}
 \sum_{x\in B(o,R)}m(x)|f(x)-f_{B(o,R)}|^2
 \le C_{\mathrm P}R^2\cE_{B(o,R)}(f).
\end{equation}
\end{theorem}
\begin{proof}
Let $\Lambda$ be as in Lemma~\ref{lem:weak-poincare}.
It suffices to consider integer radii $N\ge1$. Write $B=B(o,N)$, and set
$A=\lceil16\Lambda\rceil$ and $q=1+A^{-1}$.

Suppose first that $N<A$. For each $x\in B$, choose a geodesic
$o=x_0\sim x_1\sim\cdots\sim x_k=x$, where $k=\dist(o,x)\le N$. By Cauchy--Schwarz and
$m(x_{i-1})p(x_{i-1},x_i)\ge\eps m_-$ from \eqref{eq:geometry},
\begin{equation}
        \label{eq:pathestimate}
 |f(x)-f(o)|^2\le\frac{N}{\eps m_-}\sum_{i=1}^k
       m(x_{i-1})p(x_{i-1},x_i)|f(x_i)-f(x_{i-1})|^2\le\frac{N}{\eps m_-}\cE_B(f).
\end{equation}
It follows that
\[
 \sum_{x\in B}m(x)|f(x)-f_B|^2
 \le\sum_{x\in B}m(x)|f(x)-f(o)|^2
 \le\frac{N\vol(B)}{\eps m_-}\cE_B(f)
 \le C_1N^2\cE_B(f).
\]

Suppose now that $N\ge A$. Set $r_i=q^i$, and let $k$ be the
first index such that $r_k\ge2N/A$. For each $x\in B$, fix a
geodesic from $x$ to $o$, and let $c_i(x)$ be the vertex at distance
$\min\{\lfloor Ar_i\rfloor,\dist(x,o)\}$ from $x$ along this
geodesic. For $0\leq i\leq k$, define
\[
 B_i(x)=B(c_i(x),r_i),\quad
 D_i(x)=B(c_i(x),4\Lambda r_i),\quad a_i(x)=f_{B_i(x)}.
\]
These balls lie in $B$. Indeed, before the centre reaches $o$,
\[
 \dist(o,c_i(x))+4\Lambda r_i
 =\dist(o,x)-\lfloor Ar_i\rfloor+4\Lambda r_i\le N.
\]
When $c_i(x)=o$, $4\Lambda r_i\le4\Lambda r_k
<8\Lambda qN/A\le N$. Also $Ar_k\ge2N$, so $c_k(x)=o$
and $a_k(x)=a_*:=f_{B(o,r_k)}$ is independent of $x$.

For $i<k$, consecutive centres have distance at most
$A(r_{i+1}-r_i)+1\le2r_i$, and hence
$B_i(x)\cup B_{i+1}(x)\subset B(c_i(x),4r_i)$.
Set $W_i(x)=B(c_i(x),4r_i)$. For $j=i,i+1$, the inclusion
$B_j(x)\subseteq W_i(x)$ and Jensen's inequality give
\[
 \begin{aligned}
 |a_j(x)-f_{W_i(x)}|^2
 &=\left|\frac{1}{\vol(B_j(x))}\sum_{z\in B_j(x)}
      m(z)(f(z)-f_{W_i(x)})\right|^2\\
 &\le\frac{1}{\vol(B_j(x))}\sum_{z\in W_i(x)}
      m(z)|f(z)-f_{W_i(x)}|^2\\
 &\le\frac{1024r_i^2}{\eps\vol(B_j(x))}\cE_{D_i(x)}(f).
 \end{aligned}
\]
The last inequality is Lemma~\ref{lem:weak-poincare}, applied with
centre $c_i(x)$ and radius $4r_i$, whose enlarged ball is $D_i(x)$.
Moreover, $B_i(x)\subseteq B(c_{i+1}(x),3r_i)
\subseteq B(c_{i+1}(x),4r_{i+1})$. Thus Theorem~\ref{thm:main} gives
\[
 \vol(B_i(x))
 \le\vol(B(c_{i+1}(x),4r_{i+1}))
 \le C_V^2\vol(B_{i+1}(x)),
\]
where $C_V$ is the doubling constant. Combining these estimates yields
\[
 \begin{aligned}
 |a_i(x)-a_{i+1}(x)|^2
 &\le2|a_i(x)-f_{W_i(x)}|^2
       +2|a_{i+1}(x)-f_{W_i(x)}|^2\\
 &\le\frac{C_2r_i^2}{\vol(B_i(x))}\cE_{D_i(x)}(f).
 \end{aligned}
\]

Fix $i$ and $z\in V$. If $z\in D_i(x)$, one has
$\dist(x,z)\le(A+4\Lambda)r_i$ and
$\vol(B_i(x))\ge C_3^{-1}\vol(B(z,r_i))$ by doubling. Consequently,
\[
 \sum_{x\in B:\,z\in D_i(x)}
 \frac{m(x)}{\vol(B_i(x))}\le C_4.
\]
The above two estimates give
\[
 \|a_i-a_{i+1}\|_{\ell^2(B,m)}\le C_5r_i\cE_B(f)^{1/2}.
\]
Since $\sum_{i<k}r_i\le C_6N$, the $\ell^2(B,m)$ triangle
inequality therefore yields
\[
 \|a_0-a_*\|_{\ell^2(B,m)}
 \le\sum_{i<k}\|a_i-a_{i+1}\|_{\ell^2(B,m)}
 \le C_7N\cE_B(f)^{1/2}.
\]

To control the initial term, join $x$ to $c_0(x)$ along the chosen
geodesic and then to each vertex of $B_0(x)=B(c_0(x),1)$.
After erasing any backtracking, these paths lie in
$B\cap B(x,A+1)$ and have length at most $A+1$.
Arguing as in the derivation of \eqref{eq:pathestimate}, we obtain
\[
 |f(x)-a_0(x)|^2
 \le\frac{A+1}{\eps m_-}\cE_{B\cap B(x,A+1)}(f).
\]
The degree bound and $m(x)\le\kappa m_-$ therefore give
\[
 \|f-a_0\|_{\ell^2(B,m)}\le C_8\cE_B(f)^{1/2}.
\]
Combining these bounds with the minimizing property of $f_B$ gives
\[
 \begin{aligned}
 \|f-f_B\|_{\ell^2(B,m)}
 &\le\|f-a_*\|_{\ell^2(B,m)}\\
 &\le\|f-a_0\|_{\ell^2(B,m)}+\|a_0-a_*\|_{\ell^2(B,m)}\\
 &\le C_9N\cE_B(f)^{1/2}.
 \end{aligned}
\]
Squaring proves \eqref{eq:local-poincare} with
$C_{\mathrm P}=\max\{C_1,C_9^2\}$.
\end{proof}

\Needspace{7\baselineskip}
The volume doubling and local Poincar\'e inequalities established above
yield a parabolic Harnack inequality through the following
theorem of Delmotte.

\begin{theorem}[{\cite[Theorem~2.1]{Delmotte}}]
\label{thm:delmotte}
Suppose $V$ satisfies volume doubling with constant $C_V$.
Assume that, for every $o\in V$,
$r>0$, and real function $f$ on $B(o,2r)$, the following local
Poincar\'e inequality with dilation factor $2$ holds:
\[
 \begin{aligned}
 \sum_{x\in B(o,r)}m(x)|f(x)-f_{B(o,r)}|^2
 &\le C_{\mathrm P}r^2\cE_{B(o,2r)}(f).
 \end{aligned}
\]
Then there is $C_{\mathrm H}=C_{\mathrm H}(\eps,C_V,C_{\mathrm P})$ with the
following property. For every $o\in V$ and real $R\ge1$, let $u$
be nonnegative and continuous in time on
$[0,5R^2]\times B(o,2R+1)$ and satisfy
\[
 \partial_tu(t,x)=L u(t,\cdot)(x)
 \qquad(0<t<5R^2,\ x\in B(o,2R)).
\]
Then
\[
 \sup_{\substack{R^2\le t\le2R^2\\x\in B(o,R)}}u(t,x)
 \le C_{\mathrm H}
 \inf_{\substack{4R^2\le t\le5R^2\\x\in B(o,R)}}u(t,x).
\]
\end{theorem}

\begin{theorem}[Parabolic Harnack inequality]
\label{thm:parabolic-harnack}
Under \eqref{eq:geometry} and $\CD(0,\infty)$,
there is $C_{\mathrm{PH}}=C_{\mathrm{PH}}(\eps,\kappa)$ such that, for
every $o\in V$ and real $R\ge1$, every nonnegative function $u$
continuous in time on $[0,5R^2]\times B(o,2R+1)$ and satisfying
\[
 \partial_tu(t,x)=L u(t,\cdot)(x)
 \qquad(0<t<5R^2,\ x\in B(o,2R))
\]
obeys
\begin{equation}\label{eq:parabolic-harnack}
 \sup_{\substack{R^2\le t\le2R^2\\x\in B(o,R)}}u(t,x)
 \le C_{\mathrm{PH}}
 \inf_{\substack{4R^2\le t\le5R^2\\x\in B(o,R)}}u(t,x).
\end{equation}
\end{theorem}
\begin{proof}
Theorem~\ref{thm:main} gives the volume doubling hypothesis of
Theorem~\ref{thm:delmotte}. The Poincar\'e hypothesis follows from
Theorem~\ref{thm:local-poincare} by enlarging the energy ball from
$B(o,r)$ to $B(o,2r)$ when $r\ge1$. For $0<r<1$, $B(o,r)=\{o\}$,
so the left-hand side is zero.
Applying Theorem~\ref{thm:delmotte} gives \eqref{eq:parabolic-harnack}
with $C_{\mathrm{PH}}=C_{\mathrm H}(\eps,C_V,C_{\mathrm P})$, which
depends only on $\eps,\kappa$.
\end{proof}

\begin{corollary}
\label{cor:positive-liouville}
Under \eqref{eq:geometry} and $\CD(0,\infty)$, every nonnegative harmonic function is constant.
\end{corollary}
\begin{proof}
On finite connected graphs, this follows from the maximum principle.
Suppose $V$ is infinite. For a nonnegative function $h$ on
$B(o,2R+1)$ satisfying $Lh=0$
on $B(o,2R)$, apply \eqref{eq:parabolic-harnack} to the stationary
solution $u(t,x)=h(x)$. This gives the elliptic Harnack inequality
\begin{equation}\label{eq:elliptic-harnack}
 \sup_{B(o,R)}h\le C_{\mathrm{PH}}\inf_{B(o,R)}h
 \qquad(R\ge1).
\end{equation}

Finally, for a nonnegative global harmonic $h$, put $a=\inf_Vh$.
Applying \eqref{eq:elliptic-harnack} to $h-a$ gives
\[
 0\le h(o)-a\le C_{\mathrm{PH}}\inf_{B(o,R)}(h-a)\to0.
\]
Thus $h$ is constant.
\end{proof}
\begin{remark}[Liouville property]
In Riemannian geometry,
Yau~\cite[Section~2, Corollary~1]{Yau1975Harmonic} proved that
every positive harmonic function on a complete connected manifold with
nonnegative Ricci curvature is constant.
In the graph setting, Hua~\cite{Hua2019Liouville} established the
Liouville property for bounded harmonic functions on connected locally
finite weighted graphs satisfying $\CD(0,\infty)$, provided that the
Laplacian is bounded on $\ell^2(V,m)$.
\end{remark}

\section{Li--Yau gradient estimates}
\label{sec:li-yau}
\subsection{Li--Yau gradient estimates}
We recall Li-Yau gradient estimate. Let $(M,g)$ be a complete
$n$-dimensional Riemannian manifold without boundary, and let $u>0$
solve $\partial_tu=\Delta_Mu$ on $(0,\infty)\times M$.
If $\operatorname{Ric}_g\ge0$, Li and Yau~\cite{LiYau1986} proved
\begin{equation}\label{eq:ly-classical}
 \frac{|\nabla u|^2}{u^2}-\frac{\partial_tu}{u}
 =|\nabla\log u|^2-\partial_t\log u
 \le\frac n{2t}.
\end{equation}
The chain rule gives
$\Delta_M\log u=(\Delta_Mu)/u-|\nabla\log u|^2$, so
\eqref{eq:ly-classical} is equivalent to the Laplacian estimate
\begin{equation}\label{eq:ly-classical-laplacian}
 -\Delta_M\log u\le\frac n{2t}.
\end{equation}

Let $f$ be bounded, nonnegative, and nonzero, and put $u_t=P_tf$.
For $t>0$, we have
\[
 L(\log u_t)(x)=\partial_t\log u_t(x)-\Psi(\log u_t)(x).
\]
\Needspace{12\baselineskip}
On graphs, we prove the following strengthened form of
Theorem~\ref{thm:li-yau}:

\begin{theorem}[Li--Yau gradient estimate]
\label{thm:li-yau-full}
There are constants
$\alpha_0>1$ and $C_{LY}<\infty$, depending only on $\eps,\kappa$,
such that if $f$ is a bounded, nonnegative, nonzero function, then
\[
 \begin{aligned}
 \Gamma(\log P_tf)(x)-\alpha\,\partial_t\log P_tf(x)
 &\le\frac{C_{LY}\alpha}{t},\\[3pt]
 -L(\log P_tf)(x)-(\alpha-1)\,\partial_t\log P_tf(x)
 &\le\frac{C_{LY}\alpha}{t},
 \end{aligned}
\]
for $\alpha\geq \alpha_0,\ t>0$ and $x\in V$.
\end{theorem}

We first establish three auxiliary estimates.

\begin{lemma}
\label{lem:ly-time}
There exists a constant $C>0$, depending only on $\eps,\kappa$, such that
\[
 -\partial_t\log u_t(x)\le\frac{C}t
 \qquad(t>0).
\]
\end{lemma}
\begin{proof}
Fix $C\ge1$ such that
\eqref{eq:sharp-time-harnack} holds with $C_2=C$ when $\Lambda=2$. Since $\partial_t\log u_t=Pu_t/u_t-1\ge-1$, for every finite
$T>0$ the quantity
\[
 S=S(T):=\sup_{0<t\le T,\ z\in V}
               t[-\partial_t\log u_t(z)]_+\le T
\]
is finite. If $S=0$, the desired bound holds on $(0,T]$.

Suppose $S>0$, put $\delta=1/(4C)$, and choose $0<t\le T$ and
$x\in V$ such that
$-t\partial_t\log u_t(x)\ge(1-\delta)S$.
Set
\[
 s_0=\frac{t}{1+\delta},\qquad
 w_r=\frac{S}{s_0}u_r+\partial_ru_r.
\]
The definition of $S$ gives $w_{s_0}\ge0$ on $V$.
Moreover, $w_r=(\frac{S}{s_0}I+L)P_rf$ is a bounded heat solution, so
$w_r=P_{r-s_0}w_{s_0}\ge0$ for $r\ge s_0$.
At the time $t$,
\[
 0\le w_t(x)\le\frac{2\delta S}{t}u_t(x).
\]

Put $h=(t-s_0)/2$. For $r\in[t-h,t]$ we have
$r-s_0\in[h,2h]$; hence \eqref{eq:sharp-time-harnack}, applied to
$w_{s_0+\cdot}$, yields
\[
 w_r(x)\le C w_t(x)\le\frac{2C\delta S}{t}u_t(x).
\]
Applying \eqref{eq:duhamel} with generator $-SI/s_0$ to
$\partial_ru_r=-(S/s_0)u_r+w_r$ and solving for $u_{t-h}(x)$ yields
\[
 \begin{aligned}
 u_{t-h}(x)
 &=e^{Sh/s_0}u_t(x)-\int_{t-h}^t e^{S(r-t+h)/s_0}w_r(x)\dd r\\
 &\ge\left[e^{Sh/s_0}-\frac{2C\delta}{1+\delta}
                 (e^{Sh/s_0}-1)\right]u_t(x)\\
 &\ge\tfrac12e^{Sh/s_0}u_t(x).
 \end{aligned}
\]
Here $2C\delta/(1+\delta)\le1/2$.
Since $h<t/2$, \eqref{eq:sharp-time-harnack} also gives
$u_{t-h}(x)\le Cu_t(x)$. Consequently,
\[
 \frac{\delta S}{2}=\frac{Sh}{s_0}\le\log(2C),\qquad
 S(T)\le8C\log(2C).
\]
The right-hand side depends only on $\eps,\kappa$ and is independent
of $T$. Letting $T$ vary proves the lemma.
\end{proof}

\begin{lemma}
\label{lem:ly-sqrt}
There are $\beta\ge1$ and $C<\infty$ such that
\[
 \frac{\Gamma(\sqrt{u_t})}{u_t}
 \le\beta\partial_t\log u_t+\frac{C}{t}
 \qquad(t>0).
\]
\end{lemma}
\begin{proof}
Write $\rho=\max\{1,t_{\mathrm{reg}}\}$, with $t_{\mathrm{reg}}$
from Lemma~\ref{lem:sharp-positive-regularity}.
Let $C_0$ be the constant in Lemma~\ref{lem:ly-time}, and fix
$A\ge\max\{2C_0,4(\rho+1)\}$.
Fix $t>0$ and set
\[
 w_r=\frac At u_r+\partial_ru_r\quad(r\ge t/2).
\]
Lemma~\ref{lem:ly-time} gives $w\ge0$.

The identity
\[
 \begin{aligned}
 \frac{\Gamma(\sqrt{u_t})(x)}{u_t(x)}
 &=\tfrac12\partial_t\log u_t(x)+
   \sum_{y\sim x}p(x,y)
       \left(1-\sqrt{\frac{u_t(y)}{u_t(x)}}\right)\\
 &\le\tfrac12\partial_t\log u_t(x)+1
 \end{aligned}
\]
holds for every $t>0$.

If $w_t(x)/u_t(x)>1/\rho$, then $A>0$ gives
\[
 \frac{\Gamma(\sqrt{u_t})(x)}{u_t(x)}
 \le\frac12\frac{w_t(x)}{u_t(x)}+1
 \le\left(\rho+\frac12\right)\frac{w_t(x)}{u_t(x)}.
\]

It remains to consider $0\le w_t(x)/u_t(x)\le1/\rho$.
Since $\partial_t\log u_t(x)\ge-1$, we have
\[
 \frac At=\frac{w_t(x)}{u_t(x)}-\partial_t\log u_t(x)
 \le1+\frac1\rho,
 \qquad
 t\ge\frac{A\rho}{\rho+1}\ge4\rho.
\]
For $\rho\le h\le t/4$, \eqref{eq:duhamel} with generator $-(A/t)I$ gives
\[
 u_t=e^{-Ah/t}u_{t-h}+\int_0^h e^{-Ar/t}w_{t-r}\dd r.
\]
The gradient estimate in Lemma~\ref{lem:sharp-positive-regularity},
applied to $u(t-2h+\cdot,\cdot)$ and $w(t-2r+\cdot,\cdot)$, gives
\[
 \sqrt{\Gamma u_{t-h}(x)}\le C_1 h^{-1/2}u_t(x),\qquad
 \sqrt{\Gamma w_{t-r}(x)}\le C_1 r^{-1/2}w_t(x)
 \quad(\rho\le r\le h).
\]

For the integral $U=\int_0^\rho e^{-Ar/t}w_{t-r}\dd r$,
apply \eqref{eq:sharp-time-harnack} with $\Lambda=2$ to
$w_{t/2+\cdot}$. Since $r\le\rho\le t/4$, this gives
$w_{t-r}(x)\le Cw_t(x)$ for $0\le r\le\rho$.
Integrating the heat equation at $x$ therefore yields
\[
 \sum_{y\sim x}p(x,y)\int_0^\rho w_{t-r}(y)\dd r
 \le(1+C\rho)w_t(x).
\]
Since $w\ge0$ and $p(x,y)\ge\eps$, this implies
$|U(y)-U(x)|\le Cw_t(x)$ for $y\sim x$, and hence
$\sqrt{\Gamma U(x)}\le C_2w_t(x)$.

Splitting the Duhamel integral at $\rho$ and applying Minkowski's
inequality gives
\[
 \begin{aligned}
 \sqrt{\Gamma u_t(x)}
 &\le e^{-Ah/t}\sqrt{\Gamma u_{t-h}(x)}
       +\sqrt{\Gamma U(x)}+\int_\rho^h e^{-Ar/t}\sqrt{\Gamma w_{t-r}(x)}\dd r\\
 &\le C_1h^{-1/2}u_t(x)+C_2w_t(x)
       +C_1w_t(x)\int_\rho^h r^{-1/2}\dd r\\
 &=C_1h^{-1/2}u_t(x)
       +\bigl[C_2+2C_1(\sqrt h-\sqrt\rho)\bigr]w_t(x)\\
 &\le C_3\bigl(h^{-1/2}u_t(x)+h^{1/2}w_t(x)\bigr).
 \end{aligned}
\]
Here we used $A>0$ and $h\ge\rho\ge1$.
Dividing by $u_t(x)>0$, we obtain
\[
 \frac{\sqrt{\Gamma u_t(x)}}{u_t(x)}
 \le C_3\left(h^{-1/2}+h^{1/2}\frac{w_t(x)}{u_t(x)}\right).
\]

\Needspace{7\baselineskip}
As $\Gamma(\sqrt{u_t})/u_t\le\Gamma u_t/u_t^2$, choosing $
 h=\min\{t/4,u_t(x)/w_t(x)\}$ if $w_t(x)>0$
gives
\[
 \frac{\Gamma(\sqrt{u_t})(x)}{u_t(x)}
 \le C\left(\frac{w_t(x)}{u_t(x)}+\frac1t\right).
\]
Use $h=t/4$ if $w_t(x)=0$.

Combining the two cases and choosing $\beta\ge1$ sufficiently large gives
\[
 \begin{aligned}
 \frac{\Gamma(\sqrt{u_t})(x)}{u_t(x)}
 &\le\beta\left(\frac{w_t(x)}{u_t(x)}+\frac1t\right)=\beta\partial_t\log u_t(x)+\frac{\beta(A+1)}t.
 \end{aligned}
\]
Taking $C=\beta(A+1)$ proves the lemma.
\end{proof}

\begin{lemma}
\label{lem:sharp-positive-regularity}
There are constants $C,t_{\mathrm{reg}}<\infty$ such that every
nonnegative bounded solution of the heat equation $\partial_tu=Lu$
on $(0,\infty)\times V$ satisfies
\begin{equation*}
 \sqrt{\Gamma u(t,x)}\le Ct^{-1/2}u(2t,x),\qquad
 |Lu(t,x)|\le Ct^{-1}u(2t,x)
 \quad(t\ge t_{\mathrm{reg}}).
\end{equation*}
In addition, for every fixed $\Lambda\ge2$,
\begin{equation}\label{eq:sharp-time-harnack}
 u(a,x)\le C_\Lambda u(b,x)
 \quad\text{if }0<a<b\le\Lambda a.
\end{equation}
The constants depend only on $\eps,\kappa$, and on $\Lambda$ in
\eqref{eq:sharp-time-harnack}.
\end{lemma}
\begin{proof}
We first prove \eqref{eq:sharp-time-harnack}.
Set
\[
 \widetilde u(s,x)=u(s/4,x),\qquad
 \widetilde P=\tfrac34I+\tfrac14P,\qquad
 \widetilde L=\widetilde P-I=\tfrac14L.
\]
Then $\partial_s\widetilde u=\widetilde L\widetilde u$.
The Markov kernel $\widetilde P$ has the same support graph and
reversible measure $m$ as $P$, and satisfies \eqref{eq:geometry}
with $\eps$ replaced by $\eps/4$ and the same ratio bound $\kappa$.
Moreover, $\widetilde\Gamma_2=\Gamma_2/16$, so
$\CD_{\widetilde L}(0,\infty)$ holds.
Thus Theorem~\ref{thm:parabolic-harnack} applies to $\widetilde u$
with constant $H=C_{\mathrm{PH}}(\eps/4,\kappa)$.

First suppose $0<a<b\le2a$. If $b-a\le1$, then
$\partial_tu=Pu-u\ge-u$ gives
\[
 u(a,x)\le e^{b-a}u(b,x)\le e\,u(b,x).
\]
If $b-a>1$, put $A=4a$, $B=4b$, and choose
\[
 \rho^2=\frac{B-A}{4}=b-a>1,\qquad
 \tau=A-\rho^2\ge3a>0.
\]
Since $A=\tau+\rho^2$ and $B=\tau+5\rho^2$, applying that theorem
to $\widetilde u(\tau+\cdot,\cdot)$ with radius $\rho$ gives
$\widetilde u(A,x)\le H\widetilde u(B,x)$.
Thus $u(a,x)\le K u(b,x)$ whenever $0<a<b\le2a$, where
$K=\max\{e,H\}$.

For general $0<a<b\le\Lambda a$, set
$N=\lceil\log_2(b/a)\rceil\ge1$.
The times $s_\ell=2^\ell a$ for $0\le\ell<N$, followed by $s_N=b$,
have successive ratios in $(1,2]$. Applying the preceding comparison
along this chain yields
\[
 u(a,x)\le K^N u(b,x)
 \le K^{\lceil\log_2\Lambda\rceil}u(b,x),
\]
which proves \eqref{eq:sharp-time-harnack} with the stated dependence
of the constant.

For the derivative estimates, fix $t>0$ and $x\in V$, and write
\[
 \mathcal Q_\rho=(t-\rho^2,t]\times B(x,\rho),\qquad
 \iint_{I\times U}F=\int_I\sum_{y\in U}m(y)F(s,y)\dd s.
\]
The $L^1$ mean-value inequality obtained from the
Moser argument in \cite[Section~2]{Delmotte} and interpolation gives,
for every nonnegative heat subsolution $(\partial_s-L)w\le0$,
\begin{equation*}
 w(t,x)\le\frac{C}{r^2\vol(B(x,r))}
                  \iint_{\mathcal Q_{2r}}w
 \qquad(r\ge2,\ t>4r^2).
\end{equation*}

For $r\ge2$ and $t>(c_0r)^2$, where $c_0=32$, we claim
\begin{align}
 \iint_{\mathcal Q_{2r}}\Gamma u
   &\le Cr^{-2}\iint_{\mathcal Q_{c_0r}}u^2,
       \notag\\
 \iint_{\mathcal Q_{2r}}|Lu|^2
   &\le Cr^{-4}\iint_{\mathcal Q_{c_0r}}u^2.
       \label{eq:sharp-caccioppoli-two}
\end{align}
To this end, let $R\ge2$ and $t>16R^2$, and set
$I_R=[t-16R^2,t]$. Lemma~\ref{lem:second-cutoffs}, applied with
radius $3R/2$, gives a cutoff $0\le\eta\le1$ such that $\eta=1$ on
$B(x,3R/2)$, $\operatorname{supp}\eta\subset B(x,3R)$, and
$\Gamma\eta+|L\eta|\le CR^{-2}$.
Choose $0\le\chi\le1$ on $I_R$ with $\chi(t-16R^2)=0$,
$\chi=1$ on $[t-R^2,t]$, and $|\chi'|\le CR^{-2}$.
The functions $\eta$, $L\eta$, and $\Gamma\eta$ vanish outside
$B(x,4R)$. Testing $(\partial_s-L)u^2=-2\Gamma u$ against
$\chi\eta$ and integrating by parts gives
\[
 \begin{aligned}
 2\iint_{\mathcal Q_R}\Gamma u
 &\le2\iint_{I_R\times V}\chi\eta\Gamma u\\
 &\le\iint_{I_R\times V}(\chi'\eta+\chi L\eta)u^2
 \le CR^{-2}\iint_{\mathcal Q_{4R}}u^2,
 \end{aligned}
\]
where the initial-time term vanishes and the final-time term is nonpositive.

With the same cutoffs, put $v=\chi\eta u$. The product formula gives
\[
 F:=(\partial_s-L)v
   =\chi'\eta u-\chi uL\eta-2\chi\Gamma(\eta,u).
\]
Write $\|\cdot\|_2=\|\cdot\|_{\ell^2(V,m)}$.
Since $v$ vanishes at the left endpoint of $I_R$, reversibility gives
\[
 \begin{aligned}
 \int_{I_R}\|F\|_2^2\dd s
 &=\int_{I_R}\bigl(\|\partial_sv\|_2^2+\|Lv\|_2^2\bigr)\dd s
   +\cE(v(t))\\
 &\ge\int_{I_R}\|Lv\|_2^2\dd s.
 \end{aligned}
\]
The product formula and
$|\Gamma(\eta,u)|^2\le\Gamma\eta\,\Gamma u$ imply
\[
 \int_{I_R}\|F\|_2^2\dd s
 \le CR^{-4}\iint_{\mathcal Q_{4R}}u^2
       +CR^{-2}\iint_{\mathcal Q_{4R}}\Gamma u.
\]
Since $R\ge2$, $\eta=1$ on $B(x,R)$ and its neighbouring vertices,
so $Lv=Lu$ on $\mathcal Q_R$. If $t>(16R)^2$, applying the preceding
gradient energy bound at radius $4R$ therefore yields
\[
 \iint_{\mathcal Q_R}|Lu|^2
 \le CR^{-4}\iint_{\mathcal Q_{16R}}u^2.
\]
Taking $R=2r$ proves \eqref{eq:sharp-caccioppoli-two}.

Finally, the curvature assumption gives
\[
 (\partial_s-L)\Gamma u=-2\Gamma_2u\le0,\qquad
 (\partial_s-L)(Lu)^2=-2\Gamma(Lu)\le0.
\]
Thus both $\Gamma u$ and $(Lu)^2$ are nonnegative heat subsolutions.
Set $t_{\mathrm{reg}}=16c_0^2$. For $t\ge t_{\mathrm{reg}}$, take
$r=\sqrt t/(2c_0)\ge2$, so that
$\mathcal Q_{c_0r}=(3t/4,t]\times B(x,\sqrt t/2)$.
Theorem~\ref{thm:parabolic-harnack}, applied with radius $\sqrt{t/3}$
and initial time $t/3$, gives
\[
 \sup_{\mathcal Q_{c_0r}}u\le Cu(2t,x).
\]
The mean-value inequality for these two subsolutions, together with
\eqref{eq:sharp-caccioppoli-two} and volume doubling, yields
\[
 \begin{aligned}
 \Gamma u(t,x)&\le Cr^{-2}u(2t,x)^2\le Ct^{-1}u(2t,x)^2,\\
 |Lu(t,x)|^2&\le Cr^{-4}u(2t,x)^2\le Ct^{-2}u(2t,x)^2.
 \end{aligned}
\]
\end{proof}

\begin{lemma}
\label{lem:ly-poisson-ratios}
For every $t>0$ and $x\sim y$,
\[
 \frac{\eps}{1+\partial_t\log u_t(x)+t^{-1}}
 \le\frac{u_t(y)}{u_t(x)}
 \le\frac{1+\partial_t\log u_t(x)}{\eps}.
\]
In particular, $\partial_t\log u_t(x)>-1$.
\end{lemma}
\begin{proof}
The upper bound follows from
$Pu_t(x)=(1+\partial_t\log u_t(x))u_t(x)$.

For the lower bound, use the argument in the proof of
Lemma~\ref{lem:poisson-insertion} with a single edge.
The uniformization formula defines a probability distribution on $\mathbb{N}$
\[
 \nu_j=\frac{e^{-t}t^jP^jf(x)}{j!u_t(x)},\qquad
 \sum_{j\ge0}j\nu_j=t(1+\partial_t\log u_t(x)).
\]
Using $P^{j+1}f(y)\ge p(y,x)P^jf(x)$ and Jensen's inequality gives
\[
 \frac{u_t(y)}{u_t(x)}
 \ge\eps t\sum_{j\ge0}\frac{\nu_j}{j+1}
 \ge\frac{\eps t}{t(1+\partial_t\log u_t(x))+1}.
\]
\end{proof}

\begin{proof}[Proof of Theorem~\ref{thm:li-yau-full}]
The chain with exact one vertex is immediate, so assume that $V$ has at least
two vertices; in particular, $0<\eps\le1$.
Fix $t>0$ and $x\in V$, and write
\[
 Q=\max\{\Gamma(\log u_t)(x),\Psi(\log u_t)(x)\}.
\]
Set $\alpha_0=12\beta/\eps$. Lemma~\ref{lem:ly-poisson-ratios}
gives $\partial_t\log u_t(x)>-1$.

For $0<\delta\le1$ and $r\ge\delta$, we have
\[
 \max\left\{\frac12(\log r)^2,\ r-1-\log r\right\}
 \le\frac2\delta(\sqrt r-1)^2.
\]
These two scalar bounds follow, respectively, from the mean value
theorem for $\log s$ and Taylor's theorem for $s^2-1-2\log s$,
with $s=\sqrt r$.

If $t\ge1$ and $\partial_t\log u_t(x)\le1$,
Lemma~\ref{lem:ly-poisson-ratios} gives
$u_t(y)/u_t(x)\ge\eps/3$ for $y\sim x$. The scalar inequality
with $\delta=\eps/3$ and Lemma~\ref{lem:ly-sqrt} imply
\[
 Q\le\frac{12}{\eps}\frac{\Gamma(\sqrt{u_t})(x)}{u_t(x)}
 \le\alpha_0\partial_t\log u_t(x)+\frac{12C}{\eps t}.
\]

In the remaining cases, $t<1$ or $\partial_t\log u_t(x)>1$, put
$D=(1+\partial_t\log u_t(x)+t^{-1})/\eps>1$.
Lemma~\ref{lem:ly-poisson-ratios} gives
$D^{-1}\le u_t(y)/u_t(x)\le D$ for $y\sim x$.
Using $(\log D)^2\le D$, we obtain
\[
 \begin{aligned}
 \Gamma(\log u_t)(x)&\le\tfrac12(\log D)^2\le\frac D2,\\
 \Psi(\log u_t)(x)
 &=\partial_t\log u_t(x)
   -\sum_{y\sim x}p(x,y)\log\frac{u_t(y)}{u_t(x)}\\
 &\le\partial_t\log u_t(x)+\log D
   \le\partial_t\log u_t(x)+D.
 \end{aligned}
\]
Since $\partial_t\log u_t(x)>-1$, both bounds imply
$Q\le1+\partial_t\log u_t(x)+D$.

If $t\ge1$ and $\partial_t\log u_t(x)>1$, then
$D\le3\partial_t\log u_t(x)/\eps$ and
$1+\partial_t\log u_t(x)\le2\partial_t\log u_t(x)$, so
\[
 Q\le\frac5{\eps}\partial_t\log u_t(x)
 \le\alpha_0\partial_t\log u_t(x).
\]

If $t<1$, we instead have
\[
 Q-\alpha_0\partial_t\log u_t(x)
 \le(1-\alpha_0+\eps^{-1})\partial_t\log u_t(x)
       +1+\eps^{-1}+\frac1{\eps t}.
\]
Since $1-\alpha_0+\eps^{-1}<0$ and $\partial_t\log u_t(x)>-1$, we obtain
\[
 Q-\alpha_0\partial_t\log u_t(x)\le\alpha_0+\frac1{\eps t}
 \le\frac{\alpha_0+\eps^{-1}}t.
\]

Thus $Q-\alpha_0\partial_t\log u_t(x)\le C/t$ in all cases.
For $\alpha\ge\alpha_0$, use $Q\ge0$ to obtain
\[
 \begin{aligned}
 Q-\alpha\partial_t\log u_t(x)
 &\le\frac{\alpha}{\alpha_0}
           \bigl(Q-\alpha_0\partial_t\log u_t(x)\bigr)\le\frac{C\alpha}{t}.
 \end{aligned}
\]
Choosing $C_{LY}$ large enough proves both assertions.
\end{proof}

\begin{remark}
\label{rem:ly-literature}
Several versions of the Li--Yau gradient estimate have been
established on graphs under different curvature assumptions.
For solutions of the heat equation, Bauer et al.~\cite{BauerEtAl2015} proved
\[
 \frac{\Gamma(\sqrt u)}u-\frac12\partial_t\log u
 =-\frac{L\sqrt u}{\sqrt u}\le\frac n{2t}
\]
under $\mathrm{CDE}(n,0)$, on finite graphs and infinite
bounded-degree graphs. Dier, Kassmann, and Zacher
\cite{DierKassmannZacher2021} obtained
$-L\log u\le\varphi(t)$ on finite connected graphs under
$\CD(F;0)$, where $\varphi'+F(\varphi)=0$ and
$\varphi(0+)=\infty$. M\"unch~\cite{Munch2018Nonlinear}
obtained $-L\log u\le n/(2t)$ on finite unweighted graphs
under $\mathrm{CD}^{\log}(n,0)$.

For the modified heat flow
$\partial_t v_t=Lv_t+\Gamma(v_t)$,
M\"unch~\cite{Munch2019CD} proved
\[
 \Gamma(v_t)-\partial_t v_t=-Lv_t\le\frac n{2t}
\]
under the classical condition $\CD(0,n)$ on finite graphs.
Russ and Pajot~\cite{RussPajot} extended this estimate to infinite
weighted graphs of bounded geometry under the same curvature condition.
Guo, Huang, and Huang~\cite{GuoHuangHuang} established
\[
 \Gamma(v_t)-\partial_t v_t=-Lv_t
 \le\frac{d_*(d_*+3)}{2t}
\]
under $\CD(0,\infty)$ on connected simple
unweighted graphs of maximum degree $d_*\ge2$, using the
non-normalized Laplacian. These modified-flow results assume $\norm{\Gamma v_0}_{\infty}$ is small enough; $v_t$ is generally not $\log(e^{tL}f)$.
\end{remark}

\subsection{A global Harnack inequality}

\begin{theorem}
\label{thm:ly-global-harnack}
Under the hypotheses of Theorem~\ref{thm:li-yau-full}, put $u_t=P_tf$. For
$\alpha\ge\alpha_0$, $0<s<t$, and $x,y\in V$,
\[
 u_s(x)\le u_t(y)\left(\frac ts\right)^{C_{LY}}
 \exp\left(\frac{2\alpha\dist(x,y)^2}{\eps(t-s)}\right).
\]
\end{theorem}
\begin{proof}
For $s\le r\le t$, put $F_r=\log u_r+C_{LY}\log(r/s)$.
Theorem~\ref{thm:li-yau-full} gives
\[
 \partial_rF_r\ge\alpha^{-1}\Gamma F_r.
\]
Applying Lemma~\ref{lem:differential-harnack} with $c(r)=1/\alpha$
gives
\[
 F_s(x)-F_t(y)
 \le\frac{2\alpha\dist(x,y)^2}{\eps(t-s)}.
\]
Substituting the definition of $F_r$, we obtain
\[
 \log\frac{u_s(x)}{u_t(y)}
 \le C_{LY}\log\frac ts
   +\frac{2\alpha\dist(x,y)^2}{\eps(t-s)}.
\]
Exponentiating proves the assertion.
\end{proof}

\section{A sharp volume growth}
\label{sec:sharp-volume}

We prove Theorem~\ref{thm:sharp-volume}. Recall that $d=\sup_x\deg(x)$ and
$n_*=\lfloor d/2\rfloor$, with $d\le\lfloor\eps^{-1}\rfloor$. 

\subsection{Harmonic functions with constant gradient norm}
Fix a root $o$ and put
\begin{equation}\label{eq:flat-function-space}
 \mathcal F_o=\{f:V\to\R:f(o)=0,\ Lf=0,
                     \ \Gamma f\text{ is a finite constant}\}.
\end{equation}
\begin{proposition}
\label{prop:complete-zero-modes}
Under \eqref{eq:geometry} and $\CD(0,\infty)$,
\begin{equation}\label{eq:complete-zero-modes}
 \{u:\Gamma_2u=0\}=\R\1+\mathcal F_o.
\end{equation}
Here $\mathcal F_o$ is defined in \eqref{eq:flat-function-space}.
\end{proposition}
\begin{proof}
Suppose $\Gamma_2u=0$. Positive semidefiniteness under
$\CD(0,\infty)$ gives $\Gamma_2(u,h)=0$ for every real function $h$.
For distinct vertices $x$ and $z$,
\[
 \begin{aligned}
 0&=4m(x)\Gamma_2(u,\1_{\{z\}})(x)\\
 &=\sum_y m(x)p(x,y)p(y,z)\bigl(u(x)+u(z)-2u(y)\bigr)\\
 &\quad-2m(x)p(x,z)\bigl(Lu(z)-Lu(x)\bigr)\\
 &=\sum_y m(z)p(z,y)p(y,x)\bigl(u(z)+u(x)-2u(y)\bigr)\\
 &\quad-2m(z)p(z,x)\bigl(Lu(z)-Lu(x)\bigr)\\
 &=-4m(x)p(x,z)\bigl(Lu(z)-Lu(x)\bigr).
 \end{aligned}
\]
Hence $Lu(x)=Lu(z)$
whenever $x\sim z$, so connectedness gives $Lu=c$ for a constant $c$.

Consequently, $L\Gamma u=2\Gamma(u,Lu)=0$. By
Corollary~\ref{cor:positive-liouville}, $\Gamma u$ is a constant
$A\ge0$. The edge lower bound therefore gives
\[
 |u(v)-u(w)|\le M:=\sqrt{2A/\eps}\qquad(v\sim w).
\]

For integers $r\ge1$, we have
\[
 \begin{aligned}
 |c|\vol(B(o,r))
 &=\left|\sum_{\substack{v\in B(o,r)\\w\notin B(o,r)}}
       m(v)p(v,w)\bigl(u(w)-u(v)\bigr)\right|\\
 &\le M\bigl(\vol(B(o,r))-\vol(B(o,r-1))\bigr).
 \end{aligned}
\]
Summing over $r=R+1,\ldots,2R$ and using the volume doubling
constant $C_V$ from Theorem~\ref{thm:main}, we obtain
\[
 \begin{aligned}
 |c|R\vol(B(o,R))
 &\le M\bigl(\vol(B(o,2R))-\vol(B(o,R))\bigr)\\
 &\le MC_V\vol(B(o,R)).
 \end{aligned}
\]
Letting $R\to\infty$ gives $c=0$. Thus $u-u(o)\in\mathcal F_o$.
The reverse inclusion follows from the definition.
\end{proof}

For $f\in\mathcal F_o$, \eqref{eq:complete-zero-modes} gives
$\Gamma_2f=0$.
At each vertex $x$, $\CD(0,\infty)$ makes the form
$\Gamma_2(\cdot,\cdot)(x)$ positive semidefinite. Hence, for every
real function $u$,
\[
 |\Gamma_2(f,u)(x)|^2
 \le\Gamma_2f(x)\Gamma_2u(x)=0.
\]
Since $Lf=0$, the definition of $\Gamma_2$ then gives
\begin{equation}\label{eq:flat-polarization}
 \Gamma_2(f,u)=0,\qquad L\Gamma(f,u)=\Gamma(f,Lu).
\end{equation}

In particular, $\Gamma_2(af+bg)=0$ for $f,g\in\mathcal F_o$ and
$a,b\in\R$. Since $af+bg$ is harmonic and vanishes at $o$,
Proposition~\ref{prop:complete-zero-modes} implies
$af+bg\in\mathcal F_o$. Thus $\mathcal F_o$ is a vector space.
Polarization now gives
$\Gamma(f,g)=\tfrac14(\Gamma(f+g)-\Gamma(f-g))$, which is constant;
its value defines an inner product.
At every vertex $x$, the map $f\mapsto(f(y)-f(x))_{y\sim x}$ is
injective: a zero gradient makes the constant energy zero, so $f=0$.
Hence $\mathcal F_o$ is finite-dimensional.

Let $f\in\mathcal F_o$. For real $\theta$,  define kernels
\begin{equation}\label{eq:flat-phase-twists}
 B_f(\theta)(x,y)=p(x,y)e^{i\theta(f(x)-f(y))}
\end{equation}
The following commutator theorem for normal operators will be used to
establish pairwise commutativity of the family
$\{B_f(\theta):\theta\in\mathbb R\}$.

\begin{theorem}[{\cite[Theorem~0.5]{AckermansVanEijndhovenMartens1983}}]
\label{thm:double-commutator}
Let $D$ and $E$ be bounded operators on a complex Hilbert space, with
$D$ normal. Then if $[D,[D,E]]=0$, then $[D,E]=0$ where $[D,E]=DE-ED$.
\end{theorem}

\begin{lemma}[Commuting phase twists]
\label{lem:flat-twists}
Let $f\in\mathcal F_o$, with $\mathcal F_o$ as in
\eqref{eq:flat-function-space}. The operators $B_f(\theta)$ defined by
\eqref{eq:flat-phase-twists} are bounded and self-adjoint on
$\ell^2(V,m;\mathbb C)$ for real $\theta$, and
$[B_f(\theta),B_f(\eta)]=0$ for all real $\theta,\eta$.
\end{lemma}
\begin{proof}
Set $M_f=\sqrt{2\Gamma f(o)/\eps}$ and define
\[
 C_f(x,y)=p(x,y)(f(x)-f(y)).
\]
The edge bound gives $|C_fu|\le M_fP|u|$, so $C_f$ is bounded.
Since $Lf=0$, we have $C_fu=-2\Gamma(f,u)$.
Consequently, \eqref{eq:flat-polarization} yields
\[
 [P,C_f]u=-2\bigl(L\Gamma(f,u)-\Gamma(f,Lu)\bigr)=0.
\]
These identities hold first on finitely supported functions and
extend by boundedness.

Put $B(\theta)=B_f(\theta)$. Its differentiated kernels satisfy
$\|B^{(j)}(\theta)\|\le M_f^j$ for $j\ge1$, and the exponential
series converges in operator norm. Thus $B$ is entire in $\theta$, and
reversibility makes it self-adjoint for real $\theta$. Direct
multiplication of the kernels gives
\[
 [B(\theta),B'(\theta)](x,z)
 =i e^{i\theta(f(x)-f(z))}[P,C_f](x,z)=0.
\]

We prove by induction that $[B(\theta),B^{(j)}(\theta)]=0$ for every
real $\theta$ and every $j\ge1$. The case $j=1$ has just been proved.
Suppose the identity holds for some $j\ge1$ and all real $\theta$.
Differentiation gives $[B,B^{(j+1)}]=-[B',B^{(j)}]$. The Jacobi identity, together with $[B,B']=0$ and the induction
hypothesis, then yields
\[
 \begin{aligned}
 [B,[B,B^{(j+1)}]]
 &=-[B,[B',B^{(j)}]]\\
 &=-[[B,B'],B^{(j)}]-[B',[B,B^{(j)}]]=0.
 \end{aligned}
\]
For each real $\theta$, the operator $B(\theta)$ is self-adjoint and
hence normal. Applying Theorem~\ref{thm:double-commutator} with
$D=B(\theta)$ and $E=B^{(j+1)}(\theta)$ gives
$[B(\theta),B^{(j+1)}(\theta)]=0$, completing the induction.

Finally, for real $\theta,\eta$, the operator-norm convergent expansion
\[
 B(\eta)=\sum_{j=0}^{\infty}
 \frac{(\eta-\theta)^j}{j!}B^{(j)}(\theta)
\]
gives $[B(\theta),B(\eta)]=0$ by termwise commutation.
\end{proof}

The commuting kernels will identify the distribution of the difference of spectral embedding in the next lemma.
Let $r=\dim\mathcal F_o$. Choose an orthonormal basis in the inner
product $\Gamma(\cdot,\cdot)$ and let
$\Phi:V\to\R^r$ be its spectral embedding, with $\Phi(o)=0$.
Set
\begin{equation}\label{eq:flat-increment-set}
 S_0=\{\Phi(x)-\Phi(y):x,y\in V,\ p(x,y)>0\}.
\end{equation}

\begin{lemma}\label{lem:flat-increments}
For the embedding $\Phi:V\to\R^r$ fixed above, the set $S_0$ in
\eqref{eq:flat-increment-set} is finite and symmetric, that is, $S_0=-S_0$.
There are positive numbers $k_s$, $s\in S_0$, such that
\begin{equation}\label{eq:flat-increment-law}
 \sum_{y\in V:\,\Phi(y)-\Phi(x)=s}p(x,y)=k_s
 \qquad(x\in V,\ s\in S_0),
\end{equation}
and
\begin{equation}\label{eq:flat-increment-moments}
 \sum_{s\in S_0}k_s=1,\qquad
 \sum_{s\in S_0}k_s\,s=0,
 \qquad \frac12\sum_{s\in S_0}k_s\,ss^{\mathsf T}=I_r.
\end{equation}
The set $S_0$ spans $\R^r$. For $r=0$, these statements mean
$S_0=\{0\}$ and $k_0=1$.
\end{lemma}
\begin{proof}
By the definition of $S_0$, $S_0$ is symmetric.  For $f\in\mathcal F_o$ and real $\theta$,
Lemma~\ref{lem:flat-twists} gives $[P,B_f(\theta)]=0$.
Both product kernels have finite rows, so this identity also holds
when applied to $\1$. Hence the bounded function $g=B_f(\theta)\1$
satisfies
\[
 Pg=PB_f(\theta)\1=B_f(\theta)P\1=g.
\]
Thus $g$ is harmonic, since $P\1=\1$. Corollary~\ref{cor:positive-liouville},
applied to its real and imaginary parts after adding constants,
makes $g$ constant.

The finite probability measures
\[
 \mu_x=\sum_yp(x,y)\delta_{\Phi(y)-\Phi(x)}
\]
therefore have the same characteristic function, since
$B_{a\cdot\Phi}(-1)\1(x)=\int e^{ia\cdot s}\,\mu_x(ds)$ for
all $a\in\R^r$. Hence $\mu_x$ is independent of $x$. Thus $S_0$ is the common finite support of $\mu_x$.
Setting $k_s=\mu_x(\{s\})$
for $s\in S_0$ gives \eqref{eq:flat-increment-law},
with $k_s>0$ for $s\in S_0$. The normalization of $\mu_x$,
harmonicity, and orthonormality give \eqref{eq:flat-increment-moments},
including the spanning assertion.
\end{proof}

\begin{lemma}
\label{lem:sharp-zero-rank}
Under \eqref{eq:geometry} and $\CD(0,\infty)$,
\[
 \dim\mathcal F_o\le
 \left\lfloor\frac{\deg(x)}2\right\rfloor
 \qquad(x\in V).
\]
Here $\mathcal F_o$ is defined in \eqref{eq:flat-function-space}.
\end{lemma}
\begin{proof}
Let $r=\dim\mathcal F_o$ and apply Lemma~\ref{lem:flat-increments} to
an orthonormal basis in the inner product $\Gamma(\cdot,\cdot)$. The set $S_0$ is symmetric and spans
$\R^r$, and every nonzero increment in $S_0$ occurs at every vertex.
It therefore contains at
least $r$ nonzero opposite pairs, requiring at least $2r$ distinct
neighbours at each vertex.
\end{proof}
\begin{remark}[Dimension estimates for harmonic functions]
On complete connected Riemannian manifolds with nonnegative Ricci
curvature, Colding and
Minicozzi~\cite{ColdingMinicozzi1997} proved finite-dimensionality for
each fixed polynomial growth order.  Later, Peter Li~\cite{Li1997HarmonicSections}, Colding and
Minicozzi~\cite{ColdingMinicozzi1998}
obtained the optimal dimension estimates. Kleiner~\cite{Kleiner2010} gave a new proof of finite-dimensionality
for spaces of harmonic functions of polynomial growth
on Cayley graphs of groups of polynomial growth, as a key ingredient
in his proof of Gromov's theorem.

For infinite semiplanar graphs with nonnegative combinatorial curvature,
Hua, Jost, and Liu~\cite{HuaJostLiu2015} obtained a quadratic dimension
bound in the growth order. Hua and Jost~\cite{HuaJost2015SemiplanarII}
subsequently improved this to the optimal linear bound. Hua~\cite{Hua2021UnboundedCaloric} also established dimension comparisons
between spaces of polynomial-growth ancient caloric and harmonic
functions on graphs.
\end{remark}

Define the weighted root gradient by
\begin{equation}\label{eq:sharp-root-gradient}
 D_x f=
  \bigl(\sqrt{p(x,y)}(f(y)-f(x))\bigr)_{y\sim x},
 \qquad \|D_xf\|^2=2\Gamma f(x).
\end{equation}
For an integer $R\ge1$, let
\begin{equation}\label{eq:sharp-local-quadratic-form}
 \mathfrak q_{x,R}(f)=
 \sum_{v\in B(x,R)}\Gamma_2f(v).
\end{equation}
This is a nonnegative quadratic form on functions on
$B(x,R+2)$. 

To reduce $\mathfrak q_{x,R}$ to a quadratic form in the root
gradient $D_xf$, we first record the following linear-algebra lemma.

\begin{lemma}
\label{lem:quadratic-minimization}
Let $q$ be a nonnegative quadratic form on $\R^k\times\R^\ell$,
where $k,\ell\ge0$. For every $a\in\R^k$, the minimum
$\min_{z\in\R^\ell}q(a,z)$ is attained. Moreover, there is a unique
symmetric positive semidefinite $k\times k$ matrix $M$ such that
\[
 \min_{z\in\R^\ell}q(a,z)=a^{\mathsf T}Ma
 \qquad(a\in\R^k).
\]
\end{lemma}
\begin{proof}
Let $H$ be the positive semidefinite matrix of $q$ and factor
$H=T^{\mathsf T}T$. Partition $T=(U\ V)$ by columns according to
the variables $(a,z)$; both $U$ and $V$ retain all rows of $T$.
Then
\[
 q(a,z)=\left\|T\binom{a}{z}\right\|^2=\|Ua+Vz\|^2.
\]
Let $Q$ be the orthogonal projection onto the finite-dimensional
closed subspace $\operatorname{im}V$. Orthogonality gives
\[
 q(a,z)=\|(I-Q)Ua\|^2+\|QUa+Vz\|^2.
\]
Since $QUa\in\operatorname{im}V$, one can choose $z$ with
$Vz=-QUa$. Thus the minimum is attained and equals
$\|(I-Q)Ua\|^2$. The matrix
\[
 M=U^{\mathsf T}(I-Q)U
\]
is symmetric positive semidefinite and satisfies the stated identity.
Its uniqueness follows by polarization.
\end{proof}

\begin{lemma}
\label{lem:sharp-local-projection}
For $D_x$ and $\mathfrak q_{x,R}$ defined in
\eqref{eq:sharp-root-gradient} and \eqref{eq:sharp-local-quadratic-form},
there are $R_*\in\N$ and $C_*<\infty$, depending only on $\eps,\kappa$,
such that at every vertex $x$ one can choose an orthogonal projection
$\Pi_x:\ \R^{\deg(x)}\to \R^{\deg(x)}$ with
\begin{equation}\label{eq:sharp-projection-rank}
 \operatorname{rank}\Pi_x\le
       \left\lfloor\frac{\deg(x)}2\right\rfloor
\end{equation}
and, for all real $f$,
\begin{equation}\label{eq:sharp-local-control}
 \|(I-\Pi_x)D_xf\|^2\le C_*\cdot\mathfrak q_{x,R_*}(f).
\end{equation}
\end{lemma}
\begin{proof}
Fix $x,R$ and put $k=\deg(x)$, $r_x=\lfloor k/2\rfloor$.
The case $k=0$ is trivial. Now assume $k\ge1$. For $a\in \mathbb{R}^k$, 
Lemma~\ref{lem:quadratic-minimization} then gives a symmetric
positive semidefinite matrix $M_{x,R}$ such that
\[
 a^{\mathsf T}M_{x,R}a
  =\min\{\mathfrak q_{x,R}(f):f(x)=0,\ D_xf=a\}.
\]

We claim that there are uniform $R_*,c_*>0$ such that, with eigenvalues
listed in increasing order,
\begin{equation}\label{eq:sharp-local-gap}
 \lambda_{r_x+1}(M_{x,R_*})\ge c_*.
\end{equation}

Otherwise, for each integer $j$ there would be a rooted chain
$(P_j,m_j,x_j)$ and $r_{x_j}+1$ functions $f_{j,i}$ on its radius-$(j+2)$
ball such that
\[
 f_{j,i}(x_j)=0,\qquad
 \bigl(D_{x_j}f_{j,i}\bigr)_i
     \text{ is orthonormal},\qquad
 \mathfrak q_{x_j,j}(f_{j,i})\le j^{-1}.
\]
Passing to a subsequence makes $k=\deg(x_j)$ constant; put
$r=\lfloor k/2\rfloor$.
Dividing each $m_j$ by its positive lower bound, we may assume
$1\le m_j\le\kappa$; the probabilities and curvature forms are unchanged.
For each fixed radius $s$, the degree bound allows only finitely
many isomorphism types of rooted support balls $B(x_j,s)$.
Passing to a subsequence fixes this type. After identifying
corresponding vertices, compactness of the non-diagonal edge
probabilities in $[\eps,1]$, diagonal probabilities in $[0,1]$,
and masses in $[1,\kappa]$ gives a further subsequence on which
all these coordinates converge.

Repeating this procedure for $s=1,2,\ldots$ with compatible
identifications and taking a diagonal subsequence gives a connected
rooted graph $(V_\infty,E_\infty,x_\infty)$ with limiting weights
$P_\infty,m_\infty$. The edge lower bound preserves the
support, so $\deg(x_\infty)=k$, and $1\le m_\infty\le\kappa$.

For each vertex, a sufficiently large fixed ball contains all its
neighbours, so the finite row sums and reversibility identities
pass to the limit. The forms $\Gamma_{2,j}(\cdot)(v)$ depend only
on the radius-two neighbourhood of $v$; their matrix coefficients
converge, so positive semidefiniteness also passes to the limit.
Thus the limiting chain satisfies \eqref{eq:geometry} with $m_-=1$
and $\CD(0,\infty)$.

To extract the functions, write $f=f_{j,i}$. Nonnegativity gives
$\Gamma_{2,j}f\le j^{-1}$ on $B(x_j,j)$. The square expansion in
Remark~\ref{rem:universal-curvature-bounds} yields
\[
 \sum_{y,z}p_j(v,y)p_j(y,z)
       \bigl(f(z)-2f(y)+f(v)\bigr)^2
 \le4\Gamma_jf(v)+4j^{-1}
 \qquad(v\in B(x_j,j)).
\]
For $v\sim w$ with $v\in B(x_j,j)$, it follows that
\[
 \begin{aligned}
 \Gamma_jf(w)
 &\le\sum_zp_j(w,z)\bigl(f(z)-2f(w)+f(v)\bigr)^2
       +\bigl(f(w)-f(v)\bigr)^2\\
 &\le6\eps^{-1}\Gamma_jf(v)+4(\eps j)^{-1}.
 \end{aligned}
\]
Since $\Gamma_jf(x_j)=1/2$, the energy is uniformly bounded on each
fixed rooted ball. The edge lower bound and $f(x_j)=0$ then bound the
function values there. A second diagonal extraction gives $r+1$
global functions $f_i$ on the limiting chain with
$f_i(x_\infty)=0$ and $\Gamma_2f_i=0$, whose weighted root gradients
remain orthonormal. Proposition~\ref{prop:complete-zero-modes} places
these functions in $\mathcal F_{x_\infty}$. This contradicts
Lemma~\ref{lem:sharp-zero-rank}, and proves that
\eqref{eq:sharp-local-gap} holds.

Let $\Pi_x$ be the spectral projection of $M_{x,R_*}$ onto $[0,c_*)$.
Its rank is at most $r_x$, and
\[
 \|(I-\Pi_x)a\|^2\le c_*^{-1}a^{\mathsf T}M_{x,R_*}a.
\]
The minimizing definition proves \eqref{eq:sharp-local-control}.
\end{proof}

\subsection{Smoothing markov kernels and Fisher information matrices}

To apply Lemma~\ref{lem:entropy-volume}, set
\begin{equation}\label{eq:sharp-fixed-mixture}
 a_j=(1-2^{-1/8})2^{-j/8},\qquad
 M=\sum_{j\ge0}a_jP_{2^j},\qquad q_t=P_tM.
\end{equation}
Since $a_j\ge0$ and $\sum_ja_j=1$, $M$ is a Markov kernel preserving $m$.
The semigroup property and reversibility give
\[
 \begin{aligned}
 q_t(x,z)&=\sum_{j\ge0}a_jp_{t+2^j}(x,z),\quad\partial_tq_t=L_xq_t,\quad
 m(x)q_t(x,z)&=m(z)q_t(z,x).
 \end{aligned}
\]

\begin{lemma}
\label{lem:sharp-mixture-regularity}
Let $q_t$ be the kernel in \eqref{eq:sharp-fixed-mixture}, and let
$t_{\mathrm{reg}}$ be as in Lemma~\ref{lem:sharp-positive-regularity}.
There are $t_0\ge t_{\mathrm{reg}}$ and $C<\infty$, depending only on
$\eps,\kappa$, such that
\begin{equation}\label{eq:sharp-mixture-entropy-coarse}
 H(q_t(x,\cdot))\le C+C\log(1+t)\qquad(t\ge0).
\end{equation}
Let $\eta_t=Ct^{-3/8}$ and $\zeta_t=Ct^{-7/8}$. For $t\ge t_0$ one has $\eta_t\le1/2$ and for $y\sim x$ 
\begin{equation}\label{eq:sharp-relative-derivatives}
 |q_t(y,z)-q_t(x,z)|\le\eta_tq_t(x,z),\quad
 |L_xq_t(x,z)|\le\zeta_tq_t(x,z),
\end{equation}
and for $v\in B(x,R_*+2)$
\begin{equation}\label{eq:sharp-fixed-neighbourhood}
       \tfrac12q_t(x,z)\le q_t(v,z)\le2q_t(x,z).
\end{equation}
Here $R_*$ is the fixed radius from Lemma~\ref{lem:sharp-local-projection}.
\end{lemma}
\begin{proof}
The entropy chain rule and \eqref{eq:curved-entropy}, applied with
$M=I$ and $K=0$, give
\[
 H(q_t(x,\cdot))
 \le H((a_j)_j)+\sum_ja_jH(p_{t+2^j}(x,\cdot))
 \le C+C\log(1+t).
\]
Here $\log(1+t+2^j)\le\log(1+t)+(j+1)\log2$, and both
$H((a_j)_j)$ and $\sum_ja_j(j+1)$ are finite.

For every fixed $\lambda>1$, We claim
\begin{equation}\label{eq:sharp-mixture-time-comparison}
 q_{\lambda t}(x,z)\le C_\lambda t^{1/8}q_t(x,z)
 \qquad(t\ge1).
\end{equation}
Set $J=\lceil\log_2t\rceil$ and choose $c=c(\lambda)\in\N$ with
$2^c\ge2\lambda+4$. Compare $\lambda t+2^j$ with $t+2^{J+c}$
when $j\le J$, and with $t+2^{j+c}$ when $j>J$. In both cases
the later time is between twice and a fixed multiple of the earlier
one. Then \eqref{eq:sharp-time-harnack} gives
\[
 \begin{aligned}
 q_{\lambda t}(x,z)
 &\le H_\lambda p_{t+2^{J+c}}(x,z)
      +H_\lambda\sum_{j>J}a_jp_{t+2^{j+c}}(x,z)\\
 &\le H_\lambda\bigl(a_{J+c}^{-1}+2^{c/8}\bigr)q_t(x,z)
 \le C_\lambda t^{1/8}q_t(x,z).
 \end{aligned}
\]
We used $a_j=2^{c/8}a_{j+c}$ and $2^J\le2t$.

For fixed $z$, the function $u(t,x)=q_t(x,z)>0$ is a bounded heat
solution. Apply Lemma~\ref{lem:sharp-positive-regularity}, use
\eqref{eq:sharp-mixture-time-comparison} with $\lambda=2$, and use
$p(x,y)|u(y)-u(x)|^2\le2\Gamma u(x)$.
This gives \eqref{eq:sharp-relative-derivatives}.

Since $\eta_t\to0$, increase $t_0$ so that
$\eta_t\le[2(R_*+2)]^{-1}$ for all $t\ge t_0$.
For $v\in B(x,R_*+2)$, choose a path
$x=x_0\sim x_1\sim\cdots\sim x_\ell=v$ with $\ell\le R_*+2$.
Applying the first estimate in \eqref{eq:sharp-relative-derivatives}
along this path gives
\[
 (1-\eta_t)^\ell
 \le\prod_{i=1}^\ell\frac{q_t(x_i,z)}{q_t(x_{i-1},z)}
 =\frac{q_t(v,z)}{q_t(x,z)}
 \le(1+\eta_t)^\ell.
\]
The outer terms satisfy $(1-\eta_t)^\ell\ge1-\ell\eta_t\ge1/2$
and $(1+\eta_t)^\ell\le e^{\ell\eta_t}\le e^{1/2}<2$, proving
\eqref{eq:sharp-fixed-neighbourhood}.
\end{proof}

Define
\begin{equation}\label{eq:sharp-fisher-functionals}
 \begin{aligned}
 J_t(x)=\sum_z\frac{\Gamma(q_t(\cdot,z))(x)}{q_t(x,z)},\quad A_t(x)=\sum_z\frac{\Gamma_2(q_t(\cdot,z))(x)}{q_t(x,z)}.
 \end{aligned}
\end{equation}

\begin{lemma}\label{lem:sharp-fisher-differentiation}
Let $q_t$ and $t_0$ be as in Lemma~\ref{lem:sharp-mixture-regularity},
and define $J_t,A_t$ by \eqref{eq:sharp-fisher-functionals}.
For every $T>t_0$, we have
$J\in C^1([t_0,T];\ell^\infty(V))$, and
\begin{equation}\label{eq:sharp-fisher-time-derivative}
 \begin{aligned}
 \partial_tJ_t(x)
 ={}&2\sum_z
       \frac{\Gamma(q_t(\cdot,z),Lq_t(\cdot,z))(x)}{q_t(x,z)}-\sum_z
       \frac{\Gamma(q_t(\cdot,z))(x)L_xq_t(x,z)}{q_t(x,z)^2}.
 \end{aligned}
\end{equation}
Both series converge absolutely. The series defining $A_t(x)$ also
converges absolutely, with
$\sup_{t\in[t_0,T],\,x\in V}|A_t(x)|<\infty$.
\end{lemma}
\begin{proof}
Fix $T>t_0$. For each $z$, write $u_t=q_t(\cdot,z)$.
On $[t_0,T]$, \eqref{eq:sharp-relative-derivatives} and
\eqref{eq:sharp-fixed-neighbourhood} give, uniformly in $x,z$,
\[
 \begin{aligned}
 |\Gamma_2u_t(x)|+|\Gamma(u_t,Lu_t)(x)|+\Gamma u_t(x)
 &\le Cu_t(x)^2,\\
 |Lu_t(x)|+|L^2u_t(x)|&\le Cu_t(x),\\
 \Gamma(Lu_t)(x)+|\Gamma(u_t,L^2u_t)(x)|&\le Cu_t(x)^2.
 \end{aligned}
\]
Indeed, the bound on $L^2u_t=L(Lu_t)$ follows from
$|Lu_t|\le Cu_t$ and comparison of neighbouring values of $u_t$;
the last line follows from the same comparison.

Since $\partial_tu_t=Lu_t$,
\[
 \partial_t\left(\frac{\Gamma u_t}{u_t}\right)
 =\frac{2\Gamma(u_t,Lu_t)}{u_t}
  -\frac{\Gamma u_t\,Lu_t}{u_t^2}.
\]
Differentiating once more and using these estimates gives
$|\partial_t^j(\Gamma u_t(x)/u_t(x))|\le Cu_t(x)$ for $j=0,1,2$.
Moreover, $|\partial_tu_t|\le Cu_t$ implies
$q_t(x,z)\le e^{C(T-t_0)}q_{t_0}(x,z)$. Thus
\[
 \sup_{t\in[t_0,T]}\sum_{j=0}^2
       \left|\partial_t^j\left(\frac{\Gamma u_t(x)}{u_t(x)}\right)\right|
 \le C_Tq_{t_0}(x,z),
 \qquad \sum_zq_{t_0}(x,z)=1.
\]
Consequently, the series defining $J_t(x)$ may be differentiated
twice term by term. The first-derivative formula gives
\eqref{eq:sharp-fisher-time-derivative}, and the preceding estimates
give absolute convergence of its two series.

Summing the second-derivative bound gives
$\sup_{t\in[t_0,T],\,x\in V}|\partial_t^2J_t(x)|\le C$;
the bounds for $j=0,1$ also show that $J_t$ and $\partial_tJ_t$
are uniformly bounded. Taylor's theorem therefore gives
\[
 \begin{aligned}
 \|J_{t+h}-J_t-h\partial_tJ_t\|_\infty&\le Ch^2,\\
 \|\partial_tJ_t-\partial_sJ_s\|_\infty&\le C|t-s|,
 \end{aligned}
\]
whenever the times belong to $[t_0,T]$. Hence
$J\in C^1([t_0,T];\ell^\infty(V))$.

Local finiteness also permits interchanging $L$ with the series.
The same estimates and $\sum_zq_t(x,z)=1$ give
$\sum_z|\Gamma_2(q_t(\cdot,z))(x)|/q_t(x,z)\le C$,
so $A_t$ is well defined and bounded.
\end{proof}

For $t>0$, we also define the following matrix by
\begin{equation}\label{eq:sharp-fisher-matrix}
 \mathbf F_t(x)_{y,z}
 =\sqrt{p(x,y)p(x,z)}
   \sum_u\frac{(q_t(y,u)-q_t(x,u))(q_t(z,u)-q_t(x,u))}{q_t(x,u)},
 \quad y,z\sim x,
\end{equation}
which is a discrete analogue of the Fisher information matrix;
see, for example, \cite[Eq.~(2.32)]{PolyanskiyWu2025}.
Heuristically, $\mathbf F_t(x)$ plays the role of the Hessian of
$D(q_t(y,\cdot)\|q_t(x,\cdot))$ with respect to the starting point $y$
at $y=x$. Set
\begin{equation}\label{eq:sharp-transverse-fisher}
 b_t(x)=\operatorname{tr}\bigl((I-\Pi_x)\mathbf F_t(x)(I-\Pi_x)\bigr).
\end{equation}

\Needspace{12\baselineskip}
\begin{lemma}
\label{lem:sharp-fisher-dissipation}
For the kernel $q_t$ in \eqref{eq:sharp-fixed-mixture} and the
projections $\Pi_x$ from Lemma~\ref{lem:sharp-local-projection}, use
$J_t,A_t,\mathbf F_t,b_t$ as defined in
\eqref{eq:sharp-fisher-functionals}--\eqref{eq:sharp-transverse-fisher}.
Let $d,n_*$ be as in Theorem~\ref{thm:sharp-volume}. For $t>0$,
\begin{equation}\label{eq:sharp-fisher-bound}
 0\le\mathbf F_t(x)\le t^{-1}I,\qquad
 J_t(x)=\tfrac12\operatorname{tr}\mathbf F_t(x)\le\frac{d}{2t}.
\end{equation}
For $t\ge t_0$, with $t_0$ as in Lemma~\ref{lem:sharp-mixture-regularity},
\begin{equation}\label{eq:sharp-fisher-rank-split}
 \operatorname{tr}\mathbf F_t(x)\le\frac{n_*}t+b_t(x).
\end{equation}
Moreover,
\begin{equation}\label{eq:sharp-integrated-errors}
 \sup_{T\ge t_0}\sup_x\int_{t_0}^T
       P_{T-s}A_s(x)\dd s\le C,\qquad
 \sup_{T\ge t_0}\sup_x\int_{t_0}^T
       P_{T-s}b_s(x)\dd s\le C.
\end{equation}
\end{lemma}
\begin{proof}
With $M$ as defined in \eqref{eq:sharp-fixed-mixture}, the reverse Poincar\'e inequality \eqref{eq:variance-K} and
Jensen's inequality give, for every bounded real function $f$,
\[
 2t\Gamma(q_tf)
 \le P_t((Mf)^2)-(q_tf)^2
 \le q_t(f^2)-(q_tf)^2.
\]

Fix $x$ and $t>0$. For
$a=(a_y)_{y\sim x}\in\R^{\deg(x)}$, define the signed measure
\[
 \nu_a(z)=\sum_{y\sim x}a_y\sqrt{p(x,y)}
                  \bigl(q_t(y,z)-q_t(x,z)\bigr),
\]
which satisfies $\sum_z\nu_a(z)=0$.
Also, $\|D_x(q_tf)\|^2=2\Gamma(q_tf)(x)$, so the preceding estimate gives
\[
 \begin{aligned}
 \left|\sum_z\nu_a(z)f(z)\right|^2
 &=|\langle a,D_x(q_tf)\rangle|^2\\
 &\le\|a\|^2\|D_x(q_tf)\|^2
 \le\frac{\|a\|^2}{t}\Var_{q_t(x,\cdot)} f.
 \end{aligned}
\]
Since $\nu_a$ has zero total mass and $q_t(x,\cdot)$ has full support,
variance duality and the preceding estimate give
\[
 \begin{aligned}
 0\le a^{\mathsf T}\mathbf F_t(x)a=\sum_z\frac{\nu_a(z)^2}{q_t(x,z)}=\sup_{\substack{f\in\ell^\infty(V)\\\Var_{q_t(x,\cdot)} f\le1}}
       \left|\sum_z\nu_a(z)f(z)\right|^2
 \le\frac{\|a\|^2}{t}.
 \end{aligned}
\]
The variational identity on the countable space $V$ is justified by
bounded truncations of the density
$z\mapsto\nu_a(z)/q_t(x,z)$.
Since $a$ is arbitrary, we conclude that $0\le\mathbf F_t(x)\le t^{-1}I$.
Its trace is $2J_t(x)$ by definition.

Since $\Pi_x$ is an orthogonal projection,
$0\le\Pi_x\mathbf F_t(x)\Pi_x\le t^{-1}\Pi_x$.
Thus, by \eqref{eq:sharp-projection-rank},
\[
 \begin{aligned}
 \operatorname{tr}\mathbf F_t(x)
 &=\operatorname{tr}\bigl(\Pi_x\mathbf F_t(x)\Pi_x\bigr)+b_t(x)\\
 &\le\frac{\operatorname{rank}\Pi_x}{t}+b_t(x)
 \le\frac{n_*}{t}+b_t(x),
 \end{aligned}
\]
which proves \eqref{eq:sharp-fisher-rank-split}. Moreover,
\eqref{eq:sharp-fisher-bound} gives
$0\le J_t(x)\le d/(2t)\le d/(2t_0)$ for $t\ge t_0$ and $x\in V$,
so $J_t$ is uniformly bounded on $[t_0,\infty)\times V$.

For a positive heat solution $u$, expanding $L(\Gamma u/u)$ and using
$(\partial_t-L)\Gamma u=-2\Gamma_2u$ gives
\begin{equation}\label{eq:sharp-discrete-quotient}
 \begin{aligned}
 (\partial_t-L)\frac{\Gamma u}{u}(x)
 ={}&-\frac{2\Gamma_2u(x)}{u(x)}
 -\frac{\Gamma u(x)}{u(x)^2}
       \sum_yp(x,y)\frac{(u(y)-u(x))^2}{u(y)}\\
 &+\sum_yp(x,y)
       \frac{(\Gamma u(y)-\Gamma u(x))(u(y)-u(x))}{u(x)u(y)}.
 \end{aligned}
\end{equation}
For $u(t,x)=q_t(x,z)$ and $t\ge t_0$, the first estimate in
\eqref{eq:sharp-relative-derivatives} implies
$\Gamma u(x)\le\eta_t^2u(x)^2/2$. Neighbouring values of $u$ are comparable, so
$|\Gamma u(y)-\Gamma u(x)|\le C\eta_t^2u(x)^2$. Hence the last line of
\eqref{eq:sharp-discrete-quotient} has absolute value at most
$C_1\eta_t^3u(x)$. Dropping the other nonpositive term and summing in $z$,
with termwise differentiation justified by
Lemma~\ref{lem:sharp-fisher-differentiation}, gives
\begin{equation*}
 (\partial_t-L)J_t\le-2A_t+C_1\eta_t^3.
\end{equation*}
By \eqref{eq:duhamel}, $J_T\ge0$, and $\eta_s^3\le C_2s^{-9/8}$,
\[
 2\int_{t_0}^TP_{T-s}A_s(x)\dd s
 \le P_{T-t_0}J_{t_0}(x)-J_T(x)
       +C_1\int_{t_0}^T\eta_s^3\dd s\le C_3.
\]

To estimate $b_t(x)$, fix $t\ge t_0$.
The definition of $\mathbf F_t(x)$ and the fact that $I-\Pi_x$ is an
orthogonal projection give
\[
 b_t(x)=\sum_z\frac{\|(I-\Pi_x)D_x(q_t(\cdot,z))\|^2}{q_t(x,z)}.
\]
Applying \eqref{eq:sharp-local-control} to $q_t(\cdot,z)$ for each $z$ and using the
definition of $\mathfrak q_{x,R_*}$, we obtain
\[
 b_t(x)\le C_*\sum_{v\in B(x,R_*)}\sum_z
       \frac{\Gamma_2(q_t(\cdot,z))(v)}{q_t(x,z)}.
\]
All summands on the right are nonnegative by $\CD(0,\infty)$,
so the order of summation may be interchanged.
For $v\in B(x,R_*)$, \eqref{eq:sharp-fixed-neighbourhood} gives
$q_t(v,z)\le2q_t(x,z)$. Hence
\[
 \sum_z\frac{\Gamma_2(q_t(\cdot,z))(v)}{q_t(x,z)}
 \le2\sum_z\frac{\Gamma_2(q_t(\cdot,z))(v)}{q_t(v,z)}=2A_t(v).
\]
Combining these estimates gives
\[
 b_t(x)\le2C_*\sum_{v\in B(x,R_*)}A_t(v).
\]

Uniformisation gives
$p_1(x,v)\ge\Pois(1)(R_*)\eps^{R_*}$ for $v\in B(x,R_*)$.
Since $A_t\ge0$,
\[
 \begin{aligned}
 \sum_{v\in B(x,R_*)}A_t(v)
 &\le\Pois(1)(R_*)^{-1}\eps^{-R_*}\sum_{v\in B(x,R_*)}p_1(x,v)A_t(v)\\
 &\le\Pois(1)(R_*)^{-1}\eps^{-R_*}P_1A_t(x)\le CP_1A_t(x).
 \end{aligned}
\]
We therefore have $b_t\le CP_1A_t$. Since $P_1$ commutes with
$P_{T-s}$, the first integral bound gives the second bound in
\eqref{eq:sharp-integrated-errors}.
\end{proof}

\subsection{Sharp volume growth}

\begin{proof}[Proof of Theorem~\ref{thm:sharp-volume}]
Use the fixed kernel $M$ in
\eqref{eq:sharp-fixed-mixture} and the
entropy $\mathcal H_t^M$ from Lemma~\ref{lem:entropy-volume}.
By \eqref{eq:sharp-mixture-entropy-coarse} and $m\le\kappa m_-$,
$\mathcal H_{t_0}^M-\log m_-$ is uniformly bounded, and
\eqref{eq:heat-entropy-evolution} applies to $\mathcal H_t^M$ and $q_t$.

For $y\sim x$, \eqref{eq:sharp-relative-derivatives} gives
$|q_t(y,z)/q_t(x,z)-1|\le\eta_t\le1/2$.
Writing $a_z=q_t(y,z)/q_t(x,z)-1$, we have
$\sum_zq_t(x,z)a_z=0$. Thus the Taylor estimate
$(1+a)\log(1+a)-a\le(1/2+C_1\eta_t)a^2$ for $|a|\le\eta_t$ yields
\[
 \D(q_t(y,\cdot)\|q_t(x,\cdot))
 \le(\tfrac12+C_1\eta_t)
       \sum_z\frac{(q_t(y,z)-q_t(x,z))^2}{q_t(x,z)}.
\]
Multiplying by $p(x,y)$ and summing over $y$, we use
\eqref{eq:heat-entropy-evolution} on the left and
\eqref{eq:sharp-fisher-matrix} on the right to obtain
\[
 \begin{aligned}
 (\partial_t-L)\mathcal H_t^M(x)
 &\le(\tfrac12+C_1\eta_t)\operatorname{tr}\mathbf F_t(x)\\
 &=\tfrac12\operatorname{tr}\mathbf F_t(x)
       +C_1\eta_t\operatorname{tr}\mathbf F_t(x)\\
 &\le\tfrac12\left(\frac{n_*}{t}+b_t(x)\right)
       +C_1\eta_t\frac{d}{t}\\
 &\le\frac{n_*}{2t}+\tfrac12b_t(x)+C_2\frac{\eta_t}{t}.
 \end{aligned}
\]
Here we used \eqref{eq:sharp-fisher-rank-split} for the first trace
term and \eqref{eq:sharp-fisher-bound} for the second, absorbing
$d\le\eps^{-1}$ into $C_2$.

Since $\eta_t/t\le C_3t^{-11/8}$, \eqref{eq:duhamel} and
\eqref{eq:sharp-integrated-errors} imply
\begin{equation*}
 \mathcal H_t^M(x)-\log m_-
 \le\frac{n_*}{2}\log t+C\qquad(t\ge t_0).
\end{equation*}

\Needspace{7\baselineskip}
Applying \eqref{eq:curved-volume-entropy-transfer} with $K=0$ and
$t=R^2$ now gives
\[
 \log\frac{\vol(B(o,R))}{m_-}
 \le n_*\log R+C\qquad(R\ge\sqrt{t_0}).
\]
For $1\le R<\sqrt{t_0}$, the degree bound and
$m\le\kappa m_-$ give the same conclusion after increasing $C$.
This proves \eqref{eq:sharp-volume-main} for both finite and infinite $V$.
\end{proof}

\begin{example}[Sharpness of the exponent $\lfloor d/2\rfloor$]
\label{ex:sharp-volume-scope}
For $d=2n$, take $G=\mathbb Z^n$ with
$S=\{\pm e_1,\ldots,\pm e_n\}$; for $d=2n+1\ge3$, take
$G=\mathbb Z^n\square K_2$ with
$S=\{(\pm e_i,0):1\le i\le n\}\cup\{(0,1)\}$.
Choose $a_s=a_{-s}>0$ with $\sum_{s\in S}a_s\le1$, and set
\[
 p(x,x+s)=a_s,\qquad
 p(x,x)=1-\sum_{s\in S}a_s,\qquad m\equiv m_0\in[m_-,\kappa m_-].
\]
The kernel satisfies \eqref{eq:geometry} with
$\eps=\min_{s\in S}a_s$, and its maximum support degree is $d$.

Since the translations commute, its generator $L=P-I$ satisfies
\begin{equation}
      \label{equ:square-formula}
       \Gamma_2f(x)=\frac14\sum_{s,t\in S}a_sa_t
  \bigl(f(x+s+t)-f(x+s)-f(x+t)+f(x)\bigr)^2\ge0.
\end{equation}
The element $(0,1)=-(0,1)$ is counted once in $S$. Hence $\CD(0,\infty)$ holds.

The positive weights
leave the support graph unchanged, so
$\vol(B(o,R))=m_0|B(o,R)|\asymp m_0R^n$ for $R\ge1$.
Thus the exponent $n=\lfloor d/2\rfloor$ cannot be lowered in the
weighted setting.
\end{example}

\section{Proof of the gradient estimate}
\label{sec:moment}
\label{sec:signed-moment}

In the next three lemmas, assume \eqref{eq:geometry} and
$\CD(K,\infty)$, and normalize the measure so that $1\le m\le\kappa$.
Let $h>0$ belong to $\ell^1(m)$, let $\beta\ge0$, and use $a,M,S$
as in Proposition~\ref{prop:moment}. Put
\[
 G=\Gamma(\log h),\quad \mu_x=m(x)h_x,\quad
 \delta_*=\sup_{x\sim y}|\log h_x-\log h_y|,\quad B_2=\sum_x\mu_xG_x^2.
\]
Assume $\delta_*\le1/4$. Then
\[
 G\le\delta_*^2/2,\qquad |a|\le\beta+e^{\delta_*}-1.
\]
The same edge bound controls $L(\log h)$, $\Gamma_2(\log h)$, and all local
differences below. Thus every pairing converges absolutely, since
$h\in\ell^1(m)$, and integration by parts is justified even when
$\log h$ is unbounded.

Write
\[
 \begin{gathered}
 W_{xy}=m(x)p(x,y)(h_x+h_y),\qquad
 E=\sum_{\{x,y\}}W_{xy}(G_x-G_y)^2,\\
 X=M^{2/3}S^{1/3},\qquad Y=M^{1/3}S^{2/3}.
 \end{gathered}
\]

Orient each edge
$e=\{x,y\}$ so that $h_x\ge h_y$, and set
\[
 r=\log h_x-\log h_y\ge0,\qquad
 \theta=m(x)p(x,y)(h_x-h_y),\qquad W_e=W_{xy},\qquad
 \Delta F=F_x-F_y.
\]

By $\CD(K,\infty)$ and $K+K_-\ge0$, the form
$\mathcal Q_{K_-}=\Gamma_2+K_-\Gamma$ is positive semidefinite. Set
\[
 \mathcal B_{K_-}=\sum_x\mu_xG_x\mathcal Q_{K_-}(\log h)(x)\ge0.
\]

\Needspace{13\baselineskip}
\begin{lemma}
\label{lem:signed-elementary}
With the notation above, one has
\begin{equation}\label{eq:signed-holder}
 \begin{gathered}
 \beta^3\norm{h}_1\le M,\qquad
 \sum_x\mu_xa_x^2G_x\le X,\qquad
 \sum_x\mu_x|a_x|G_x^2\le Y,\\
 \beta B_2\le Y,\qquad
 \beta\sum_x\mu_x|a_x|G_x\le X,\qquad
 B_2\le\norm{h}_1^{1/3}S^{2/3}.
 \end{gathered}
\end{equation}
Moreover, for every bounded real function $F$,
\begin{equation}\label{eq:signed-current}
 \begin{split}
 \sum_e\theta\Delta F&=\sum_x\mu_x(a_x-\beta)F_x,
 \qquad \theta=W_e\tanh(r/2)\le W_er/2,\\
 L(\log h)&=\beta-a-G-T_h,\qquad
 T_h=\Psi(\log h)-G,\qquad |T_h|\le\delta_* G/2.
 \end{split}
\end{equation}
\begin{equation}\label{eq:vertex-edge-comparison}
 2\sum_x\mu_xG_xF_x^2
 \le\sum_eW_er^2(F_x^2+F_y^2)
 \le\frac{14}{3}\sum_x\mu_xG_xF_x^2.
\end{equation}
\end{lemma}
\begin{proof}
Conservativity gives $\sum_x\mu_xa_x=\beta\norm{h}_1$, so
H\"older's inequality yields \eqref{eq:signed-holder}.
For the penultimate bound, use
$\sum\mu|a|G\le M^{1/3}S^{1/3}\norm{h}_1^{1/3}$ and
$\beta\norm{h}_1^{1/3}\le M^{1/3}$.

Integration by parts and Taylor's formula give \eqref{eq:signed-current}.
Grouping edge terms at their endpoints gives
\eqref{eq:vertex-edge-comparison}, because
$h_x\le h_x+h_y\le(7/3)h_x$ for either endpoint $x$ of an edge.
\end{proof}

\begin{lemma}
\label{lem:signed-weighted-bochner}
With the notation above,
\begin{equation}\label{eq:signed-bochner-bound}
 \mathcal B_{K_-}+E
 \le128(X+Y+\delta_*^2S+{K_-} B_2).
\end{equation}
\end{lemma}
\begin{proof}
Set
\[
 Q_F=\sum_e\theta(G_x+G_y)\Delta F.
\]
Define $R_F$ by
\[
 \sum_e m(x)p(x,y)r(h_xG_x+h_yG_y)\Delta F=Q_F+R_F.
\]
Since $0\le W_er/2-\theta\le W_er^3/24$, one has
\begin{equation}\label{eq:edge-remainder}
 |R_F|\le\sum_eW_e
 \left(\frac{r^3}{24}(G_x+G_y)
       +\frac{r^2}{4}|\Delta G|\right)|\Delta F|.
\end{equation}
Indeed, the difference of the two edge coefficients is
$(W_er/2-\theta)(G_x+G_y)+(\theta r/2)\Delta G$.
By \eqref{eq:signed-current},
\begin{equation}\label{eq:Q-identity}
 Q_F=2\sum_x\mu_x(a_x-\beta)G_xF_x
       -\sum_e\theta(F_x+F_y)\Delta G,
 \qquad Q_G=\sum_x\mu_x(a_x-\beta)G_x^2.
\end{equation}

Multiply $2\Gamma_2(\log h)=LG+2\Gamma(\log h,a+G+T_h)$ by $\mu G$ and sum.
Using
\[
 \sum_x\mu_xG_xLG_x
 =-\tfrac12E-\tfrac12\sum_x\mu_x(a_x-\beta)G_x^2
\]
gives
\begin{equation}\label{eq:signed-bochner-identity}
 \begin{split}
 2\mathcal B_{K_-}+\tfrac12E
 ={}&2{K_-} B_2+\tfrac12\sum_x\mu_x(a_x-\beta)G_x^2\\
    &+Q_a+Q_{T_h}+R_a+R_G+R_{T_h}.
 \end{split}
\end{equation}

The vertex term is at most $Y/2$. The bounds in
Lemma~\ref{lem:signed-elementary}, together with
\eqref{eq:edge-remainder}, \eqref{eq:Q-identity}, and Cauchy--Schwarz,
give
\[
 Q_a\le4X+\sqrt{7/3}\sqrt{EX},\qquad
 |R_a|\le\frac{7\delta_*}{18}Y
                  +\frac{\delta_*}{2}\sqrt{7/3}\sqrt{EX}.
\]
Here one uses $|\Delta a|\le|a_x|+|a_y|$ and
$\sum\mu(a^2+\beta|a|)G\le2X$.
Since $\delta_*\le1/4$, Young's inequality implies
\[
 Q_a+R_a\le\tfrac18E+12X+\tfrac{\delta_*}{2}Y.
\]

The same estimates, with $|T_h|\le\delta_* G/2$, yield
\[
 \begin{aligned}
 |Q_{T_h}|&\le2\delta_* Y+\frac{\delta_*}{4}\sqrt{\frac{28}{3}ES},\\
 |R_G|&\le\frac{\delta_*^2}{4}E
                  +\frac{\delta_*^2}{24}\sqrt{\frac{28}{3}ES},\\
 |R_{T_h}|&\le\frac{7\delta_*^2}{36}S
                  +\frac{\delta_*^2}{8}\sqrt{\frac{28}{3}ES}.
 \end{aligned}
\]
For example, the first follows from \eqref{eq:Q-identity} and
$\sum\mu(|a|+\beta)G^2\le2Y$; the others follow from
\eqref{eq:edge-remainder} and \eqref{eq:vertex-edge-comparison}.
Their sum is at most $E/4+2\delta_* Y+3\delta_*^2S$, by Young's
inequality. Substitution into \eqref{eq:signed-bochner-identity} gives
\[
 2\mathcal B_{K_-}+\tfrac18E
 \le12X+(\tfrac12+\tfrac{5\delta_*}{2})Y
       +3\delta_*^2S+2{K_-} B_2.
\]
This proves \eqref{eq:signed-bochner-bound}.
\end{proof}

\begin{lemma}
\label{lem:signed-pointmass-trace}
With the notation above, assume in addition that $K_-\le1$ and
$\delta_*\le\eps/(12\kappa)$. Then, with $c=\eps/(3\kappa)$,
\begin{equation}\label{eq:signed-trace-bound}
 S\le\frac{64}{c^2}(\mathcal B_{K_-}+E)+\frac{14}{c}Y
 \le C_T(\mathcal B_{K_-}+E+Y),
 \qquad C_T=1024\kappa^2\eps^{-2}.
\end{equation}
\end{lemma}
\begin{proof}
Set
\[
 c_{xy}=p(x,y)p(y,x)\quad(x\ne y),\quad c_{xx}=0,
 \quad \mathcal A(\log h)(x)=\sum_yc_{xy}(\log h_y-\log h_x),\quad
 \mathcal H=hG^2.
\]
The defining formulas for $\Gamma$ and $\Gamma_2$ give
\begin{equation}\label{eq:signed-pointmass}
 \begin{split}
 4\mathcal Q_{K_-}(\log h,\1_{\{x\}})(x)
   &=L^2(\log h)(x)-2\mathcal A(\log h)(x)-2{K_-} L(\log h)(x),\\
 \mathcal Q_{K_-}(\1_{\{x\}})(x)
   &=\tfrac14(1-p(x,x))^2+\tfrac34\sum_{y\ne x}c_{xy}
                               +\tfrac{K_-}{2}(1-p(x,x))\le\tfrac32.
 \end{split}
\end{equation}
Positive semidefiniteness therefore implies
$|\mathcal Q_{K_-}(\log h,\1_{\{x\}})(x)|^2
 \le(3/2)\mathcal Q_{K_-}(\log h)(x)$.

Multiply the first identity by $\mathcal H_x$, sum, and use
$L^2(\log h)=-L(a+G+T_h)$ to obtain
\begin{equation}\label{eq:signed-trace-identity}
 \begin{split}
 2\cE_c(\mathcal H,\log h)
 ={}&4\sum_x\mathcal H_x\mathcal Q_{K_-}(\log h,\1_{\{x\}})(x)
        +2{K_-}\sum_x\mathcal H_xL(\log h)(x)\\
    &-\cE(\mathcal H/m,a+G+T_h).
 \end{split}
\end{equation}
Here $\cE_c$ is the edge form with weights $c_{xy}$ and counting
measure. Cauchy--Schwarz and $m\ge1$ bound the first term by
$4\sqrt{(3/2)\mathcal B_{K_-} S}$. Since
$L(\log h)=\beta-a-\Psi(\log h)\le\beta+|a|$, the curvature term is at most
$2{K_-}(\beta B_2+Y)\le4Y$.

On each edge,
\[
 2c_{xy}r\Delta\mathcal H
 =c_{xy}r(h_x+h_y)(G_x+G_y)\Delta G
  +c_{xy}r(h_x-h_y)(G_x^2+G_y^2).
\]
Using $\eps/\kappa\le c_{xy}/(m(x)p(x,y))\le1$,
$\tanh(r/2)\ge r/3$, and \eqref{eq:vertex-edge-comparison}, the
second term sums to at least $2\eps S/(3\kappa)$, whereas the
absolute sum of the first terms is at most
$\sqrt{(28/3)ES}\le4\sqrt{ES}$. Thus
\begin{equation}\label{eq:signed-positive-trace}
 \frac{2\eps}{3\kappa}S
 \le2\cE_c(\mathcal H,\log h)+4\sqrt{ES}.
\end{equation}

For the remaining energies in \eqref{eq:signed-trace-identity},
$|\Delta(\mathcal H/m)|\le h_xG_x^2+h_yG_y^2$ gives
\[
 |\cE(\mathcal H/m,F)|
 \le\sum_{x,y}\mu_xp(x,y)G_x^2|F_x-F_y|.
\]
Reversibility and $h_x/h_y\le4/3$ imply
$\sum_{x,y}\mu_xp(x,y)|a_y|^3\le(4/3)M$.
H\"older's inequality, $G\le\delta_*^2/2$, and
$|T_h|\le\delta_* G/2$ consequently yield
\begin{equation}\label{eq:signed-trace-remainders}
 \begin{split}
 |\cE(\mathcal H/m,a)|&\le3Y,\\
 |\cE(\mathcal H/m,G)|&\le\delta_*\sqrt{ES}/\sqrt2,\\
 |\cE(\mathcal H/m,T_h)|&\le2\delta_* S.
 \end{split}
\end{equation}
For the middle inequality, apply Cauchy--Schwarz on edges and use
\[
 \frac{(h_xG_x^2+h_yG_y^2)^2}{h_x+h_y}
 \le h_xG_x^4+h_yG_y^4,
 \qquad \sum_x\mu_x(1-p(x,x))G_x^4\le\delta_*^2S/2.
\]

Combining the preceding bounds gives
\[
 \frac{2\eps}{3\kappa}S
 \le4\sqrt{\tfrac32\mathcal B_{K_-} S}
       +(4+\delta_*/\sqrt2)\sqrt{ES}+7Y+2\delta_* S.
\]
Because $\delta_*\le\eps/(12\kappa)$, it follows that
\[
 \frac\eps{3\kappa}S
 \le8\sqrt{(\mathcal B_{K_-}+E)S}+7Y.
\]
Young's inequality gives \eqref{eq:signed-trace-bound}.
\end{proof}

\begin{proof}[Proof of Proposition~\ref{prop:moment}]
Normalize the measure so that $1\le m\le\kappa$; both sides of
\eqref{eq:moment} and \eqref{eq:source-flux} scale by the same factor.
Take $\delta=\eps/(1024\kappa)$
and use the notation introduced above, so that
$\delta_*\le\delta\le1/4$.

If $K_-\ge1$, then
$S\le\delta_*^6\norm{h}_1/8\le K_-^3\norm{h}_1$, which proves
\eqref{eq:moment}.

For $K_-\le1$, the additional assumptions of
Lemma~\ref{lem:signed-pointmass-trace} hold.
Since $\delta_*\le\delta$, our choice of $\delta$ gives
$128C_T\delta_*^2\le128C_T\delta^2=1/8$. Combining Lemmas~\ref{lem:signed-weighted-bochner} and
\ref{lem:signed-pointmass-trace}, absorbing $\delta_*^2S$, and using
Lemma~\ref{lem:signed-elementary}, we obtain
\[
 S\le C\bigl(M^{2/3}S^{1/3}+M^{1/3}S^{2/3}
                    +{K_-}\norm{h}_1^{1/3}S^{2/3}\bigr).
\]
Young's inequality bounds the right-hand side by
$S/2+C'(M+{K_-}^3\norm{h}_1)$, proving \eqref{eq:moment} also in this case.

For every $K\in\R$ and bounded $1$-Lipschitz $f$, the identities
$\sum_eW_e\le\norm{h}_1$ and $\sum_eW_er^2=2\sum_x\mu_xG_x$ give
\[
 \begin{aligned}
 |\cE(h,f)|
 &\le\frac12\sum_eW_er
 \le\frac1{\sqrt2}\norm{h}_1^{1/2}\Bigl(\sum_x\mu_xG_x\Bigr)^{1/2}\\
 &\le\frac1{\sqrt2}\norm{h}_1^{5/6}S^{1/6}\\
 &\le C\norm{h}_1^{5/6}(M+{K_-}^3\norm{h}_1)^{1/6}.
 \end{aligned}
\]
This proves \eqref{eq:source-flux}.
\end{proof}

\section*{Acknowledgments}
This work was supported by the Scientific Research Innovation Capability Support Project for Young Faculty (Grant No. SRICSPYF-ZY2025160) and the National Natural Science Foundation of China (Grant No. 12431004).

\IfFileExists{references.bib}{%
  \bibliographystyle{amsplain}
  \bibliography{references}

@article {MunchRemarks,
    AUTHOR = {M\"unch, Florentin},
     TITLE = {Remarks on curvature dimension conditions on graphs},
   JOURNAL = {Calc. Var. Partial Differential Equations},
  FJOURNAL = {Calculus of Variations and Partial Differential Equations},
    VOLUME = {56},
      YEAR = {2017},
    NUMBER = {1},
     PAGES = {Paper No. 11, 8},
}

@misc{Zhang2026CDpsiDoubling,
      title={Gradient estimates and volume doubling for locally finite weighted graphs with {$CD\psi(n,-K)$} condition}, 
      author={Qianwei Zhang},
      year={2026},
      eprint={2608.30409},
      note={arXiv:2608.30409},
      archivePrefix={arXiv},
      primaryClass={math.DG},
      url={https://arxiv.org/abs/2608.30409}, 
}

@article {WeberZacherJFA,
    AUTHOR = {Weber, Frederic and Zacher, Rico},
     TITLE = {The entropy method under curvature-dimension conditions in the
              spirit of {B}akry-\'Emery in the discrete setting of {M}arkov
              chains},
   JOURNAL = {J. Funct. Anal.},
  FJOURNAL = {Journal of Functional Analysis},
    VOLUME = {281},
      YEAR = {2021},
    NUMBER = {5},
     PAGES = {Paper No. 109061, 81},
      ISSN = {0022-1236,1096-0783},
   MRCLASS = {60J27 (47D07 49Q22 53C21)},
  MRNUMBER = {4253929},
       DOI = {10.1016/j.jfa.2021.109061},
       URL = {https://doi.org/10.1016/j.jfa.2021.109061},
}

@article {RocknerWang2010LogHarnack,
    AUTHOR = {R\"ockner, Michael and Wang, Feng-Yu},
     TITLE = {Log-{H}arnack inequality for stochastic differential equations
              in {H}ilbert spaces and its consequences},
   JOURNAL = {Infin. Dimens. Anal. Quantum Probab. Relat. Top.},
  FJOURNAL = {Infinite Dimensional Analysis, Quantum Probability and Related
              Topics},
    VOLUME = {13},
      YEAR = {2010},
    NUMBER = {1},
     PAGES = {27--37},
      ISSN = {0219-0257,1793-6306},
   MRCLASS = {60J60 (58J65 60H15)},
  MRNUMBER = {2646789},
MRREVIEWER = {Cristina\ Zucca},
       DOI = {10.1142/S0219025710003936},
       URL = {https://doi.org/10.1142/S0219025710003936},
}

@article {BakryGentilLedoux2015Harnack,
    AUTHOR = {Bakry, Dominique and Gentil, Ivan and Ledoux, Michel},
     TITLE = {On {H}arnack inequalities and optimal transportation},
   JOURNAL = {Ann. Sc. Norm. Super. Pisa Cl. Sci. (5)},
  FJOURNAL = {Annali della Scuola Normale Superiore di Pisa. Classe di
              Scienze. Serie V},
    VOLUME = {14},
      YEAR = {2015},
    NUMBER = {3},
     PAGES = {705--727},
      ISSN = {0391-173X,2036-2145},
   MRCLASS = {58J60 (35B45 35B65 35K08 35R01 49Q20 58J35)},
  MRNUMBER = {3445198},
MRREVIEWER = {Fabrice\ Baudoin},
}

@article {Cheng1975,
    AUTHOR = {Cheng, Shiu Yuen},
     TITLE = {Eigenvalue comparison theorems and its geometric applications},
   JOURNAL = {Math. Z.},
  FJOURNAL = {Mathematische Zeitschrift},
    VOLUME = {143},
      YEAR = {1975},
    NUMBER = {3},
     PAGES = {289--297},
      ISSN = {0025-5874,1432-1823},
   MRCLASS = {58G99 (53C20 58E15)},
  MRNUMBER = {378001},
MRREVIEWER = {J.\ R.\ Vanstone},
       DOI = {10.1007/BF01214381},
       URL = {https://doi.org/10.1007/BF01214381},
}

@article {Delmotte,
    AUTHOR = {Delmotte, Thierry},
     TITLE = {Parabolic {H}arnack inequality and estimates of {M}arkov
              chains on graphs},
   JOURNAL = {Rev. Mat. Iberoamericana},
  FJOURNAL = {Revista Matem\'atica Iberoamericana},
    VOLUME = {15},
      YEAR = {1999},
    NUMBER = {1},
     PAGES = {181--232},
      ISSN = {0213-2230},
   MRCLASS = {35K10 (05C99 35B45)},
  MRNUMBER = {1681641},
MRREVIEWER = {P.\ Lindqvist},
       DOI = {10.4171/RMI/254},
       URL = {https://doi.org/10.4171/RMI/254},
}

@misc{GuoHuangHuang,
      title={Nonnegative {Bakry--{\'E}mery} curvature on bounded-degree graphs implies volume doubling and {Poincar\'e} inequalities}, 
      author={Qi Guo and Xueping Huang and Yi C. Huang},
      year={2026},
      eprint={2607.15522},
      note={arXiv:2607.15522},
      archivePrefix={arXiv},
      primaryClass={math.DG},
      url={https://arxiv.org/abs/2607.15522}, 
}

@misc{RussPajot,
      title={Infinite graphs satisfying the {Bakry--{\'E}mery} curvature condition {$CD(0,n)$}: The modified heat equation and applications to geometric analysis}, 
      author={Emmanuel Russ and Hervé Pajot},
      year={2025},
      eprint={2507.19235},
      note={arXiv:2507.19235},
      archivePrefix={arXiv},
      primaryClass={math.DG},
      url={https://arxiv.org/abs/2507.19235}, 
}

@article {KellerLenz2012,
    AUTHOR = {Keller, Matthias and Lenz, Daniel},
     TITLE = {Dirichlet forms and stochastic completeness of graphs and
              subgraphs},
   JOURNAL = {J. Reine Angew. Math.},
  FJOURNAL = {Journal f\"ur die Reine und Angewandte Mathematik. [Crelle's
              Journal]},
    VOLUME = {666},
      YEAR = {2012},
     PAGES = {189--223},
      ISSN = {0075-4102,1435-5345},
   MRCLASS = {60J45 (05C22 31C20)},
  MRNUMBER = {2920886},
MRREVIEWER = {Fr\'ed\'eric\ Math\'eus},
       DOI = {10.1515/CRELLE.2011.122},
       URL = {https://doi.org/10.1515/CRELLE.2011.122},
}

@article {LinLiu2018,
    AUTHOR = {Lin, Yong and Liu, Shuang},
     TITLE = {Equivalent properties of {CD} inequalities on graphs},
   JOURNAL = {Acta Math. Sinica (Chinese Ser.)},
  FJOURNAL = {Acta Mathematica Sinica. Chinese Series},
    VOLUME = {61},
      YEAR = {2018},
    NUMBER = {3},
     PAGES = {431--440},
      ISSN = {0583-1431},
   MRCLASS = {58J35 (53C23)},
  MRNUMBER = {4545901},
}

@article {Wojciechowski2009,
    AUTHOR = {Wojciechowski, Rados\l aw K.},
     TITLE = {Heat kernel and essential spectrum of infinite graphs},
   JOURNAL = {Indiana Univ. Math. J.},
  FJOURNAL = {Indiana University Mathematics Journal},
    VOLUME = {58},
      YEAR = {2009},
    NUMBER = {3},
     PAGES = {1419--1441},
      ISSN = {0022-2518,1943-5258},
   MRCLASS = {35K08 (05C50 35P05)},
  MRNUMBER = {2542093},
       DOI = {10.1512/iumj.2009.58.3575},
       URL = {https://doi.org/10.1512/iumj.2009.58.3575},
}

@article {HuaLin2017,
    AUTHOR = {Hua, Bobo and Lin, Yong},
     TITLE = {Stochastic completeness for graphs with curvature dimension
              conditions},
   JOURNAL = {Adv. Math.},
  FJOURNAL = {Advances in Mathematics},
    VOLUME = {306},
      YEAR = {2017},
     PAGES = {279--302},
      ISSN = {0001-8708,1090-2082},
   MRCLASS = {58J65 (05C63)},
  MRNUMBER = {3581303},
MRREVIEWER = {Alexander\ G.\ Losev},
       DOI = {10.1016/j.aim.2016.10.022},
       URL = {https://doi.org/10.1016/j.aim.2016.10.022},
}

@article {KellerMunch2024,
    AUTHOR = {Keller, Matthias and M\"unch, Florentin},
     TITLE = {Gradient estimates, {B}akry-{E}mery {R}icci curvature and
              ellipticity for unbounded graph {L}aplacians},
   JOURNAL = {Comm. Anal. Geom.},
  FJOURNAL = {Communications in Analysis and Geometry},
    VOLUME = {32},
      YEAR = {2024},
    NUMBER = {2},
     PAGES = {343--364},
      ISSN = {1019-8385,1944-9992},
   MRCLASS = {53C21 (05C63 58J65)},
  MRNUMBER = {4828188},
       DOI = {10.4310/cag.241015004626},
       URL = {https://doi.org/10.4310/cag.241015004626},
}

@book {CoverThomas2006,
    AUTHOR = {Cover, Thomas M. and Thomas, Joy A.},
     TITLE = {Elements of information theory},
   EDITION = {Second},
 PUBLISHER = {Wiley-Interscience [John Wiley \& Sons], Hoboken, NJ},
      YEAR = {2006},
     PAGES = {xxiv+748},
      ISBN = {978-0-471-24195-9; 0-471-24195-4},
   MRCLASS = {00-01 (94-01)},
  MRNUMBER = {2239987},
}

@book {EngelNagel2000,
    AUTHOR = {Engel, Klaus-Jochen and Nagel, Rainer},
     TITLE = {One-parameter semigroups for linear evolution equations},
    SERIES = {Graduate Texts in Mathematics},
    VOLUME = {194},
 PUBLISHER = {Springer-Verlag, New York},
      YEAR = {2000},
     PAGES = {xxii+586},
      ISBN = {0-387-98463-1},
   MRCLASS = {47D06 (34G10 35K90 47N20)},
  MRNUMBER = {1721989},
MRREVIEWER = {Charles\ Batty},
}

@book {PolyanskiyWu2025,
    AUTHOR = {Polyanskiy, Yury and Wu, Yihong},
     TITLE = {Information theory: from coding to learning},
 PUBLISHER = {Cambridge University Press},
      YEAR = {2025},
      ISBN = {978-1-108-83290-8},
       DOI = {10.1017/9781108966351},
       URL = {https://doi.org/10.1017/9781108966351},
}

@article {AckermansVanEijndhovenMartens1983,
    AUTHOR = {Ackermans, S. T. M. and van Eijndhoven, S. J. L. and Martens,
              F. J. L.},
     TITLE = {On almost commuting operators},
   JOURNAL = {Nederl. Akad. Wetensch. Indag. Math.},
  FJOURNAL = {Koninklijke Nederlandse Akademie van Wetenschappen.
              Indagationes Mathematicae},
    VOLUME = {45},
      YEAR = {1983},
    NUMBER = {4},
     PAGES = {385--391},
      ISSN = {0019-3577},
   MRCLASS = {47B47 (47B15)},
  MRNUMBER = {731821},
MRREVIEWER = {Takayuki\ Furuta},
}

@article {CushingKamtueLiuPeyerimhoff2022,
    AUTHOR = {Cushing, David and Kamtue, Supanat and Liu, Shiping and
              Peyerimhoff, Norbert},
     TITLE = {{Bakry--\'Emery} curvature on graphs as an eigenvalue problem},
   JOURNAL = {Calc. Var. Partial Differential Equations},
  FJOURNAL = {Calculus of Variations and Partial Differential Equations},
    VOLUME = {61},
      YEAR = {2022},
    NUMBER = {2},
     PAGES = {Paper No. 62, 33},
      ISSN = {0944-2669,1432-0835},
   MRCLASS = {51K10 (05C50 05C76 53C23)},
  MRNUMBER = {4376541},
       DOI = {10.1007/s00526-021-02179-z},
       URL = {https://doi.org/10.1007/s00526-021-02179-z},
}

@book {SchoenYau1994Lectures,
    AUTHOR = {Schoen, R. and Yau, S.-T.},
     TITLE = {Lectures on differential geometry},
    SERIES = {Conference Proceedings and Lecture Notes in Geometry and
              Topology},
    VOLUME = {I},
      NOTE = {Lecture notes prepared by Wei Yue Ding, Kung Ching Chang [Gong
              Qing Zhang], Jia Qing Zhong and Yi Chao Xu,
              Translated from the Chinese by Ding and S. Y. Cheng,
              With a preface translated from the Chinese by Kaising Tso},
 PUBLISHER = {International Press, Cambridge, MA},
      YEAR = {1994},
     PAGES = {v+235},
      ISBN = {1-57146-012-8},
   MRCLASS = {53-01 (53-02 53C21 58G30)},
  MRNUMBER = {1333601},
MRREVIEWER = {Man\ Chun\ Leung},
}

@book {Petersen2016,
    AUTHOR = {Petersen, Peter},
     TITLE = {Riemannian geometry},
    SERIES = {Graduate Texts in Mathematics},
    VOLUME = {171},
   EDITION = {Third},
 PUBLISHER = {Springer, Cham},
      YEAR = {2016},
     PAGES = {xviii+499},
      ISBN = {978-3-319-26652-7; 978-3-319-26654-1},
   MRCLASS = {53-01 (53C20 53C21 53C23)},
  MRNUMBER = {3469435},
       DOI = {10.1007/978-3-319-26654-1},
       URL = {https://doi.org/10.1007/978-3-319-26654-1},
}

@article {LiYau1986,
    AUTHOR = {Li, Peter and Yau, Shing-Tung},
     TITLE = {On the parabolic kernel of the {S}chr\"odinger operator},
   JOURNAL = {Acta Math.},
  FJOURNAL = {Acta Mathematica},
    VOLUME = {156},
      YEAR = {1986},
    NUMBER = {3-4},
     PAGES = {153--201},
      ISSN = {0001-5962,1871-2509},
   MRCLASS = {58G11 (35J10)},
  MRNUMBER = {834612},
MRREVIEWER = {Harold\ Donnelly},
       DOI = {10.1007/BF02399203},
       URL = {https://doi.org/10.1007/BF02399203},
}

@article {LiXu2011DifferentialHarnack,
    AUTHOR = {Li, Junfang and Xu, Xiangjin},
     TITLE = {Differential {H}arnack inequalities on {R}iemannian manifolds
              {I}: linear heat equation},
   JOURNAL = {Adv. Math.},
  FJOURNAL = {Advances in Mathematics},
    VOLUME = {226},
      YEAR = {2011},
    NUMBER = {5},
     PAGES = {4456--4491},
      ISSN = {0001-8708,1090-2082},
   MRCLASS = {53C21 (58J35)},
  MRNUMBER = {2770456},
MRREVIEWER = {Jia-Yong\ Wu},
       DOI = {10.1016/j.aim.2010.12.009},
       URL = {https://doi.org/10.1016/j.aim.2010.12.009},
}

@article {BakryBolleyGentil2017LiYau,
    AUTHOR = {Bakry, Dominique and Bolley, Fran\c{c}ois and Gentil, Ivan},
     TITLE = {The {L}i-{Y}au inequality and applications under a
              curvature-dimension condition},
   JOURNAL = {Ann. Inst. Fourier (Grenoble)},
  FJOURNAL = {Universit\'e{} de Grenoble. Annales de l'Institut Fourier},
    VOLUME = {67},
      YEAR = {2017},
    NUMBER = {1},
     PAGES = {397--421},
      ISSN = {0373-0956,1777-5310},
   MRCLASS = {58J35 (46L53 58J65 60J25)},
  MRNUMBER = {3612336},
MRREVIEWER = {Fabrice\ Baudoin},
       DOI = {10.5802/aif.3086},
       URL = {https://doi.org/10.5802/aif.3086},
}

@article {Grigoryan1992,
    AUTHOR = {Grigor'yan, A. A.},
     TITLE = {The heat equation on noncompact {R}iemannian manifolds},
   JOURNAL = {Mat. Sb.},
  FJOURNAL = {Matematicheski\u i\ Sbornik},
    VOLUME = {182},
      YEAR = {1991},
    NUMBER = {1},
     PAGES = {55--87},
      ISSN = {0368-8666},
   MRCLASS = {58G11 (35K05)},
  MRNUMBER = {1098839},
MRREVIEWER = {J\'ozef\ Dodziuk},
}

@article {SaloffCoste1992,
    AUTHOR = {Saloff-Coste, L.},
     TITLE = {A note on {P}oincar\'e, {S}obolev, and {H}arnack inequalities},
   JOURNAL = {Internat. Math. Res. Notices},
  FJOURNAL = {International Mathematics Research Notices},
      YEAR = {1992},
    NUMBER = {2},
     PAGES = {27--38},
      ISSN = {1073-7928,1687-0247},
   MRCLASS = {58G11 (31C12 35B05)},
  MRNUMBER = {1150597},
MRREVIEWER = {Thierry\ Coulhon},
       DOI = {10.1155/S1073792892000047},
       URL = {https://doi.org/10.1155/S1073792892000047},
}

@incollection {BakryEmery1985,
    AUTHOR = {Bakry, D. and \'Emery, Michel},
     TITLE = {Diffusions hypercontractives},
 BOOKTITLE = {S\'eminaire de probabilit\'es, {XIX}, 1983/84},
    SERIES = {Lecture Notes in Math.},
    VOLUME = {1123},
     PAGES = {177--206},
 PUBLISHER = {Springer, Berlin},
      YEAR = {1985},
      ISBN = {3-540-15230-X},
   MRCLASS = {60J60 (58C40 58G32)},
  MRNUMBER = {889476},
MRREVIEWER = {Jacques\ Vauthier},
       DOI = {10.1007/BFb0075847},
       URL = {https://doi.org/10.1007/BFb0075847},
}

@article {LinYau2010,
    AUTHOR = {Lin, Yong and Yau, Shing-Tung},
     TITLE = {Ricci curvature and eigenvalue estimate on locally finite
              graphs},
   JOURNAL = {Math. Res. Lett.},
  FJOURNAL = {Mathematical Research Letters},
    VOLUME = {17},
      YEAR = {2010},
    NUMBER = {2},
     PAGES = {343--356},
      ISSN = {1073-2780},
   MRCLASS = {05C10 (05C50)},
  MRNUMBER = {2644381},
       DOI = {10.4310/MRL.2010.v17.n2.a13},
       URL = {https://doi.org/10.4310/MRL.2010.v17.n2.a13},
}

@article {BauerEtAl2015,
    AUTHOR = {Bauer, Frank and Horn, Paul and Lin, Yong and Lippner, Gabor
              and Mangoubi, Dan and Yau, Shing-Tung},
     TITLE = {Li-{Y}au inequality on graphs},
   JOURNAL = {J. Differential Geom.},
  FJOURNAL = {Journal of Differential Geometry},
    VOLUME = {99},
      YEAR = {2015},
    NUMBER = {3},
     PAGES = {359--405},
      ISSN = {0022-040X,1945-743X},
   MRCLASS = {35R02 (35A23 58J35)},
  MRNUMBER = {3316971},
MRREVIEWER = {Thierry\ Coulhon},
       DOI = {10.4310/jdg/1424880980},
       URL = {http://projecteuclid.org/euclid.jdg/1424880980},
}

@article {DierKassmannZacher2021,
    AUTHOR = {Dier, Dominik and Kassmann, Moritz and Zacher, Rico},
     TITLE = {Discrete versions of the {L}i-{Y}au gradient estimate},
   JOURNAL = {Ann. Sc. Norm. Super. Pisa Cl. Sci. (5)},
  FJOURNAL = {Annali della Scuola Normale Superiore di Pisa. Classe di
              Scienze. Serie V},
    VOLUME = {22},
      YEAR = {2021},
    NUMBER = {2},
     PAGES = {691--744},
      ISSN = {0391-173X,2036-2145},
   MRCLASS = {35K05 (35K10 35R02)},
  MRNUMBER = {4288669},
}

@article {WeberZacher2023HeatKernel,
    AUTHOR = {Weber, Frederic and Zacher, Rico},
     TITLE = {Li-{Y}au inequalities for general non-local diffusion
              equations via reduction to the heat kernel},
   JOURNAL = {Math. Ann.},
  FJOURNAL = {Mathematische Annalen},
    VOLUME = {385},
      YEAR = {2023},
    NUMBER = {1-2},
     PAGES = {393--419},
      ISSN = {0025-5831,1432-1807},
   MRCLASS = {58J35 (35K08 35R11 60J27)},
  MRNUMBER = {4542719},
MRREVIEWER = {Feng-Yu\ Wang},
       DOI = {10.1007/s00208-021-02350-z},
       URL = {https://doi.org/10.1007/s00208-021-02350-z},
}

@article {Munch2018Nonlinear,
    AUTHOR = {M\"unch, Florentin},
     TITLE = {Li-{Y}au inequality on finite graphs via non-linear curvature
              dimension conditions},
   JOURNAL = {J. Math. Pures Appl. (9)},
  FJOURNAL = {Journal de Math\'ematiques Pures et Appliqu\'ees. Neuvi\`eme
              S\'erie},
    VOLUME = {120},
      YEAR = {2018},
     PAGES = {130--164},
      ISSN = {0021-7824,1776-3371},
   MRCLASS = {35R02 (05C10 35K05 53C21 58J35)},
  MRNUMBER = {3906157},
MRREVIEWER = {Thierry\ Coulhon},
       DOI = {10.1016/j.matpur.2018.10.006},
       URL = {https://doi.org/10.1016/j.matpur.2018.10.006},
}

@article {HornLinLiuYau2019,
    AUTHOR = {Horn, Paul and Lin, Yong and Liu, Shuang and Yau, Shing-Tung},
     TITLE = {Volume doubling, {P}oincar\'e{} inequality and {G}aussian heat
              kernel estimate for non-negatively curved graphs},
   JOURNAL = {J. Reine Angew. Math.},
  FJOURNAL = {Journal f\"ur die Reine und Angewandte Mathematik. [Crelle's
              Journal]},
    VOLUME = {757},
      YEAR = {2019},
     PAGES = {89--130},
      ISSN = {0075-4102,1435-5345},
   MRCLASS = {58J35 (05C22 35A23 35K08 35R01 53C21)},
  MRNUMBER = {4036571},
MRREVIEWER = {Melanie\ Pivarski},
       DOI = {10.1515/crelle-2017-0038},
       URL = {https://doi.org/10.1515/crelle-2017-0038},
}

@article {CushingLiuPeyerimhoff2020,
    AUTHOR = {Cushing, David and Liu, Shiping and Peyerimhoff, Norbert},
     TITLE = {{Bakry-\'Emery} curvature functions on graphs},
   JOURNAL = {Canad. J. Math.},
  FJOURNAL = {Canadian Journal of Mathematics. Journal Canadien de
              Math\'ematiques},
    VOLUME = {72},
      YEAR = {2020},
    NUMBER = {1},
     PAGES = {89--143},
      ISSN = {0008-414X,1496-4279},
   MRCLASS = {05C50 (53C21)},
  MRNUMBER = {4045968},
MRREVIEWER = {Nelia\ Charalambous},
       DOI = {10.4153/cjm-2018-015-4},
       URL = {https://doi.org/10.4153/cjm-2018-015-4},
}

@misc{Munch2019CD,
      title={{Li--Yau} inequality under {$CD(0,n)$} on graphs}, 
      author={Florentin M\"unch},
      year={2019},
      eprint={1909.10242},
      note={arXiv:1909.10242},
      archivePrefix={arXiv},
      primaryClass={math.DG},
      url={https://arxiv.org/abs/1909.10242}, 
}

@misc{BlacharPajotSalez2026,
      title={Edge-regular graphs with non-negative curvature have polynomial growth}, 
      author={Guy Blachar and Hervé Pajot and Justin Salez},
      year={2026},
      eprint={2606.11094},
      note={arXiv:2606.11094},
      archivePrefix={arXiv},
      primaryClass={math.CO},
      url={https://arxiv.org/abs/2606.11094}, 
}

@article {Wang2010LogHarnack,
    AUTHOR = {Wang, Feng-Yu},
     TITLE = {Harnack inequalities on manifolds with boundary and
              applications},
   JOURNAL = {J. Math. Pures Appl. (9)},
  FJOURNAL = {Journal de Math\'ematiques Pures et Appliqu\'ees. Neuvi\`eme
              S\'erie},
    VOLUME = {94},
      YEAR = {2010},
    NUMBER = {3},
     PAGES = {304--321},
      ISSN = {0021-7824,1776-3371},
   MRCLASS = {58J32 (58J35 58J65 60J60)},
  MRNUMBER = {2679029},
MRREVIEWER = {Ming\ Liao},
       DOI = {10.1016/j.matpur.2010.03.001},
       URL = {https://doi.org/10.1016/j.matpur.2010.03.001},
}

@article {Yau1975Harmonic,
    AUTHOR = {Yau, Shing Tung},
     TITLE = {Harmonic functions on complete {R}iemannian manifolds},
   JOURNAL = {Comm. Pure Appl. Math.},
  FJOURNAL = {Communications on Pure and Applied Mathematics},
    VOLUME = {28},
      YEAR = {1975},
     PAGES = {201--228},
      ISSN = {0010-3640,1097-0312},
   MRCLASS = {53C20 (31C05)},
  MRNUMBER = {431040},
MRREVIEWER = {Yoshiaki\ Maeda},
       DOI = {10.1002/cpa.3160280203},
       URL = {https://doi.org/10.1002/cpa.3160280203},
}

@article {ColdingMinicozzi1997,
    AUTHOR = {Colding, Tobias H. and Minicozzi, II, William P.},
     TITLE = {Harmonic functions on manifolds},
   JOURNAL = {Ann. of Math. (2)},
  FJOURNAL = {Annals of Mathematics. Second Series},
    VOLUME = {146},
      YEAR = {1997},
    NUMBER = {3},
     PAGES = {725--747},
      ISSN = {0003-486X,1939-8980},
   MRCLASS = {53C21 (31B05 58G30)},
  MRNUMBER = {1491451},
MRREVIEWER = {Tanya\ J.\ Christiansen},
       DOI = {10.2307/2952459},
       URL = {https://doi.org/10.2307/2952459},
}

@article {Li1997HarmonicSections,
    AUTHOR = {Li, Peter},
     TITLE = {Harmonic sections of polynomial growth},
   JOURNAL = {Math. Res. Lett.},
  FJOURNAL = {Mathematical Research Letters},
    VOLUME = {4},
      YEAR = {1997},
    NUMBER = {1},
     PAGES = {35--44},
      ISSN = {1073-2780},
   MRCLASS = {53C21 (53C20)},
  MRNUMBER = {1432808},
MRREVIEWER = {Linghai\ Zhang},
       DOI = {10.4310/MRL.1997.v4.n1.a4},
       URL = {https://doi.org/10.4310/MRL.1997.v4.n1.a4},
}

@article {Kleiner2010,
    AUTHOR = {Kleiner, Bruce},
     TITLE = {A new proof of {G}romov's theorem on groups of polynomial
              growth},
   JOURNAL = {J. Amer. Math. Soc.},
  FJOURNAL = {Journal of the American Mathematical Society},
    VOLUME = {23},
      YEAR = {2010},
    NUMBER = {3},
     PAGES = {815--829},
      ISSN = {0894-0347,1088-6834},
   MRCLASS = {20F65 (20F67 20F69)},
  MRNUMBER = {2629989},
MRREVIEWER = {Fran\c cois\ Dahmani},
       DOI = {10.1090/S0894-0347-09-00658-4},
       URL = {https://doi.org/10.1090/S0894-0347-09-00658-4},
}

@article {HuaJostLiu2015,
    AUTHOR = {Hua, Bobo and Jost, J\"urgen and Liu, Shiping},
     TITLE = {Geometric analysis aspects of infinite semiplanar graphs with
              nonnegative curvature},
   JOURNAL = {J. Reine Angew. Math.},
  FJOURNAL = {Journal f\"ur die Reine und Angewandte Mathematik. [Crelle's
              Journal]},
    VOLUME = {700},
      YEAR = {2015},
     PAGES = {1--36},
      ISSN = {0075-4102,1435-5345},
   MRCLASS = {05C10 (53C45)},
  MRNUMBER = {3318509},
       DOI = {10.1515/crelle-2013-0015},
       URL = {https://doi.org/10.1515/crelle-2013-0015},
}

@article {HuaJost2015SemiplanarII,
    AUTHOR = {Hua, Bobo and Jost, J\"urgen},
     TITLE = {Geometric analysis aspects of infinite semiplanar graphs with
              nonnegative curvature {II}},
   JOURNAL = {Trans. Amer. Math. Soc.},
  FJOURNAL = {Transactions of the American Mathematical Society},
    VOLUME = {367},
      YEAR = {2015},
    NUMBER = {4},
     PAGES = {2509--2526},
      ISSN = {0002-9947,1088-6850},
   MRCLASS = {53C21 (05C81 31C20)},
  MRNUMBER = {3301872},
MRREVIEWER = {Norbert\ Peyerimhoff},
       DOI = {10.1090/S0002-9947-2014-06167-9},
       URL = {https://doi.org/10.1090/S0002-9947-2014-06167-9},
}

@incollection {Delmotte1998Harmonic,
    AUTHOR = {Delmotte, Thierry},
     TITLE = {Harnack inequalities on graphs},
 BOOKTITLE = {S\'eminaire de {T}h\'eorie {S}pectrale et {G}\'eom\'etrie,
              {V}ol.\ 16, {A}nn\'ee 1997--1998},
    SERIES = {S\'emin. Th\'eor. Spectr. G\'eom.},
    VOLUME = {16},
     PAGES = {217--228},
 PUBLISHER = {Univ. Grenoble I, Saint-Martin-d'H\`eres},
      YEAR = {[1998]},
   MRCLASS = {58J65},
  MRNUMBER = {1666463},
MRREVIEWER = {Alberto\ G.\ Setti},
       DOI = {10.5802/tsg.202},
       URL = {https://doi.org/10.5802/tsg.202},
}

@article {Hua2019Liouville,
    AUTHOR = {Hua, Bobo},
     TITLE = {Liouville theorem for bounded harmonic functions on manifolds
              and graphs satisfying non-negative curvature dimension
              condition},
   JOURNAL = {Calc. Var. Partial Differential Equations},
  FJOURNAL = {Calculus of Variations and Partial Differential Equations},
    VOLUME = {58},
      YEAR = {2019},
    NUMBER = {2},
     PAGES = {Paper No. 42, 8},
      ISSN = {0944-2669,1432-0835},
   MRCLASS = {31C05 (05C81)},
  MRNUMBER = {3910593},
       DOI = {10.1007/s00526-019-1485-4},
       URL = {https://doi.org/10.1007/s00526-019-1485-4},
}

@article {Hua2021UnboundedCaloric,
    AUTHOR = {Hua, Bobo},
     TITLE = {Ancient caloric functions on graphs with unbounded
              {L}aplacians},
   JOURNAL = {Int. Math. Res. Not. IMRN},
  FJOURNAL = {International Mathematics Research Notices. IMRN},
      YEAR = {2021},
    NUMBER = {4},
     PAGES = {3056--3077},
      ISSN = {1073-7928,1687-0247},
   MRCLASS = {05C22 (58J35)},
  MRNUMBER = {4218346},
MRREVIEWER = {Mingqing\ Zhai},
       DOI = {10.1093/imrn/rnaa045},
       URL = {https://doi.org/10.1093/imrn/rnaa045},
}

@incollection {elworthy,
    AUTHOR = {Elworthy, K. D.},
     TITLE = {Manifolds and graphs with mostly positive curvatures},
 BOOKTITLE = {Stochastic analysis and applications ({L}isbon, 1989)},
    SERIES = {Progr. Probab.},
    VOLUME = {26},
     PAGES = {96--110},
 PUBLISHER = {Birkh\"auser Boston, Boston, MA},
      YEAR = {1991},
      ISBN = {0-8176-3567-X},
   MRCLASS = {58G25 (39A12)},
  MRNUMBER = {1168070},
MRREVIEWER = {J\'ozef\ Dodziuk},
}

@incollection {michael,
    AUTHOR = {Schmuckenschl\"ager, Michael},
     TITLE = {Curvature of nonlocal {M}arkov generators},
 BOOKTITLE = {Convex geometric analysis ({B}erkeley, {CA}, 1996)},
    SERIES = {Math. Sci. Res. Inst. Publ.},
    VOLUME = {34},
     PAGES = {189--197},
 PUBLISHER = {Cambridge Univ. Press, Cambridge},
      YEAR = {1999},
      ISBN = {0-521-64259-0},
   MRCLASS = {60J75 (52C99 60J35)},
  MRNUMBER = {1665591},
MRREVIEWER = {Joerg-Uwe\ Loebus},
}

@article {FathiShu2018,
    AUTHOR = {Fathi, Max and Shu, Yan},
     TITLE = {Curvature and transport inequalities for {M}arkov chains in
              discrete spaces},
   JOURNAL = {Bernoulli},
  FJOURNAL = {Bernoulli. Official Journal of the Bernoulli Society for
              Mathematical Statistics and Probability},
    VOLUME = {24},
      YEAR = {2018},
    NUMBER = {1},
     PAGES = {672--698},
      ISSN = {1350-7265,1573-9759},
   MRCLASS = {60J27 (53C21 60E15 60J45)},
  MRNUMBER = {3706773},
MRREVIEWER = {Feng-Yu\ Wang},
       DOI = {10.3150/16-BEJ892},
       URL = {https://doi.org/10.3150/16-BEJ892},
}

@article {HornJGT,
    AUTHOR = {Horn, Paul and Purcilly, Adam and Stevens, Alex},
     TITLE = {Graph curvature and local discrepancy},
   JOURNAL = {J. Graph Theory},
  FJOURNAL = {Journal of Graph Theory},
    VOLUME = {108},
      YEAR = {2025},
    NUMBER = {2},
     PAGES = {337--360},
      ISSN = {0364-9024,1097-0118},
   MRCLASS = {05C50 (05C63 53C23)},
  MRNUMBER = {4839884},
MRREVIEWER = {Vladimir\ Mityushev},
       DOI = {10.1002/jgt.23176},
       URL = {https://doi.org/10.1002/jgt.23176},
}

@article {HM2024,
    AUTHOR = {Hua, Bobo and M\"unch, Florentin},
     TITLE = {Graphs with nonnegative curvature outside a finite subset,
              harmonic functions, and number of ends},
   JOURNAL = {J. Lond. Math. Soc. (2)},
  FJOURNAL = {Journal of the London Mathematical Society. Second Series},
    VOLUME = {110},
      YEAR = {2024},
    NUMBER = {6},
     PAGES = {Paper No. e70034, 19},
      ISSN = {0024-6107,1469-7750},
   MRCLASS = {53C21 (31C20 53A70 53C23)},
  MRNUMBER = {4829680},
MRREVIEWER = {Taiki\ Yamada},
       DOI = {10.1112/jlms.70034},
       URL = {https://doi.org/10.1112/jlms.70034},
}

@article {SalezJEMS,
    AUTHOR = {Salez, Justin},
     TITLE = {Cutoff for non-negatively curved {M}arkov chains},
   JOURNAL = {J. Eur. Math. Soc. (JEMS)},
  FJOURNAL = {Journal of the European Mathematical Society (JEMS)},
    VOLUME = {26},
      YEAR = {2024},
    NUMBER = {11},
     PAGES = {4375--4392},
      ISSN = {1435-9855,1435-9863},
   MRCLASS = {60J27},
  MRNUMBER = {4780485},
MRREVIEWER = {Leonid\ V.\ Bogachev},
       DOI = {10.4171/jems/1348},
       URL = {https://doi.org/10.4171/jems/1348},
}

@article {KMY2021,
    AUTHOR = {Kempton, Mark and M\"unch, Florentin and Yau, Shing-Tung},
     TITLE = {A homology vanishing theorem for graphs with positive
              curvature},
   JOURNAL = {Comm. Anal. Geom.},
  FJOURNAL = {Communications in Analysis and Geometry},
    VOLUME = {29},
      YEAR = {2021},
    NUMBER = {6},
     PAGES = {1449--1473},
      ISSN = {1019-8385,1944-9992},
   MRCLASS = {53C21 (57M15)},
  MRNUMBER = {4367431},
MRREVIEWER = {Marko\ \v Zivkovi\'c},
       DOI = {10.4310/CAG.2021.v29.n6.a5},
       URL = {https://doi.org/10.4310/CAG.2021.v29.n6.a5},
}

@article {MuenchRose2020,
    AUTHOR = {M\"unch, Florentin and Rose, Christian},
     TITLE = {Spectrally positive {Bakry-\'Emery} {R}icci curvature on
              graphs},
   JOURNAL = {J. Math. Pures Appl. (9)},
  FJOURNAL = {Journal de Math\'ematiques Pures et Appliqu\'ees. Neuvi\`eme
              S\'erie},
    VOLUME = {143},
      YEAR = {2020},
     PAGES = {334--344},
      ISSN = {0021-7824,1776-3371},
   MRCLASS = {53C21 (05C99 35K10 35R02)},
  MRNUMBER = {4163132},
MRREVIEWER = {Thierry\ Coulhon},
       DOI = {10.1016/j.matpur.2020.03.008},
       URL = {https://doi.org/10.1016/j.matpur.2020.03.008},
}

@article {ColdingMinicozzi1998,
    AUTHOR = {Colding, Tobias H. and Minicozzi, II, William P.},
     TITLE = {Liouville theorems for harmonic sections and applications},
   JOURNAL = {Comm. Pure Appl. Math.},
  FJOURNAL = {Communications on Pure and Applied Mathematics},
    VOLUME = {51},
      YEAR = {1998},
    NUMBER = {2},
     PAGES = {113--138},
      ISSN = {0010-3640,1097-0312},
   MRCLASS = {53C21 (53C20 58E15)},
  MRNUMBER = {1488297},
MRREVIEWER = {Man\ Chun\ Leung},
       DOI = {10.1002/(SICI)1097-0312(199802)51:2<113::AID-CPA1>3.0.CO;2-E},
       URL =
              {https://doi.org/10.1002/(SICI)1097-0312(199802)51:2<113::AID-CPA1>3.0.CO;2-E},
}

@article {Salez2022Expanders,
    AUTHOR = {Salez, Justin},
     TITLE = {Sparse expanders have negative curvature},
   JOURNAL = {Geom. Funct. Anal.},
  FJOURNAL = {Geometric and Functional Analysis},
    VOLUME = {32},
      YEAR = {2022},
    NUMBER = {6},
     PAGES = {1486--1513},
      ISSN = {1016-443X,1420-8970},
   MRCLASS = {53C23},
  MRNUMBER = {4536468},
       DOI = {10.1007/s00039-022-00618-3},
       URL = {https://doi.org/10.1007/s00039-022-00618-3},
}

@article {MunchSalez2023Expansion,
    AUTHOR = {M\"unch, Florentin and Salez, Justin},
     TITLE = {Mixing time and expansion of non-negatively curved {M}arkov
              chains},
   JOURNAL = {J. \'Ec. polytech. Math.},
  FJOURNAL = {Journal de l'\'Ecole polytechnique. Math\'ematiques},
    VOLUME = {10},
      YEAR = {2023},
     PAGES = {575--590},
      ISSN = {2429-7100,2270-518X},
   MRCLASS = {60J10},
  MRNUMBER = {4567745},
       DOI = {10.5802/jep.226},
       URL = {https://doi.org/10.5802/jep.226},
}
}{%
}
\end{document}